\documentclass[11pt,a4paper]{article}

\usepackage{amsmath}    % need for subequations
\usepackage{amssymb}
\usepackage{graphicx}   % need for figures
\usepackage{verbatim}   % useful for program listings
\usepackage{color}      % use if color is used in text
\usepackage{enumerate}
\usepackage{listings}
\usepackage{xcolor}
\usepackage{slashed}
\usepackage{adjustbox}
\usepackage{tikz}
\usepackage{amsthm}    % need for subequations
\usepackage{amssymb}
\usepackage{mathrsfs}
\usepackage{bbold}
\usepackage[british]{babel}
\usepackage[utf8x]{inputenc}
\usepackage{array}
\usepackage{subcaption}
\usepackage[dvipsnames]{xcolor}
\usepackage{tabularx}
\usepackage{graphicx,wrapfig,lipsum}
\usepackage{ifplatform}
\usepackage{preview}
\usepackage{environ}
\usepackage{graphics}
\usepackage{epstopdf}
\usepackage{amsmath}
\usepackage{amssymb}
\usepackage{amsfonts}	
\usepackage{moreverb}
\usepackage{dsfont}
\usepackage{grffile}
\usepackage{bm}
\usepackage{multirow}
\usepackage{frcursive}
\usepackage{mathrsfs}
\usepackage{nicefrac}
\usepackage{url}
\usepackage[toc]{appendix}
\usepackage{booktabs}

\usepackage[textsize=footnotesize backgroundcolor=yellow!70, bordercolor=orange]{todonotes}

\newcommand\BibTeX{{\rmfamily B\kern-.05em \textsc{i\kern-.025em b}\kern-.08em
T\kern-.1667em\lower.7ex\hbox{E}\kern-.125emX}}

\usepackage{lmodern}  % for bold teletype font
\usepackage{xcolor}   % for \textcolor
\usepackage{listings}
\usepackage{hyperref} 
\hypersetup{
    colorlinks=true,                          
    linkcolor=blue, % Couleur des liens internes
    citecolor=red, % Couleur des numéros de la biblio dans le corps
    urlcolor=blue  } % Couleur des url
\usepackage[square,numbers]{natbib}
\RequirePackage[hyperpageref]{backref}
\backreffrench
\renewcommand*{\backref}[1]{}  % Disable standard
\renewcommand*{\backrefalt}[4]{% Detailed backref
  \ifcase #1 %
  \relax%(Not cited.)%
  \or
  (Cited page~#2.)%
  \else
  (Cited page~#2.)%
  \fi}

\renewcommand{\epsilon}{\varepsilon }

\newfont{\numerikEleven}{ecrm1000}
\newfont{\numerikTen}{cmss10}
\newfont{\numerikNine}{cmss9}
\newfont{\numerikEight}{cmss8}
\newfont{\numerikSeven}{cmss7}
\newfont{\numerikSix}{cmss6}
\everymath{\displaystyle}

\newcommand{\p}{\partial}

\newtheorem{theorem}{Theorem}

\newtheorem{proposition}[theorem]{Proposition}

\newtheorem{remark}[theorem]{Remark}

\definecolor{teal}{rgb}{0.0, 0.5, 0.5}

\usepackage{hyperref}
\numberwithin{equation}{section}

\newcommand{\dt}{\Delta t}
\newcommand{\dx}{\Delta x}

\newcommand{\ha}{\frac{1}{2}}
\renewcommand{\epsilon}{\varepsilon}
\DeclareMathOperator{\sech}{sech}

\usepackage[final]{showlabels}

\title{Semi-implicit strategies for the Serre-Green-Naghdi equations in hyperbolic form. Is hyperbolic relaxation really a good idea?} 
\author{Emanuele Macca$^{\times*}$, Walter Boscheri$^{\dagger}$ and Mario Ricchiuto$^{+}$\\[20pt]
$^{\times}$ Università degli Studi di Catania, Catania, Italy 
$^*$ Corresponding author \\
$^\dagger$ Laboratoire de Math\'{e}matiques UMR 5127 CNRS, Universit\'{e} Savoie Mont Blanc, France\\
$^{+}$ Centre Inria de l’Universit\'{e} de Bordeaux, CNRS UMR 5251, Talence, France}
\date{}
\begin{document}

\maketitle

 \begin{abstract}
The Serre-Green-Naghdi (SGN) equations provide a valuable framework for modelling fully nonlinear and weakly dispersive shallow-water flows. However, their elliptic formulation can considerably increase the computational cost compared to the Saint-Venant equations. To overcome this difficulty, hyperbolic models (hSGN) have been proposed that replace the elliptic operators with first-order hyperbolic formulations augmented by relaxation terms, which recover the original elliptic formulation in the stiff limit. Yet, as the relaxation parameter~$\lambda$ increases, explicit schemes face restrictive stability constraints that may offset these advantages.  

To mitigate this limitation, we introduce a semi-implicit (SI) integration strategy for the hSGN system, where the stiff acoustic terms are treated implicitly within an IMEX Runge--Kutta framework, while the advective components remain explicit. The proposed approach mitigates the CFL stability restriction and maintains dispersive accuracy at a moderate computational cost. Numerical results confirm that the combination of hyperbolization and semi-implicit time integration provides an efficient and accurate alternative to both classical SGN and fully explicit hSGN solvers.
\end{abstract}
% \keywords{}


% \tableofcontents

\section{Introduction}

%\begin{itemize}
%\item dispersive wave equations relevance in applications
%\item dispersive wave models and time-stepping
%\item larger overlook on parabolic equations with higher derivatives
%\item hyperbolization : a short review
%\item focus here on models involving relazation constant which affect the eqigenvalues 
%\item imex techniques: short overview and   contrast with original aim of hyperbolization 
%\item  main question: is hyperbolization a good idea when accuracy requires small time steps ? 
%\item scope of main contributios
%\begin{enumerate}
%\item novel and robust splitting allowing SI time integration of HSGN based on  augmented Lagrangian method
%\item combined with second order  SI time stepping using classical RK-ARS(2,2,2) method 
%\item invesitgation of cost/accuracy trade off for  imex   dispersive shallow water equations  
%\item some recommendations are provided 
%\end{enumerate}
%\item structure of the paper
%\end{itemize}

Dispersive shallow-water models are fundamental tools in the mathematical and numerical study of free-surface flows, with applications that range from coastal engineering and tsunami propagation to laboratory experiments and reduced-order modelling \cite{Lannes_2020,MULLER199389}. Among these, the Serre–Green–Naghdi (SGN) equations stand out as a reliable framework for describing fully nonlinear and weakly dispersive wave motion in shallow to intermediate water depths. Unlike the classical shallow-water (Saint–Venant) equations, the SGN system incorporates dispersive corrections, enabling a more accurate representation of wave dynamics in regimes where dispersion plays a significant role \cite{Klein2021,Duchene,BONNETON2011589}.

From a numerical perspective, 
two alternative approaches can be found in litteraure. One is based on
the classical fomrulation of the equations
which embedds mixed space and time derivatives \cite{Lannes_2020}. More precisely such formulation embeds 
an ellipic nonlienar operator 
whose inversion is necessary to recover the
depth averaged velocity \cite{wei1995fully,bonneton2011splitting,Kazolea,GS24}.
This  leads  to
%
%uses some formulation of 
%the  is 
%
%current approaches to solve the SGN equations
%lead to %exhibit 
%two contrasting potential bottlenecks. % %inherent difficulties. 
%When using  the classical formulation %leads to a mixed elliptic–hyperbolic character, since 
%the velocity update typically involves the solution of an elliptic problem with time dependent coefficients. This %requirement 
%gives 
a dramatic  increase in 
computational cost, especially compared to the underlying shallow water equations for which efficient %  and hinders the construction of 
fully explicit conservative methods can be used. 
%In contrast, alternative reformulations often give rise to fast acoustic-like modes associated with the dispersive correction terms. When large relaxation parameters are introduced, these modes impose severe constraints on the maximum admissible time step for explicit schemes. These considerations motivate the development of alternative formulations and integration techniques that preserve the accuracy of the SGN model while mitigating its computational complexity.

One of the most common strategies to circumvent the elliptic step in the SGN equations is to resort to relaxation systems, often referred to as “hyperbolization”. 
In this case, rather than solving directly the dispersive equations one solves
an augmented system of  first order PDEs,
which includes several 
 auxiliary variables and  appropriately designed relaxation sources which allow to 
 approximate the dispersive operators
 when  a relaxation parameter $\epsilon=1/\lambda$
 tends to zero. Several such approaches have been proposed over the years for dispersive equations
\cite{Alireza,Favrie,escalante_2020,Busto,Sergey,rr25}, as well as for other models involving high order derivatives \cite{firas24,firas25,Busto-a,Firas_2019}.
The first order part of the system is usually designed to be hiperbolic, so we can speak of a hyperbolic Serre-Green-Naghdi model (by some abuse of language).
 This hSGN system, is thus a hyperbolic set of balance laws, with stiff relaxation terms. In the asymptotic relaxation limit, the original parabolic equations are retrieved. In this way, the difficulties associated with nonlocal elliptic operators are removed, and the resulting system can be tackled with standard high-order methods for hyperbolic balanced laws. This advantage, however, comes with a significant trade-off: the augmented wave speeds typically grow with $\lambda$, so that for large relaxation parameters they can far exceed the natural shallow-water celerity $\sqrt{gh}$, leading to a much more restrictive CFL condition in explicit schemes.

The core question that this work addresses is therefore practical and quantitative: when accuracy demands push the model towards large $\lambda$, does hyperbolization still pay off? In other words, although hyperbolization removes elliptic inversions, it also increases system size and possibly the stiffness of acoustic modes, so the net computational gain is not automatically guaranteed. To answer this question we pursue a complementary strategy: instead of treating the entire hyperbolic system explicitly, we design semi-implicit (SI) time-stepping schemes that treat the stiff ($\lambda$-dependent) acoustic coupling implicitly while handling advective and non-stiff terms explicitly. This splitting aims to combine the structural advantages of the hSGN model with a time discretization that relaxes the severe hyperbolic CFL constraint and leads, in practice, to an implicit solver that is affordable and subject to a material-type CFL rather than to the overly restrictive acoustic CFL.

In this paper we present a new splitting of the hSGN equations tailored to SI integration, and we couple this splitting with a second-order IMEX Runge--Kutta time stepping technique. The combination yields a method that damps acoustic modes through the implicit sub-step and produces a linear or mildly nonlinear elliptic-type problem for the depth update that can be solved efficiently. We discuss stability and consistency properties of the splitting and we provide a careful assessment of the cost/accuracy trade-off between (i) the classical elliptic--hyperbolic SGN solver, (ii) fully explicit discretizations of the hyperbolic relaxation, and (iii) the proposed SI-hSGN approach.

Numerical experiments illustrate the practical performance of the method and delimit the regimes where SI-hSGN is advantageous. Benchmark tests include solitary-wave propagation, dispersive wave breaking scenarios and other standard tests that highlight dispersive accuracy, nonlinear interactions and long-time behaviour. Based on these results we provide pragmatic recommendations for choosing the stiffness parameter $\lambda$ and tuning the IMEX time integrator in applications of interest.

The remainder of the paper is organized as follows. We begin by recalling the standard SGN formulation and the hyperbolic relaxation variants used in our study. Then we introduce the semi-implicit splitting and detail the IMEX time discretization, together with stability considerations and an asymptotic analysis. Spatial discretization options for both explicit and SI-hSGN solvers are discussed next. Afterwards we present numerical results that assess accuracy, stability and computational cost, and we conclude with final remarks and practical guidelines for users interested in adopting the proposed approach.

\section{Serre-Green-Naghdi equations in hyperbolic and standard form}

In this section we discuss the different forms of the Serre-Green-Naghdi  (SGN) system used in this paper.
The interested reader is referred to \cite{Favrie,Lannes_2020,escalante_2020,Busto,kazolea2024,rr25}  
and references therein for more details on the SGN equations, and more generally  
on dispersive  wave models.

\subsection{Standard SGN equations}
Let us consider a one-dimensional domain $\Omega$ defined by the position $x \in \mathbb{R}$ and let $t \in \mathbb{R}_+$ be the time coordinate. The propagation of non-linear dispersive waves on flat bathymetry  can be modeled by
the Serre-Green-Naghdi equations which write
\begin{equation} \label{eq:SGN_flat}
\begin{split}
h_t + & (hu)_x =0,\\
(hu)_t + & (hu^2 + gh^2/2 + p_{\textsf{nh}})_x =0,\\
\end{split}
\end{equation}
where  $h$ denotes the water depth, $u$ the dept averaged velocity, $g$ the gravity acceleration,
and with $p_{\textsf{nh}}$ the non-hydrostatic pressure
\begin{equation}\label{eq:SGN_p}
p_{\textsf{nh}}= - \dfrac{h^3}{3}\big(\dfrac{Du}{Dt}\big)_x +  2\dfrac{h^3}{3} (u_x)^2,
\end{equation}
with $D/Dt$ denoting the material or Lagrangian derivative:
\begin{equation}
    \dfrac{D}{Dt}(\cdot) = (\cdot)_t + u (\cdot)_x\,.
\end{equation} 
The above system admits several equivalent formulations, more or less suitable for numerical approximation. 
To update the velocity, all of these formulations require inverting some equivalent form of the operator
\begin{equation}
    (hu)_t -(\dfrac{h^3}{3} u_{tx})_x = \text{RHS}
\end{equation}
which is embedded in    \eqref{eq:SGN_flat}-\eqref{eq:SGN_p}.
In this work, we investigate  hyperbolic reformulations which allow in principle to avoid this inversion. These are discussed in the next sections.

\subsection{Hyperbolic SGN (hSGN) equations without bathymetry}
The hyperbolized SGN model (hSGN)  obtained using an extended Lagrangian
formulation can be written in the flat bathymetry case as (see \cite{Favrie,Busto,TKACHENKO2023111901})
\begin{equation}
        \label{hSGN_eq}
    \begin{cases}
        & h_t + (hu)_x = 0,\\
        & (hu)_t + \left(  hu^2 + p  \right)_x = 0,\\
        & (h\eta)_t + (hu\eta)_x = hw,\\
        & (hw)_t + (huw)_x = - \lambda \left( \frac{\eta }{h} - 1 \right),
    \end{cases}
\end{equation}
with
\begin{equation}\label{eq:hSGN_p}
p:= g\dfrac{h^2}{2} + h \, \dfrac{\lambda}{3}\dfrac{\eta}{h} \left(1-\dfrac{\eta}{h}\right).
\end{equation}
In addition to the physical variables already defined, $w$ and $\eta$ are here  relaxation variables and $\lambda\gg 1$ a relaxation parameter. The analysis of \cite{duchene2019rigorous} 
shows that for large values of $\lambda$ then the above system provides an order $\mathcal{O}(1/\lambda)$
approximation of the SGN model, with $\eta=h+\mathcal{O}(1/\lambda)$.\\

The above system can be shown to be fully hyperbolic, with  eigenvalues
\begin{equation}\label{eq:hSGN_eigenvalues}
\lambda_1 = u - c,\; \lambda_2= \lambda_3 = u,\; 
\lambda_4= u+ c\;,
\end{equation}
and linearly independent eigenvectors.  This allows in principle to use standard high order \emph{explicit} numerical methods avoiding the inversion of second order operators.\\

Note that the extended celerity $c$ in the eigenvalues is defined 
by 
\begin{equation}\label{eq:hSGN_c}
c^2 = c_h^2 + \dfrac{\lambda}{3}\dfrac{\eta^2}{h^2}\;,\quad
c_h^2 = gh,
\end{equation}
where $c_h$ is the usual shallow water celerity. Clearly, since $\eta/h=1 +\mathcal{O}(1/\lambda)$, one has  $c\gg c_h$ in general, and especially 
for large $\lambda$. \\

For later use, it is quite interesting to note that the pressure derivative can be written as 
\begin{equation}\label{eq:hSGN_px}
    \begin{aligned}
        p_x & = (c^2 - \lambda \eta \alpha) h_x + \lambda \alpha (h\eta)_x \\
        & = a^2\, h_x + \lambda\alpha (h\eta)_x 
    \end{aligned}
\end{equation}
where 
\begin{equation}\label{eq:hSGN_alpha}
a^2 := c^2 - \lambda \eta \alpha\;,\quad 
\alpha := - \frac{1}{3h} \left(2 \frac{\eta}{h} - 1 \right).
\end{equation}
Because for $\lambda \gg 1 $ we have $\eta/h=1 +\mathcal{O}(1/\lambda)$, it is reasonable to assume $\alpha <0$.

\subsection{Including bathymetric effects}

The effects of bathymetry in the hyperbolic system can be accounted for via two simple modifications. A   hyperbolized SGN model with batrhymetry can be written according to \cite{Busto,rr25} as 
\begin{equation}
        \label{hSGNb_eq}
    \begin{cases}
        & h_t + (hu)_x = 0,\\
        & (hu)_t + \left(  hu^2 + \frac{1}{2} gh^2 +   \eta\tilde p   \right)_x +(gh + \dfrac{3}{2}\tilde p)b_x = 0,\\
        & (h\tilde \eta)_t + (hu\tilde\eta)_x = hw,\\
        & (hw)_t + (huw)_x = 3 \tilde p ,\\
        &\tilde p =  \dfrac{\lambda}{3} (1-\frac{\eta }{h}),\\
        &\tilde \eta = \eta + \dfrac{3}{2} b,
    \end{cases}
\end{equation}
where $b=b(x)$ is the topography. Note that, as remarked in \cite{rr25}, the asymptotic limit of the above system does not coincide with the usual SGN model 
but with  a slightly different approximation which can be obtained under classical mild-slope hypotheses.  In particular, following the reasoning of \cite{rr25}, one can show that the limit of system \eqref{hSGNb_eq} writes
\begin{equation}
        \label{SGNb_eq}
       \begin{cases}
        & h_t + (hu)_x = 0,\\
        & (hu)_t + \left(  hu^2 + \frac{1}{2} gh^2 +   h \tilde p_{\mathsf{nh}}     \right)_x +(gh + \dfrac{3}{2}\tilde p_{\mathsf{nh}} )b_x = 0,
        \end{cases}
\end{equation}
where now
\begin{equation}
        \label{SGNb_p}
       \tilde p_{\mathsf{nh}} = -\dfrac{h^2}{3} \big( \dfrac{D u}{Dt} \big)_x + \dfrac{hb_x}{2}  \dfrac{D u}{Dt} + \dfrac{2}{3}h^2 u_x^2 + h^2u^2 b_{xx}.
        \end{equation}
Comparison with the expressions of the original SGN model, reported e.g. in \cite{Filippini,Lannes_2020,GS24}, exhibits some differences.
Numerical tests   on classical benchmarks show that the  original SGN and the 
above mild-slope approximation  provide identical results even for rapidly varying bathymetries \cite{rr25}.  We will use thus \eqref{hSGNb_eq} and \eqref{SGNb_eq}   with no further discussion on the topic, and we refer to this system with the abbreviation hSGN. The interested reader can refer to  \cite{Busto,rr25}  for more details.

\section{Semi-implicit time integration   for the hSGN model}

% WB: too friendly for my taste... but up to you :-)
%As stated in the Introduction, the main motivation of this work  is essentially spelled in  the question contained in the title: is hyperbolic relaxation really a good idea? 

\begin{figure}[!ht]
\includegraphics[width=0.5\textwidth]{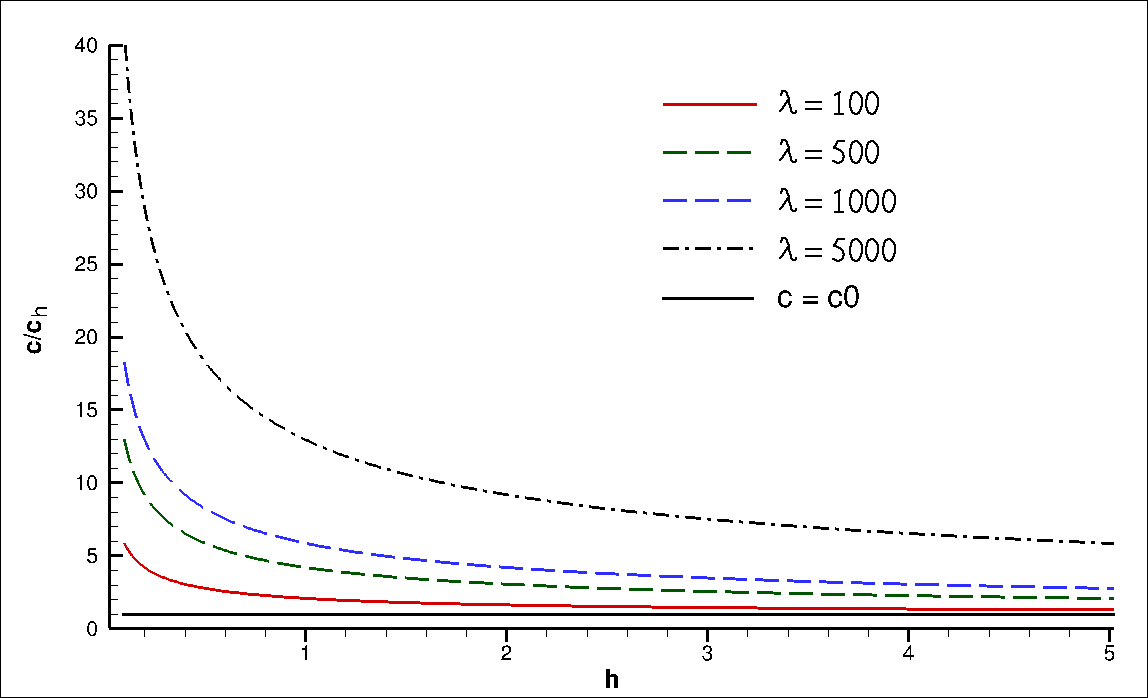}\includegraphics[width=0.5\textwidth]{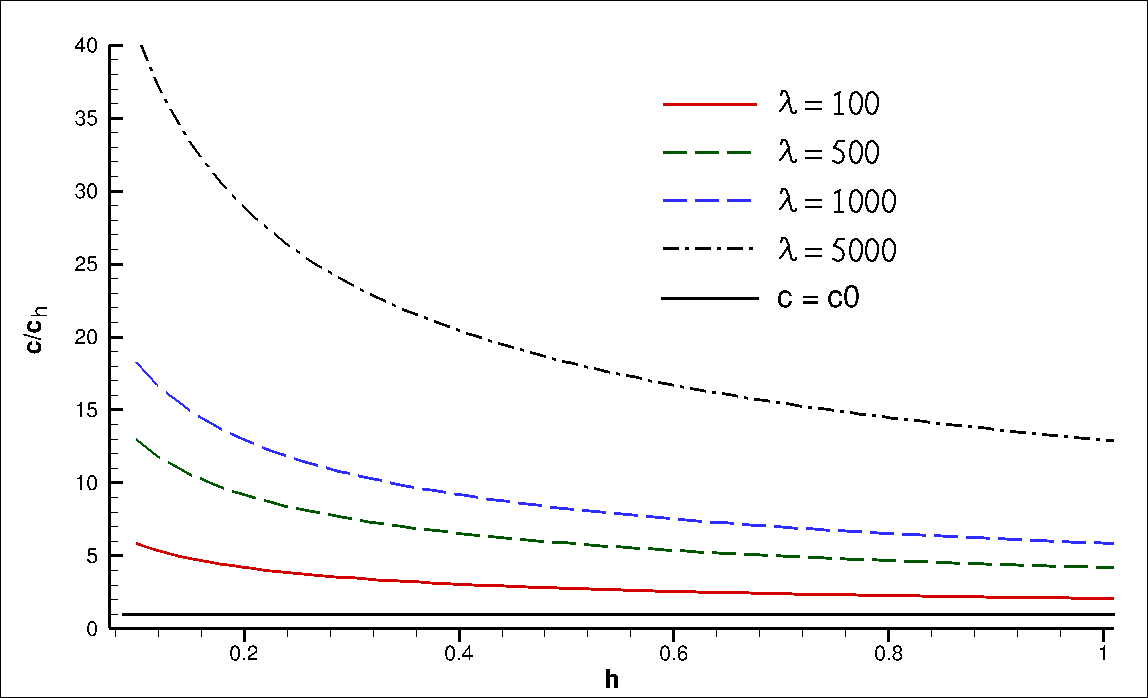}
    \caption{Celerity ratio $c/c_h$ as a function of the water depth for different values of $\lambda$.\label{fig:c_c0}}
\end{figure}

As previously stated in the Introduction, we aim at investigating whether hyperbolic relaxation is actually an effective modeling approach. Indeed, the main rationale behind the development  of hyperbolic approximations
is to avoid the inversion of the elliptic operators appearing in \eqref{eq:SGN_flat}-\eqref{eq:SGN_p} and
\eqref{SGNb_eq}-\eqref{SGNb_p}.   Indeed, this is true. However, there are two  points which may attenuate critically the advantage of using a hyperbolic formulation.  One is the number of equations which is substantially larger (double in one dimension, and larger of a factor 5/3 in two dimensions). The second is the time step limitation. Several previous works show that dispersive terms  do not necessarily affect the stability of the schemes. This is for example
the case when using centered approximations, and certain residual based approaches   (see e.g. \cite{rf14,Filippini,Kazolea} for spectral  stability results).
In this case, even in the presence of dispersive terms, the relevant time step limitation is the usual shallow water one, associated to the speed $|u| +c_h$. Conversely, the hyperbolic reformulation will require time-steps  inversely proportional to $|u| +c$. As we saw, we expect $c$ to be much larger than $c_h$, and this is indeed the case, especially in intermediate to shallow depths where weakly dispersive corrections become important, as shown in Figure~\ref{fig:c_c0}. We can see that, for large values of $\lambda$, this ratio may exceed one order of magnitude, and even beyond that. This means that, if the overhead in the inversion
of the elliptic operator  is less then a factor 10 or 20,  the hyperbolic reformulation may actually  be slower, as reported in some works e.g. \cite{rr25,GR25}.

A natural question is whether, by some efficient semi-implicit (SI) time discretization, we may offset this issue in favor of the hyperbolic formulation. Indeed, designing and implementing SI and IMplicit-EXplicit (IMEX) 
schemes can be done in a very efficient way    for first-order hyperbolic systems (see e.g. \cite{BR24,BPR-SIAMbook-24} and references therein). 
This is what we propose to study here.  For the ODE 
\begin{equation}\label{eq:ODE}
\dfrac{dQ}{dt} = R(Q,x),  
\end{equation}
the simplest form of SI schemes involves some splitting of the   right hand side  evolution operator $R$ in  a stiff  part $R_I$, treated implicitly, and a non-stiff part $R_E$, treated explicitly: 
\begin{equation}\label{eq:basicSI}
\dfrac{Q^{n+1} - Q^n}{\Delta t}  = R_I(Q^{n+1},x) +   R_E(Q^{n},x) \;,\quad
R_I(Q,x) + R_E(Q,x) = R(Q,x).
\end{equation}
Our first task is to derive a suitable splitting for the hSGN systems \eqref{hSGN_eq} and  \eqref{hSGNb_eq}.

\subsection{Splitting for the hSGN equations}~\label{sec:splitting}

Let us start from the model \eqref{hSGN_eq}.  The acoustic part of the eigenstructure is associated to the coupling  between the pressure,  mass conservation, and the wave equation defined by the last two equations. Following the general guideline \eqref{eq:basicSI}, we thus split the system as follows:
\begin{equation}
        \label{hSGN_eq-split}
        \dfrac{\partial}{\partial t} 
        \left(\begin{array}{c} h\\hu \\ h\eta \\ hw
        \end{array}\right)= R_E^{\text{hSGN}} + R_I^{\text{hSGN}}\;,\
\end{equation}
where
\begin{equation}
        \label{hSGN_eq-splitb}
        R_E^{\text{hSGN}}:= -\left(\begin{array}{c} 0\\hu^2 \\ hu \eta  \\ hu w
        \end{array}\right)_x \;,\;\;
        R_I^{\text{hSGN}}:= - \left(\begin{array}{c} hu\\p \\ 0  \\  0
        \end{array}\right)_x -  \left(\begin{array}{c} 0\\ 0 \\  -hw \\   \lambda \left( \frac{\eta }{h} - 1 \right)
        \end{array}\right).  
\end{equation}
The above splitting can be further developed to show how its solution implies the inversion
of a second order operator, allowing to dump all acoustic modes and increasing the stability range. 
To begin with, we introduce the intermediate update 
\begin{equation}
        \label{hSGN_eq-ustar}
        \begin{aligned}
        h^* = & h^n, \\
        (hu)^* = &(hu)^n -\Delta t (hu^2)_x^n,\\
        (h\eta)^* =&(h\eta)^n - \Delta t (hu\eta)_x^n, \\
        (hw)^* = &(hw)^n -\Delta t (huw)_x^n,
        \end{aligned}
\end{equation}
that essentially solves the convective contribution. This allows us to write   the implicit sub-system as
\begin{equation}
        \label{hSGN_eq-split-I}
        \begin{aligned}
            h^{n+1} =  & h^* -\Delta t  (hu)^{n+1}_x,\\
            (hu)^{n+1} = & (hu)^{*} -\Delta t p^{n+1}_x, \\
            (h\eta)^{n+1} = &(h\eta)^* + \Delta t (hw)^{n+1},\\
            (hw)^{n+1} = & (hw)^* -\lambda \Delta t \left(\dfrac{\eta}{h}-1\right)^{n+1}.
        \end{aligned}
\end{equation}
Combining the last two equations above, we can easily find that
\begin{equation}\label{eq:b1}
(h\eta)^{n+1} = \widetilde{(h\eta)} \,b_1(h^{n+1})\;,\quad b_1(h):= \dfrac{h^2}{h^2+\lambda\Delta t^2},
\end{equation}
having set
\begin{equation}
    \widetilde{(h\eta)} := (h\eta)^* + \Delta t (hw)^* + \lambda\Delta t^2.
\end{equation}
We now move on to the first equation in \eqref{hSGN_eq-split-I}, which we can combine with the second one to obtain
\begin{equation}
h^{n+1} = h^* -\Delta t q_x^* +\Delta t^2 p^{n+1}_{xx}.    
\end{equation}
We now use \eqref{eq:hSGN_px} to  evaluate the last term, hence obtaining
\begin{equation}
h^{n+1} -\Delta t^2 (a^2 h_x^{n+1})_x = h^* -\Delta t q_x^*  + \lambda \Delta t^2(\alpha (h\eta)_x^{n+1})_x.
\end{equation}
To close the problem,   we have to use \eqref{eq:b1}. In particular, we can write
\begin{equation}
    (h\eta)^{n+1}_x = b_1\, \widetilde{(h\eta)}_x +  \widetilde{(h\eta)} \partial_h b_1 h_x^{n+1}.
\end{equation}
Simple computations show that
\begin{equation}
    \partial_h b_1= 2\lambda\Delta t^2\dfrac{ h}{(h^2+\lambda\Delta t^2)^2} =  2\dfrac{\lambda\Delta t^2 }{h^3} b_1^2 \ge 0,
\end{equation}
so we can write
\begin{equation}
    (\alpha (h\eta)_x^{n+1})_x = \big[
\alpha b_1\, \widetilde{(h\eta)}_x +  2 \alpha \widetilde{(h\eta)}  \dfrac{\lambda\Delta t^2 }{h^3} b_1^2  h_x^{n+1}\big]_x.
\end{equation}
We now propose to use the following linearized second order problem to evaluate the water depth at the next time step:
\begin{equation}\label{eq:hnew}
\begin{aligned}
h^{n+1} - \Delta t^2 (\kappa(h^n,u^n,\eta^n,w^n)\, h_x^{n+1})_x  =&\;
 h^* -\Delta t q_x^* +\lambda \Delta t^2 (\alpha b_1\, \widetilde{(h\eta)}_x)_x\;,\\[5pt]
 \kappa(h^n,u^n,\eta^n,w^n):= &\; a^2(h^n,\eta^n) +   2 \alpha(h^n,\eta^n) \widetilde{(h\eta)}  \dfrac{\lambda^2\Delta t^2 }{(h^n)^3} b_1^2(h^n) 
 \end{aligned}
\end{equation}
The nature of the last equation is not clear due to the undecided sign of the second term in $\kappa$.
A more detailed study allows us to prove what follows.

\begin{proposition}[Asymptotic coercivity]  Provided $\alpha^n\le 0$,  $h^n> 0$, and $\eta^n >0$,  and provided that the timestep $\Delta t$ is independent of $\lambda$, the second order equation \eqref{eq:hnew} is asymptotically coercive
in the sense that
\begin{equation}
    \Delta t^2 \kappa =
 \Delta t^2 gh^n +\Delta t^2 \lambda ( \Lambda_0 + \mathcal{O}(1/\lambda))\;,\quad
 \Lambda_0 = \dfrac{1}{3}\dfrac{(\eta^n)^2}{(h^n)^2} - \alpha^n h^n > 0.
\end{equation}
\begin{proof}
Using the definition of $a^2$  (eq. \eqref{eq:hSGN_alpha})  we can easily write  (we omit the superscript $^n$ to lighten the notation)
\begin{equation}
    \Delta t^2\kappa= \Delta t^2 gh +  \Delta t^2\dfrac{\lambda }{3} \dfrac{\eta^2}{h^ 2}-  \alpha\lambda \eta \Gamma,
\end{equation}
where
\begin{equation}
    \Gamma:= \Delta t^2 - \dfrac{ \widetilde{(h\eta)} }{h\eta}
      \dfrac{2h^2}{\lambda}   \dfrac{\lambda^2 \Delta t^4}{(h^2+ \lambda\Delta t^2)^2}.
\end{equation}
A Taylor   development  in terms of $1/\lambda$   shows that 
\begin{equation}
    \Gamma =  \Delta t^2 - \dfrac{1}{\lambda}   2h^2 \dfrac{ \widetilde{(h\eta)} }{h\eta} ( 1  - \dfrac{2h^ 2}{\Delta t^2 } \dfrac{1}{\lambda}+ \dfrac{ 6 h^4}{ \Delta t^{4}} \dfrac{1}{\lambda^2 } + \mathcal{O}( \dfrac{1}{\lambda^3})  ).
\end{equation}
Since $\Delta t$ is independent on $\lambda$,
we  have $\Gamma/\Delta t^2 = 1 + \mathcal{O}(1/\lambda)$ from which we can
deduce the estimate on $\Delta t^2\kappa$.
\end{proof}
\end{proposition}

\begin{remark} 
Note that  
for the semi-implicit scheme we  will use $\Delta t \approx \Delta x/ \max (|u| + \sqrt{gh})$  which is indeed independent of $\lambda$. 
%$$
%\overline{\Delta x} 
%  =  \dfrac{h\max_{x} |u|}{ \sqrt{\lambda}} ( 1 + \dfrac{c^2}{|\alpha| \lambda \eta})^{1/4}
%$$
%The above $\lambda^{-1/2}$ estimate is in reality quite rough, due to
%the  use  of two bounds from below in the proof. We tried  keeping 
%the term $b_1^{-2}$ in, however this adds and additional dependence on the time/mesh
%step in the proof which is much harder to study under reasonable assumptions. 
In practice, we have found $\kappa >0$ on all meshes used,  for all benchmarks, and for all the values of $\lambda$ tested.
\end{remark}

%\subsubsection{Including the bathymetry}
%\red{Not at the moment!}

\subsection{High order generalization in time} \label{sec:high_order_time}
As done in the first-order case and following the approach of \cite{MaccaRussoBumi}, we rewrite system \eqref{hSGN_eq} in the ODE form \eqref{eq:ODE}--\eqref{eq:basicSI}. The semi-implicit semi-discrete scheme can be written as
\begin{equation}
    \label{ode_bosca_or2}
    \frac{\text{d}U}{\text{d}t} = \tilde{H}(U_E, U_I),
\end{equation}
with
%\begin{equation}
%    \label{ode_form_part_or2}
%    \tilde{H}(U_E, U_I) = 
%    \begin{bmatrix}
%        & - D_x((hu)_I) \\
%       -\hat{D}_x((hu^2)_E)  & -D_x(p_I) \\
%       -\hat{D}_x((hu\eta)_E) &  -(hw)_I \\
%       -\hat{D}((huw)_E) & \lambda\Big(\frac{h\eta}%{h^2}\Big)_I
%    \end{bmatrix}.
%\end{equation}
\begin{equation}
    \label{ode_form_part_or2}
    \tilde{H}(U_E, U_I) = 
    \begin{bmatrix}
        & - \partial_x((hu)_I) \\
       -\partial_x((hu^2)_E)  & -\partial_x(p_I) \\
       -\partial_x((hu\eta)_E) &  -(hw)_I \\
       -\partial_x((huw)_E) & \lambda\Big(\frac{h\eta}{h^2}\Big)_I
    \end{bmatrix}.
\end{equation}
Based on this formulation, we apply an IMEX scheme to system \eqref{ode_form_part_or2}.  
The general procedure to advance the numerical solution from time $t^n$ to $t^{n+1}$ using an $s$-stage Runge–Kutta IMEX method is as follows:
\begin{itemize}
    \item \textbf{Stage values:} For $i = 1, \ldots, s$, compute
    \begin{align}
        U^{(i)}_E & = U^n + \Delta t \sum_{j=1}^{i-1} a_{i,j}^E H\left(U_E^{(j)}, U_I^{(j)}\right), \\
        U^{(i)}_I & = U^n + \Delta t \left( \sum_{j=1}^{i-1} a_{i,j}^I H\left(U_E^{(j)}, U_I^{(j)}\right) + a_{i,i}^I H\left(U_E^{(i)}, U_I^{(i)}\right) \right).
    \end{align}
    
    \item \textbf{Numerical solution:} \\
    $U^{n+1} = U_I^{(s)}.$
\end{itemize}

\begin{remark}
Note that if a system of the form \eqref{ode_bosca_or2} is autonomous (i.e., the right-hand side does not explicitly depend on time), and the $s$-stage double Butcher tableau has identical $b$ coefficients, then the evolution requires only $s$ evaluations of the function $H$. Furthermore, if the last row of the matrix $A$ coincides with the weights $b$ (i.e., the implicit tableau defines a \emph{stiffly accurate} scheme), then the numerical solution corresponds to the last stage of the implicit part \cite{Boscarino-Filbet}.
\end{remark}

Here, we consider the IMEX scheme defined by the following double Butcher tableau \cite{Macca2024}:
\begin{equation}
\label{tableau}
\begin{array}{c|cc}
     & 0 &  \\
    c & c & 0\\ \hline
     & 1 - \gamma & \gamma
\end{array}
\hspace{3cm}
\begin{array}{c|cc}
    \gamma & \gamma &  \\
    1 & 1 - \gamma & \gamma\\ \hline
     & 1 - \gamma & \gamma
\end{array}
\end{equation}
where $\gamma = 1 - \frac{1}{\sqrt{2}}$ and $c = \frac{1}{2\gamma}$.

In our case, applying the scheme defined in \eqref{tableau}, we obtain:
\begin{enumerate}
    \item $U_E^{(1)} = U^n;$
    \item $U_I^{(1)} = U^n + \Delta t \gamma H(U_E^{(1)}, U_I^{(1)});$
    \item $U_E^{(2)} = U^n + \Delta t c H(U_E^{(1)}, U_I^{(1)});$
    \item $U_I^{(2)} = U^n + \Delta t (1 - \gamma) H(U_E^{(1)}, U_I^{(1)}) + \Delta t \gamma H(U_E^{(2)}, U_I^{(2)});$
    \item $U^{n+1} = U_I^{(2)}.$
\end{enumerate}

\begin{remark}
We observe that $U_E^{(2)}$, $U_I^{(2)}$, and $U_I^{(1)}$ share a common term. Hence, steps 3 and 4 can be equivalently rewritten in efficient form as:
\begin{align}
    U_E^{(2)} &= \left(1 - \frac{c}{\gamma}\right) U^n + \frac{c}{\gamma} U_I^{(1)}, \\
    U_I^{(2)} &= \left(1 - \frac{1 - \gamma}{\gamma} \right) U^n + \frac{1 - \gamma}{\gamma} U_I^{(1)} + \Delta t \gamma H(U_E^{(2)}, U_I^{(2)}).
\end{align}
\end{remark}

\subsection{Stability conditions} \label{ssec_stability}

% \section{Stability conditions} \label{ssec_stability}

% \subsection{Linear stability analysis in time}

% \subsection{Fully discrete time step constraints}

For an explicit scheme, the time step $\Delta t$ is determined at stage $t^n$ by the Courant-Friedrichs-Lewy (CFL) condition as $\Delta t = {\rm CFL}\,\Delta x/\max_\Omega\mu$, where the CFL restriction is governed by the maximum spectral radius $\mu$ of the matrix $A(U)$ defining the hyperbolic system 
\begin{equation}
    \textrm{CFL} = \frac{\mu^n_{\rm max} \Delta t^n}{\Delta x} \leq C_{\textrm{ex}},
    \label{CFL_ex}
\end{equation}
being $\mu_{\rm max} = \max_{\Omega}(u + \sqrt{gh + \lambda/3\eta^2/h^2})$, and $C_{\textrm{ex}}$ is a constant close to one.

In the case of our semi-implicit scheme, we empirically derive the following stability condition based on the velocity and not the celerity:
\begin{equation}\label{eq:MCFL}
    \textrm{MCFL} = \frac{u^n_{\rm max} \Delta t^n}{\Delta x} \leq C_{\textrm{im}},
\end{equation}
where $u_{\rm max} = \max_\Omega{|u|}$ and $C_{\textrm{im}} \approx 0.75$.

This condition is considerably less restrictive than \eqref{CFL_ex}, as the condition for the classical CFL can be expressed as
\begin{equation}
    \textrm{CFL} \ =\  \frac{\mu^n_{\rm max}}{u^n_{\rm max}}\, \textrm{MCFL} \ \leq\  \frac{\mu^n_{\rm max}}{u^n_{\rm max}} C_{\rm im},
\end{equation}
and $\mu^n_{\rm max}/u^n_{\rm max}\gg 1$ for small Froude numbers.

In the classical Shallow Water equations, the CFL condition is typically governed by the propagation speed of gravity waves, leading to a constraint of the form
\begin{equation}
    \Delta t \leq \mathrm{CFL} \frac{\Delta x}{u+\sqrt{g h_{\max}}}.
\end{equation}

However, for the hSGN model, the presence of higher-order dispersive terms modifies the characteristic wave speeds, potentially leading to a more restrictive CFL condition. The additional dispersive effects introduce new characteristic velocities that must be taken into account when selecting the time step. Consequently, the stability constraint is generally stricter than in the classical Shallow Water framework. This difference is particularly evident in numerical experiments, where maintaining stability in hSGN formulations often requires a lower effective CFL number than in Shallow Water models. As a result, the semi-implicit approach provides a significant advantage by allowing larger time steps while preserving stability and accuracy.

\section{Spatial discretization}
\subsection{Explicit hSGN}
The first-order explicit numerical scheme for system~\eqref{hSGN_eq} reads as follows:

\begin{equation}
\label{GN_hyp_eq_explicit}
\begin{cases}
h_i^{n+1} = h_i^n - \Delta t \, D_x(q^n)_i, \\

q_i^{n+1} = q_i^n 
- \Delta t \, \hat{D}_x(q^n u^n)_i 
- \Delta t \, a^2(h^n,\eta^n) D_x(h^n)_i 
- \Delta t \, \lambda \alpha(h^n,\eta^n) D_x(h^n\eta^n)_i, \\

(h\eta)_i^{n+1} = (h\eta)_i^n 
- \Delta t \, \hat{D}_x(h^n u^n \eta^n)_i 
+ \Delta t \, h_i^n w_i^n, \\

(hw)_i^{n+1} = (hw)_i^n 
- \Delta t \, \hat{D}_x(h^n u^n w^n)_i 
- \Delta t \, \lambda \left( \frac{(h\eta)_i^{n}}{(h_i^{n})^2} - 1 \right).
\end{cases}
\end{equation}
The differential operators $D_x$ and $\hat{D}_x$ applied to a given flux function $F(U)$ are respectively defined as:
\begin{itemize}
    \item $D_x(F)_i = \frac{F_{i+\ha} - F_{i-\ha}}{\dx},$ in which $F_{i\pm\ha} = \ha(F(U_{i\pm1}) + F(U_i));$
    \item $\hat{D}_x(F)_i = \frac{\hat F_{i+\ha} - \hat F_{i-\ha}}{\dx},$ where  $\hat F_{i+\ha} = \ha\Bigl( F(U_{i+\ha}^{-}) + F(U_{i+\ha}^{+}) - \alpha_{i+\ha}\big(U_{i+\ha}^{+} - U_{i+\ha}^{-}\big)\Bigr)$ is the Rusanov flux and $\alpha_{i+\ha}$ is related to the eigenvalues of the explicit sub system. In our case,  $\alpha \approx |u|.$
\end{itemize}

The  values $U_{i+\frac{1}{2}}^\pm$ are computed via 
a  conservative  component-wise  linear reconstruction  \cite{Toro2009}:
\begin{equation}
    v_{i\pm\frac{1}{2}}^{\mp} = {\mathbf{v}}_i \pm v_i' \frac{\Delta x}{2},
\end{equation}
where the reconstructed slope $v_i'$ is given by
\begin{equation}
    v_i' = \frac{\mathbf{v}_{i+1} - \mathbf{v}_{i-1}}{2\Delta x},
\end{equation}
and $\mathbf{v}_i$ represents the cell average $\int_{I_i} v(x,t)dx$. 

The second order in time has been obtained by Heun time procedure \cite{Toro2009}.

\subsection{Spatial discretization of the SI-hSGN}
In the semi-implicit case, we use the 
 semi-discrete first-order in time system\footnote{Since we are focused on the spatial discretization and in order to simplify the readability  $\bar{U} =U(x,t^{n+1})$, while $U = U(x,t^n)$.} as:
\begin{equation}
    \label{GN_hyp_eq_semi}
    \begin{cases}
        \bar{h} = h - \dt D_x(\bar{q}) \\ 
        \bar{q} = q - \dt \hat{D}_x(qu) - \dt a^2D_x(\bar{h}) - \dt\lambda\alpha D_x(\bar{h\eta})  \\
        \bar{h\eta} = h\eta - \dt \hat{D}_x(hu\eta) + \dt\bar{hw} \\
        \bar{hw} = hw - \dt \hat{D}_x(huw) - \dt\lambda\Bigl(\frac{\bar{h\eta}}{\bar{h}^2}-1\Bigr),
    \end{cases}
\end{equation}
where the differential operators $D_x$ and $\hat{D}_x$ have been defined above. 

The fully-discrete second order scheme can be derived combining the previous Sec.\ref{sec:splitting}-\ref{sec:high_order_time}. In particular, the fully-discrete first order scheme is: 
\begin{equation}
\label{GN_hyp_eq_semi_fully_discrete}
\begin{cases}hw_i^\dagger  = hw_i - \frac{\dt}{\dx}\Bigl(F_{i+\ha}^{huw} -F_{i-\ha}^{huw}\Bigr) + \lambda\dt \\ 
h\eta_i^\dagger  = h\eta_i - \frac{\dt}{\dx}\Bigl(F_{i+\ha}^{hu\eta} -F_{i-\ha}^{hu\eta}\Bigr) + \dt hw_i^\dagger \\
q_i^\dagger = q_i -\frac{\dt}{\dx}\Bigl(F_{i+\ha}^{qu} -F_{i-\ha}^{qu}\Bigr) \\
\beta_1  := 1 + \frac{\lambda\dt^2}{h^2} \\
\beta_2  := \frac{2\lambda\dt^2}{\beta_1^2h^3} \\
\beta  := a^2 + \lambda\alpha\beta_2h\eta^\dagger \\ 
\delta  := \frac{\alpha}{\beta_1} \\ 
 \Delta[h]_i  :=  
 h_i\Bigl(1 + \frac{\dt^2}{\dx^2}(\beta_{i+\ha} + \beta_{i-\ha})\Bigr) - \frac{\dt^2}{\dx^2}\beta_{i+\ha}h_{i+1} - \frac{\dt^2}{\dx^2}\beta_{i-\ha}h_{i-1}\\
h_i^\dagger  = h_i - \frac{\dt}{2\dx}\Bigl(q_{i+1}^\dagger - q_{i-1}^\dagger \Bigr) + \frac{\lambda\dt^2}{2\dx^2}\Bigl(\delta_{i+\ha}(h\eta_{i+1}^\dagger - h\eta_{i}^\dagger) -\delta_{i-\ha}(h\eta_{i}^\dagger - h\eta_{i-1}^\dagger) \Bigr) \\  
\Delta[\bar{h}]_i  = h_i^\dagger \\
\bar{q}_i = q_i^\dagger - \frac{\dt}{2\dx}\beta_i\Bigl( \bar{h}_{i+1} - \bar{h}_{i-1}\Bigr) - \frac{\lambda\dt}{2\dx}\delta_i\Bigl( h\eta_{i+1}^\dagger - h\eta_{i-1}^\dagger\Bigr) \\
\bar{h\eta}_i  = \frac{h\eta_i^\dagger}{ 1 + \frac{\lambda\dt^2}{\bar{h}_i^2}} \\
\bar{hw}_i  = hw_i^\dagger - \lambda\dt\frac{\bar{h\eta}_i}{\bar{h}_i^2}.
\end{cases}
\end{equation}
% \begin{equation}
%     \label{GN_hyp_eq_semi_fully_discrete}
%     \begin{cases}
%     \bar{h}_i+\frac{\dt}{2\dx}\Bigl(\bar q_{i+1}-\bar q_{i-1}\Bigr) \\
%         \bar{q}_i = q_i - \frac{\dt}{\dx}\Bigl(\hat F^{qu}_{i+\ha} - \hat F^{qu}_{i-\ha} \Bigr) - \frac{\dt}{2\dx} a^2\textcolor{red}{(h_i\eta_i)}\Bigl(\bar{h}_{i+1}-\bar{h}_{i-1}\Bigr) - \frac{\dt}{2\dx}\lambda\alpha\textcolor{red}{(h_i\eta_i)}\Bigl(\bar{h\eta}_{i+1}-\bar{h\eta}_{i-1}\Bigr)  \\
%         \bar{h\eta}_i = h\eta_i - \frac{\dt}{\dx} \Bigl(\hat F^{hu\eta}_{i+\ha}-\hat F^{hu\eta}_{i-\ha}\Bigr) + \dt\bar{hw}_i \\
%         \bar{hw}_i = (hw)_i - \frac{\dt}{\dx}\Bigl(\hat F^{huw}_{i+\ha}-\hat F^{huw}_{i-\ha}\Bigr) - \dt\lambda\Bigl(\frac{\bar{h\eta}_i}{\bar{h}^2_i}-1\Bigr).
%     \end{cases}
% \end{equation}

% prima risolviamo
% $$
% \begin{aligned}
%  \bar{h}_i - \Delta t^2 \Bigl( \kappa_{i+1/2} \dfrac{\bar{h}_{i+1} -\bar{h}_i}{\Delta x} -
%        \kappa_{i-1/2} \dfrac{\bar{h}_{i} -\bar{h}_{i-1}}{\Delta x} \Bigr)= & h_i -  \frac{\dt}{2\dx}\Bigl(q_{i+1}-q_{i-1}\Bigr) \\
%        + &\lambda \Delta t^2 \Bigl( (\alpha b_1)_{i+1/2} \dfrac{  \widetilde{h\eta}_{i+1} -\widetilde{h\eta}_i}{\Delta x} -
%        (\alpha b_1)_{i-1/2} \dfrac{\widetilde{h\eta}_{i} -\widetilde{h\eta}_{i-1}}{\Delta x} \Bigr)
%        \end{aligned}
% $$
% con appropriate definitions of $\widetilde{h\eta}$,
% poi si fanno  gli  altri update  partendo da  $\bar{h\eta}$ etc. 

\subsection{Standard SGN equations}
To approximate \eqref{eq:SGN_flat} numerically,  we will exploit the elliptic-hyperbolic form initially due to \cite{Filippini}, and used for the Green-Naghdi and other dispersive
systems in \cite{Torlo,Cauquis,Jouy,Kazolea,GR25}. For solution purposes, we recast the system as follows 
\begin{equation}
\begin{split}
h_t + & (hu)_x =0\\
(hu)_t + & (hu^2 + gh^2/2)_x  = h\phi, 
\end{split}
\end{equation}
and now we look at the constraint that $\phi$ must verify to satisfy the Serre-Green-Naghdi equation.
We note first that, using the mass equation we have
\begin{equation}
    \phi =  \dfrac{Du}{Dt} + g  h_x,
\end{equation}
so using the last formula of the non-hydrostatic pressure we see that
\begin{equation}
    p_{\textsf{nh}}=  - \dfrac{h^3}{3} \left(  \phi - gh_x\right)_x +  2\dfrac{h^3}{3} (u_x)^2.
\end{equation}
For the last system to be  equivalent to the Serre-Green-Naghdi equations
we must have
\begin{equation}
    h\phi    = -(p_{\textsf{nh}} )_x = ( \dfrac{h^3}{3} \left(  \phi - gh_x\right)_x - 2\dfrac{h^3}{3} (u_x)^2)_x,
\end{equation}
which is equivalent to the constraint
\begin{equation}
    h\phi  -  ( \dfrac{h^3}{3}    \phi_x )_x =- ( \dfrac{h^3}{3}  gh_{xx}+ 2\dfrac{h^3}{3} (u_x)^2)_x.
\end{equation}
So we can now solve the problem as follows
\begin{equation}
\label{eq:SGN_standard}
\begin{split}
h\phi  - & ( \dfrac{h^3}{3}    \phi_x )_x =- ( \dfrac{h^3}{3}  gh_{xx}+ 2\dfrac{h^3}{3} (u_x)^2)_x,\\
h_t + & (hu)_x =0,\\
(hu)_t + & (hu^2 + gh^2/2)_x  = h\phi,  
\end{split}
\end{equation}
marching in time by solving first for $\phi$ and then solving the shallow water equations with the additional source.\\

The above procedure can be easily generalized to \eqref{SGNb_eq}-\eqref{SGNb_p}, and leads to the coupled system 
\begin{equation}
\label{eq:SGNb_standard}
\begin{split}
h\phi  - & ( \dfrac{h^3}{3}    \phi_x )_x  
+ \kappa(h,b) \phi =- ( \dfrac{h^3}{3}  g\zeta_{xx}
+ 2\dfrac{h^3}{3} (u_x)^2 + \rho(h,u,b) )_x  -h\mathcal{Q}(h,u,b)b_x,\\
h_t + & (hu)_x =0,\\
(hu)_t + & (hu^2 + gh^2/2)_x  = h\phi,   
\end{split}
\end{equation}
having set
\begin{equation}
    \begin{aligned}
\kappa(h,b):= & \big(\dfrac{h^2b_x}{2}\big)_x + 3\dfrac{h b_x^2}{4},\\
\rho(h,u,b):= &-\dfrac{h^2b_x}{2}g\zeta_x  + h^3 u^2b_{xx},& \\
\mathcal{Q}(h,u,b):= & 
\dfrac{h}{2}g\zeta_{xx} - 3\dfrac{g b_x}{4}  \zeta_x + h u_x^2 +\dfrac{3}{2} h u^2 b_{xx}.
\end{aligned}
\end{equation}
The above can be marched in time by first inverting the second order PDE for $\phi$, and then performing the shallow water updated enhanced with the $h\phi$ source.
The interested reader is referred to \cite{Torlo,Cauquis,Jouy,Kazolea,GR25} for more details.

\subsection{Explicit Scheme for the standard SGN equations}\label{ssec:explicit_SGN_original}
The fully-discrete explicit scheme for the standard Serre-Green-Naghdi is given by
\begin{equation}
\left\{
\begin{aligned}
& h^{n+1}_i = h_i^n - \frac{\dt}{\dx}\Bigl(\hat F^{hu}_{i+\ha}-\hat F^{hu}_{i-\ha}\Bigr), \\
& hu^{n+1}_i = hu^{n+1}_i - \frac{\dt}{\dx}\left(\hat F^{hu^2 + g h^2/2}_{i+\ha} - \hat F^{hu^2 + g h^2/2}_{i-\ha}\right) + \dt h\phi_i^n, \\
& h\phi_i^{n} - \frac{1}{3\dx^2}\left(h_{i+\ha}^3(\phi_{i+1}^n-\phi_{i}^n) - h_{i-\ha}^3(\phi_{i}^n-\phi_{i-1}^n) \right) = h_i^*,
\end{aligned}
\right.
\end{equation}
where 
\begin{align}
     h_i^* = -\frac{1}{6\dx^3}\Bigl( & gh_{i+1}^n\left(h_{i+2}^n - 2h_{i+1}^n + h_{i}^n\right) + \frac{(h^3)_{i+1}^n}{2}\left(u_{i+2}^n - u_{i}^n\right)^2 + \\ & - gh_{i-1}^n\left(h_{i}^n - 2h_{i-1}^n + h_{i-2}^n\right) - \frac{(h^3)_{i-1}^n}{2}\left(u_{i}^n - u_{i-2}^n\right)^2   \Bigr).  
\end{align}
In order to achieve second order in time, the Heun method has been applied.

\section{Numerical tests}
This section presents numerical tests to assess the reliability of the proposed model and its semi-implicit numerical treatment, highlighting its accuracy, ability to capture dispersive effects, and overall performance. First, we verify accuracy using a solitary wave (soliton) with a known exact solution, allowing us to quantify numerical errors. Then, we examine the model’s ability to capture dispersive effects, a crucial aspect in wave propagation problems. Finally, we consider the evolution of Favre waves to further evaluate the method’s effectiveness in complex scenarios. These experiments demonstrate the advantages of the proposed approach in terms of accuracy, stability, and computational efficiency.

The time step $ \Delta t^n $ is determined by
\begin{equation}
    \Delta t^n = \mathrm{CFL}\frac{\Delta x}{\mu_{\max}^n}, 
\end{equation}
where $ \mu_{\max}^n $ denotes the maximum spectral radius of the Jacobian matrices $ \partial F / \partial U $ across cells. However, $ \Delta t^n $ must also satisfy the material CFL condition (MCFL, see \eqref{eq:MCFL}), which we have empirically determined to be not greater than 0.5 \cite{Bosca2024}. This introduces variability in the CFL condition, allowing it to differ from one time step to another.

The time step is thus computed in accordance with both CFL conditions, one associated with the free surface waves and the other with the material wave. Specifically, for each $ n $,
\begin{equation}
\bar{\Delta t}^n = \mathrm{CFL}\frac{\Delta x}{\lambda_{\max}^n},    
\end{equation}
after which the $\bar{\Delta t}$-material CFL condition is evaluated as $ {\rm MCFL^*} = \frac{u_{\max} \bar{\Delta t}^n}{\Delta x} $. If $ {\rm MCFL^*} \le {\rm MCFL} $, then $ \Delta t = \bar{\Delta t} $; otherwise, $ \Delta t = {\rm MCFL} \frac{\Delta x}{u_{\max}} $. This approach means that the CFL condition may vary within each time step. Consequently, when the CFL condition is not constant across iterations, we include a figure to display its evolution when significant.

Furthermore, we introduce distinct CFL numbers tailored to the characteristic speeds of each numerical scheme. For both the semi-implicit and explicit treatment of the hyperbolic Serre–Green–Naghdi model, the CFL condition is expressed as the ratio \(\Delta t/\Delta x\) multiplied by the maximal eigenvalue \(\nu_{\max}\), which itself depends on the parameter \(\lambda\). Concretely, we set

\begin{equation}
\mathrm{CFL}_{\mathrm{IMEX}} = \frac{\Delta t}{\Delta x}\,\nu_{\max},
\qquad
\mathrm{CFL}_{\mathrm{EX}} = \frac{\Delta t}{\Delta x}\,\nu_{\max},
\end{equation}
with
\begin{equation}
\nu_{\max} = u + \sqrt{g h + \frac{\lambda}{3}\left( \frac{\eta}{h} - 1 \right)}.
\end{equation}

Although the definitions of $\mathrm{CFL}_{\mathrm{IMEX}}$ and $\mathrm{CFL}_{\mathrm{EX}}$ are formally identical, they are subject to different stability constraints. The explicit CFL, $\mathrm{CFL}_{\mathrm{EX}}$, must satisfy the classical condition $\mathrm{CFL}_{\mathrm{EX}} < 1$ to ensure numerical stability. In contrast, the IMEX CFL, $\mathrm{CFL}_{\mathrm{IMEX}}$, is governed by the material stability condition and can typically exceed one. This distinction reflects the enhanced stability properties of IMEX schemes in the presence of stiff terms.

This formulation captures how the maximum wave speed in the hyperbolic system varies with \(\lambda\) and the flow state, ensuring that both semi-implicit and explicit schemes respect stability constraints dictated by the most restrictive local propagation speed.

By contrast, for the standard Serre–Green–Naghdi equations we employ a more classical CFL number defined through the simpler characteristic speed \(\sigma_{\max}\). Specifically,

\begin{equation}
    \mathrm{CFL}_{\rm SGN} = \frac{\Delta t}{\Delta x}\,\sigma_{\max},
\qquad
\sigma_{\max} = u + \sqrt{g h}.
\end{equation}

Here, \(\sigma_{\max}\) corresponds to the maximal signal velocity in the absence of the extra dispersive correction associated with \(\lambda\). By distinguishing these CFL definitions, we align the time-step restrictions with the physics and numerical treatment of each variant of the Serre–Green–Naghdi model, thereby promoting both stability and efficiency across the different schemes.

In all the numerical tests, $g$ is the gravitational constant $g=9.81\ m^2/s$ and the numerical errors have been computed in L$^1$ norm as: 
\begin{equation}
     \textrm{Error} = \frac{||u_{\rm num} - u_{\rm exact}||_1}{||u_{\rm num}||_1}.
\end{equation}

\section{Breakdown of a Gaussian bell}

\subsection{Flat bathymetry}
Dispersive waves play a crucial role in many physical systems, particularly in shallow water models where higher-order terms introduce significant dispersive effects. These small-amplitude oscillations, commonly referred to as dispersive ripples, typically emerge behind a primary wave due to the interaction between nonlinearity and dispersion. The numerical resolution of such phenomena presents several challenges: high-frequency components may lead to spurious oscillations if not properly captured, whereas excessive numerical dissipation can suppress these effects entirely, distorting the physical behavior of the solution.

\begin{figure}[!ht]
    \centering
    \begin{subfigure}[b]{0.85\textwidth}
        \centering
        \includegraphics[trim={2cm 0 2.5cm 0},clip,width=0.99\textwidth]{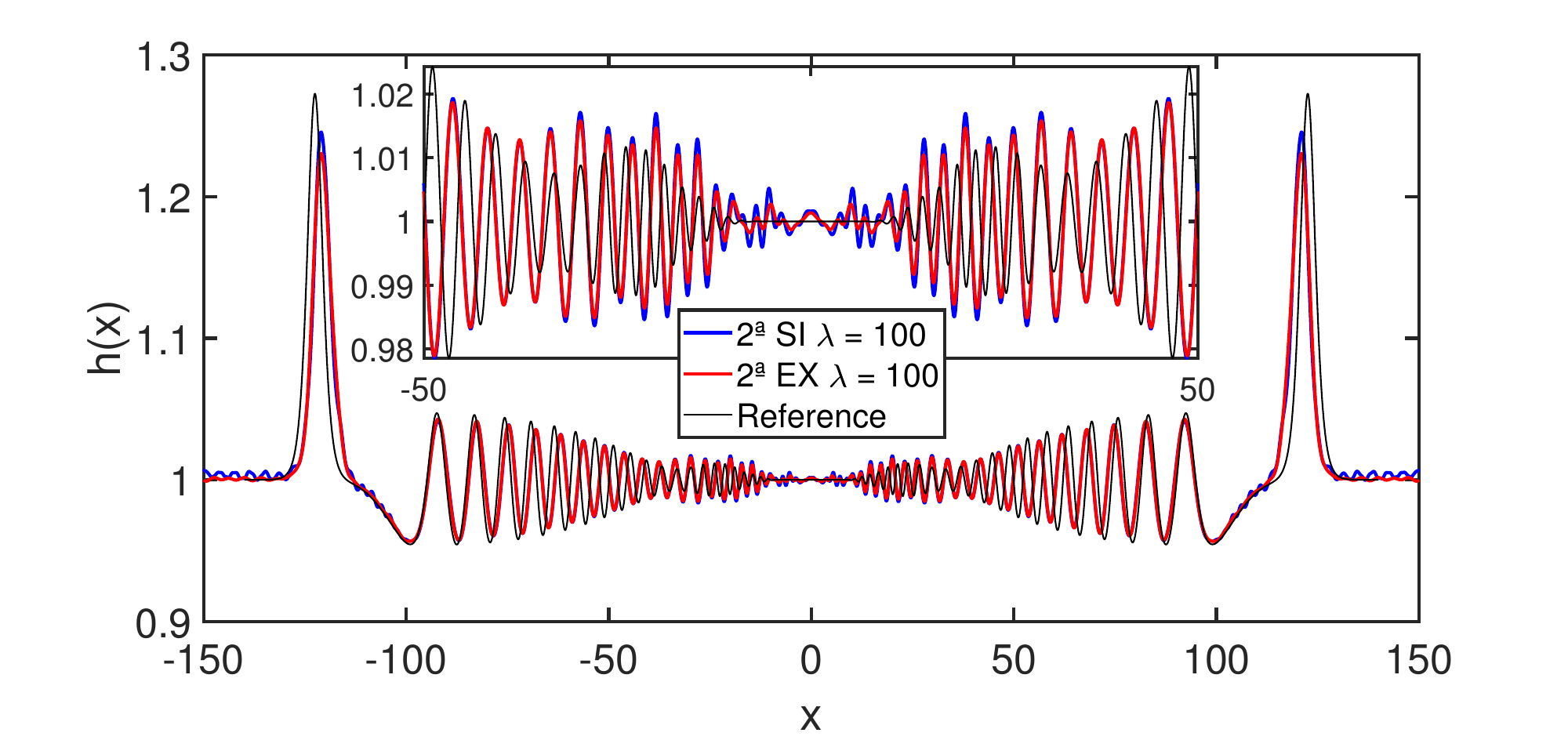}
        \caption{$\lambda = 100$.}      
        \label{ssec:sfig:dispersion:lambda100}
    \end{subfigure}
    \\[-1.5pt]
    %\hfill
    \begin{subfigure}[b]{0.85\textwidth}
        \centering
        \includegraphics[trim={2cm 0 2.5cm 0},clip,width=0.99\textwidth]{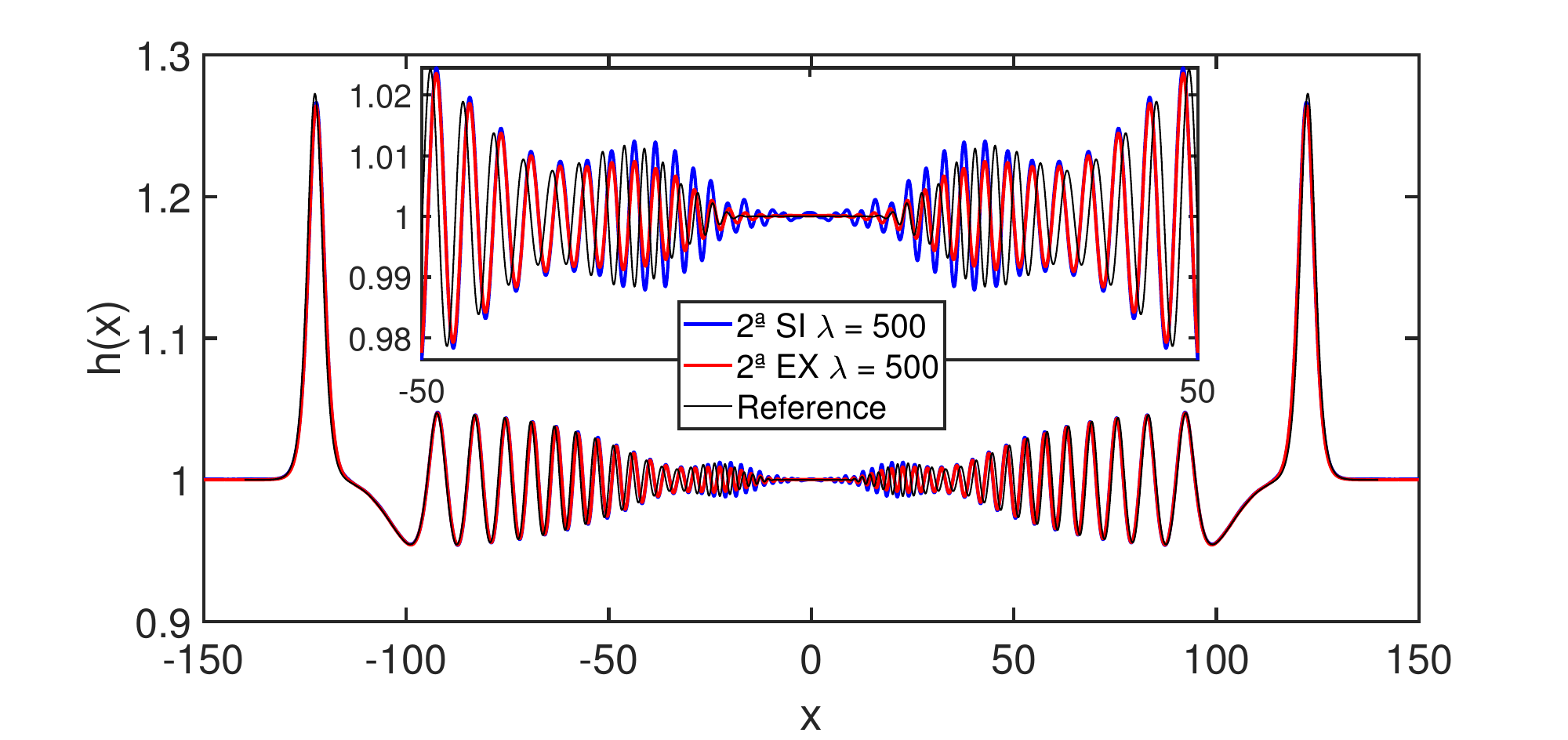}
        \caption{$\lambda = 500$.}
\label{ssec:sfig:dispersion:lambda500}
    \end{subfigure}    
        \caption{Numerical solutions for dispersive effects using second-order semi-implicit and explicit schemes for different values of $\lambda$, respectively 100 and 500. The final time is $t = 35$. The reference solution has been obtained with the fourth-order energy preserving scheme for the 
    standard SGN model\eqref{eq:SGN_standard}  developed in \cite{rr25}.}
    \label{ssec:fig:dispersion}
\end{figure}
\begin{figure}[!ht]
    \centering
    \begin{subfigure}[b]{0.85\textwidth}
        \centering
        \includegraphics[trim={2cm 0 2.5cm 0},clip,width=0.99\textwidth]{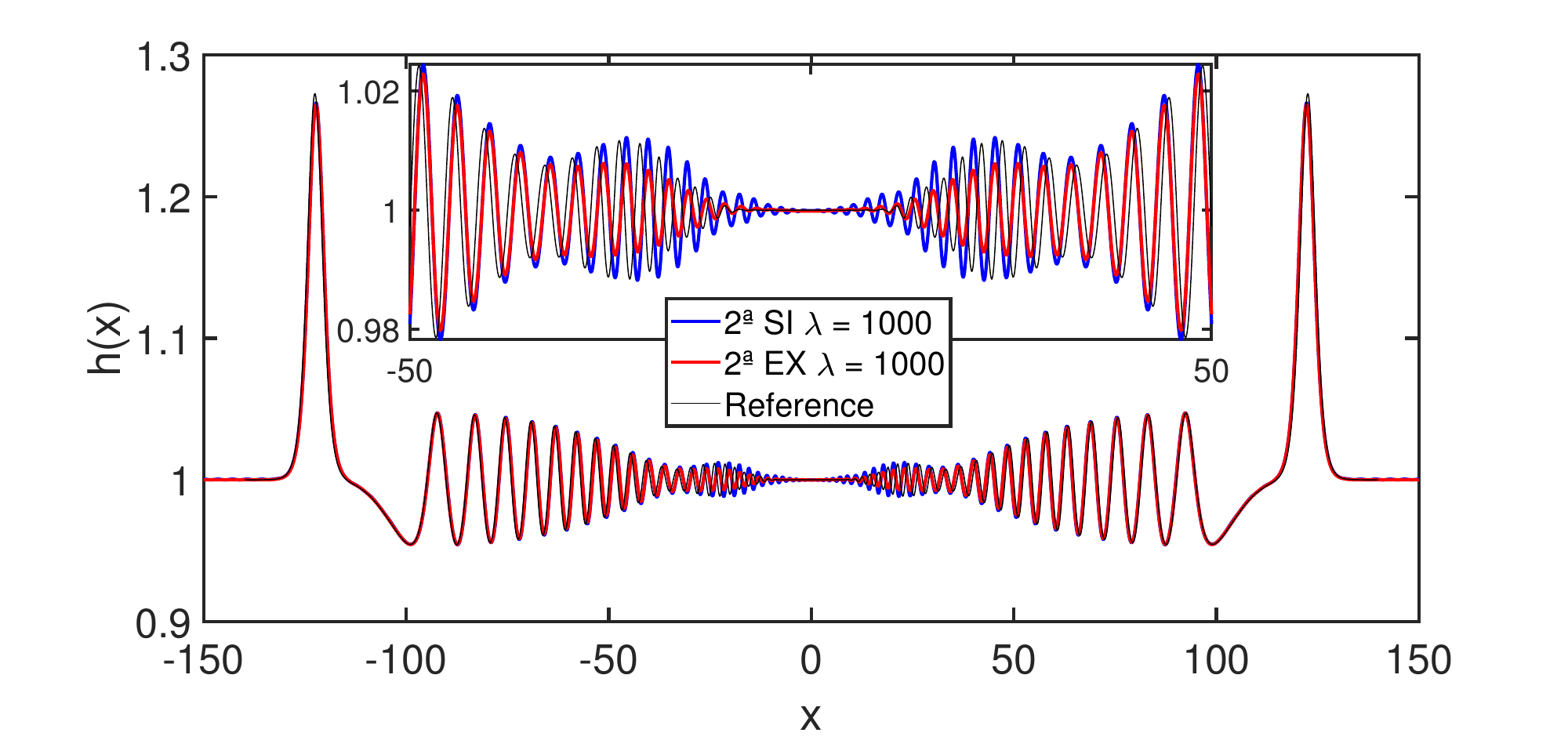}
       \caption{$\lambda = 1000$.}
\label{ssec:sfig:dispersion:lambda1000}
    \end{subfigure} \\[-1.5pt]
    %\hfill
    \begin{subfigure}[b]{0.85\textwidth}
        \centering
        \includegraphics[trim={2cm 0 2.5cm 0},clip,width=0.99\textwidth]{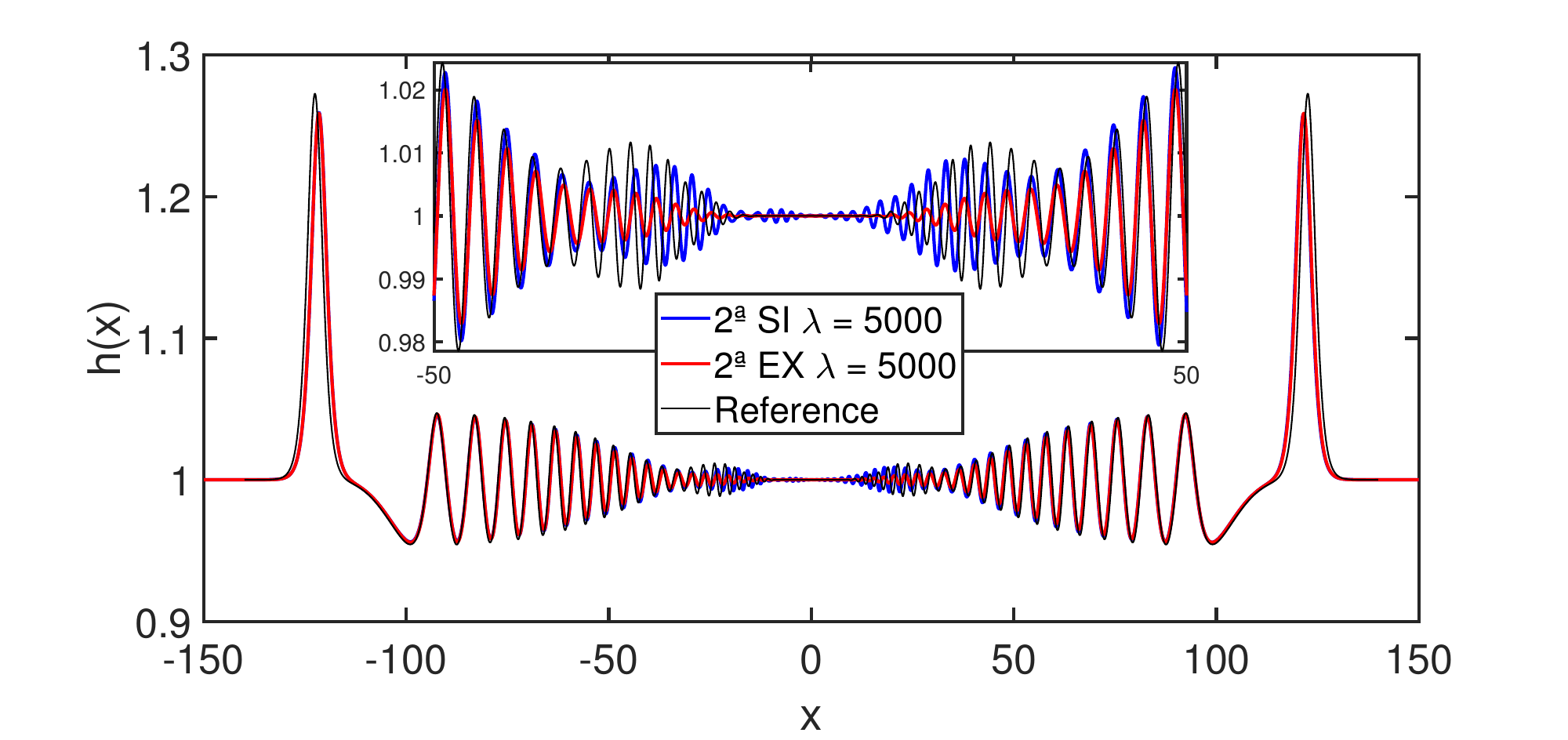}\\[-1.5pt]
        \caption{$\lambda = 5000$.}
\label{ssec:sfig:dispersion:lambda5000}
    \end{subfigure}    
    \caption{Numerical solutions for dispersive effects using second-order semi-implicit and explicit schemes for different values of $\lambda$, respectively 1000 and 5000. The final time is $t = 35$. The reference solution has been obtained with the fourth-order energy preserving scheme for the 
    standard SGN model\eqref{eq:SGN_standard}  developed in \cite{rr25}.}
    \label{ssec:fig:dispersion_2}
\end{figure}

In this test, we evaluate the ability of the first- and second-order semi-implicit schemes to accurately resolve dispersive effects. The initial condition consists of a localized perturbation over a constant background state, given by:
\begin{equation}
    h(x,0) = h_\infty + A e^{-x^2/20}, \quad h_\infty = 1, \quad A = 1.
\end{equation}
The computational domain is set to $[-200,200]$ with a uniform grid of $N = 5000$ points, and the final time of the simulation is $t = 35$. The numerical stability is ensured with a CFL$_{\rm IMEX}$ number of 2 for $\lambda \in\{100,500,1000,5000\}$. Simultaneously, the CFL$_{\rm SW}$ is, respectively, 1.34, 0.66, 0.48 and 0.22.

Figures~\ref{ssec:sfig:dispersion:lambda100}-\ref{ssec:sfig:dispersion:lambda5000} show the numerical solutions obtained with the and second-order semi-implicit and explicit schemes for different values of $\lambda$ ($500, 1000, 5000, 10000$). The results highlight the improved accuracy of the second-order method in capturing the dispersive ripples while maintaining numerical stability. The reference solution has been obtained with a fourth order energy preserving scheme  by \cite{rr25} applied to the classical SGN equations on 15000 points.

\begin{table}[!ht]
    \numerikNine
    \centering
    \begin{tabular}{|l|cccc|}
        \hline \multicolumn{5}{|c|}{\textbf{CPU time for Dispersive effects}} \\
        \hline
        Scheme & $\lambda = 100$ & $\lambda = 500$ & $\lambda = 1000$ & $\lambda = 5000$ \\
        \hline
        IMEX  & 7.68  & 16.18 & 25.59 & 63.34  \\
        Explicit  & 19.51  & 39.938 & 75.86 & 169.73 \\ 
        SGN & 12.43 &- & - &-\\
        \hline
    \end{tabular}
    \caption{Computational cost (in seconds) of semi-implicit, explicit and Serre-Green-Naghdi scheme  for $\lambda \in\{100,500, 1000,5000\}$ adopting $N=5000$ uniform mesh. The CFL$_{\rm IMEX} = 2$ while the CFL$_{\rm Exp} = 0.4$ and the CFL$_{\rm SGN} = 0.9$.}
    \label{ssec:tab:dispersive:cpu}
\end{table}

As shown in Table \ref{ssec:tab:dispersive:cpu} indicates that the semi-implicit IMEX scheme offers clear computational advantages in the presence of strong dispersive effects (large $\lambda$), while the explicit scheme becomes rapidly less efficient as the required time-step diminishes. The comparison with the standard SGN model at the baseline case highlights the trade-off between model complexity (inclusion of dispersive corrections) and computational cost: even though SGN avoids the parameter $\lambda$, its stability constraint leads to higher CPU time than IMEX at $\lambda=100$, and it cannot adapt to scenarios with varying $\lambda$ in the same framework.

\subsection{Smooth bathymetric hump}

In this test, we assess the performance of the semi-implicit (IMEX) and explicit schemes for the hyperbolic Serre–Green–Naghdi model with smooth (non-flat) bathymetry. The initial water depth and discharge are prescribed by
\[
h_0 = 1, \quad A = 1, 
\]
\[
b(x) = \frac{1}{20}\bigl(h_0 + A e^{-2(x+50)^2}\bigr),
\quad
h(x,0) = h_0 + A e^{-x^2} - b(x),
\quad
u(x,0) = 10^{-2}.
\]
\begin{figure}[!ht]
    \centering
    %\begin{subfigure}[b]{0.495\textwidth}
    %    \centering
    %    \includegraphics[width=\textwidth]{Figures/Dispersive/IC_Dispersive_bathy.png}
    %    \caption{Initial condition.}
    %    \label{ssec:sfig:dispersion:IC_bathy}
   % \end{subfigure}
   % \hfill
    \begin{subfigure}[b]{0.72\textwidth}
        \centering
        \includegraphics[width=\textwidth]{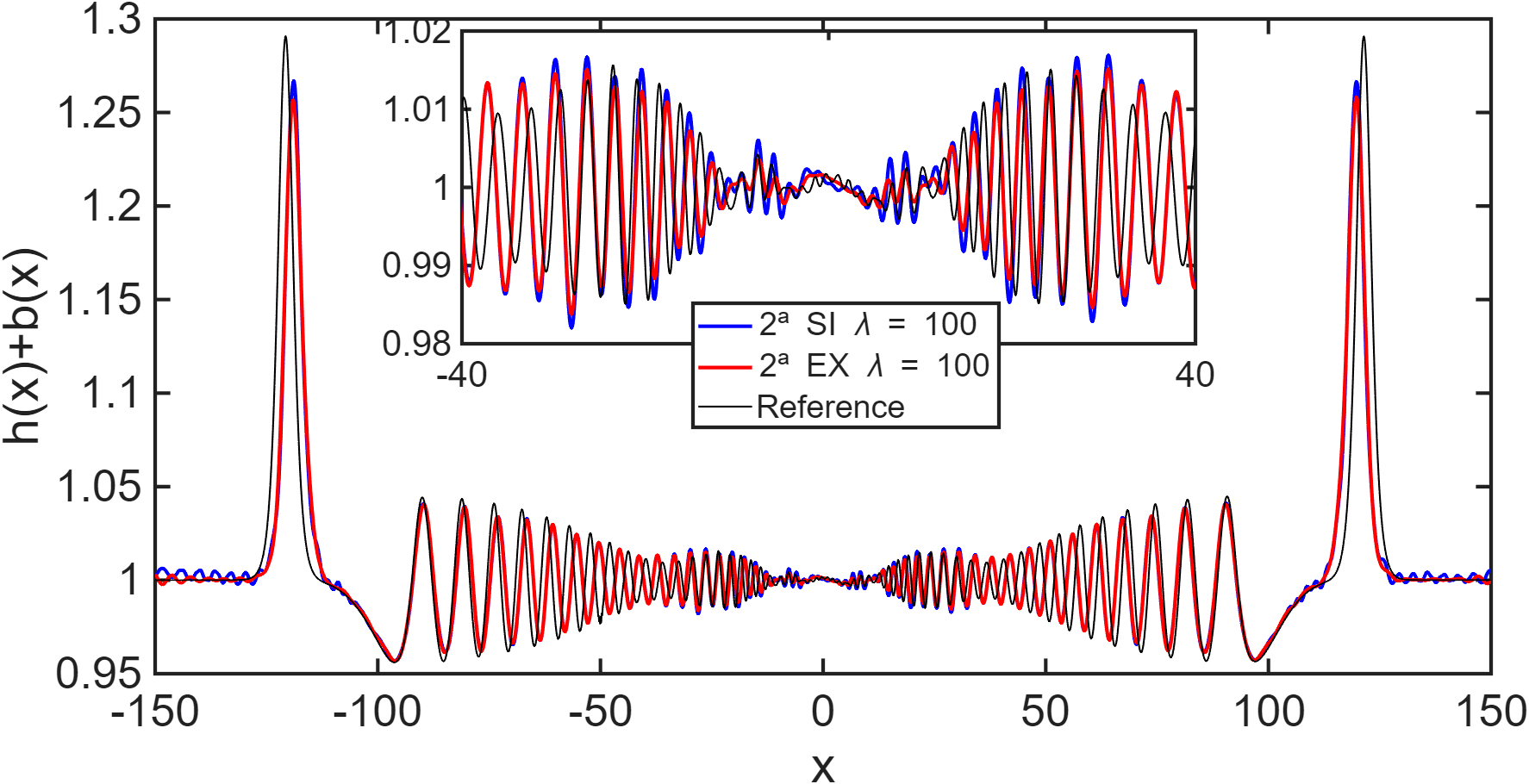}
        \caption{$\lambda = 100$.}
        \label{ssec:sfig:dispersion:lambda100_bathy}
    \end{subfigure}    
    \\[-2.5pt]
    \begin{subfigure}[b]{0.72\textwidth}
        \centering
        \includegraphics[width=\textwidth]{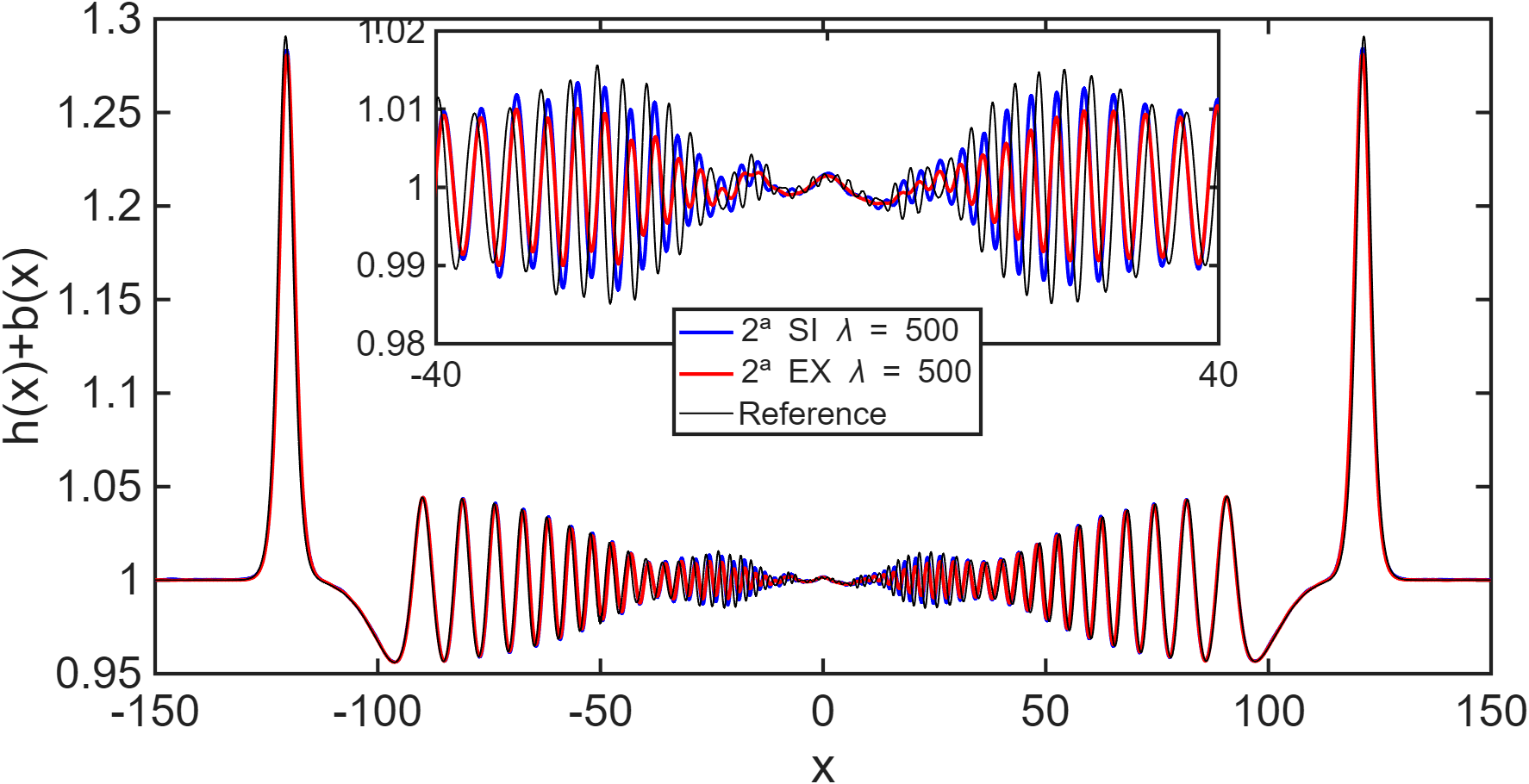}
        \caption{$\lambda = 500$.}
        \label{ssec:sfig:dispersion:lambda500_bathy}
    \end{subfigure}
    \\ [-2.5pt]
    \begin{subfigure}[b]{0.72\textwidth}
        \centering
        \includegraphics[width=\textwidth]{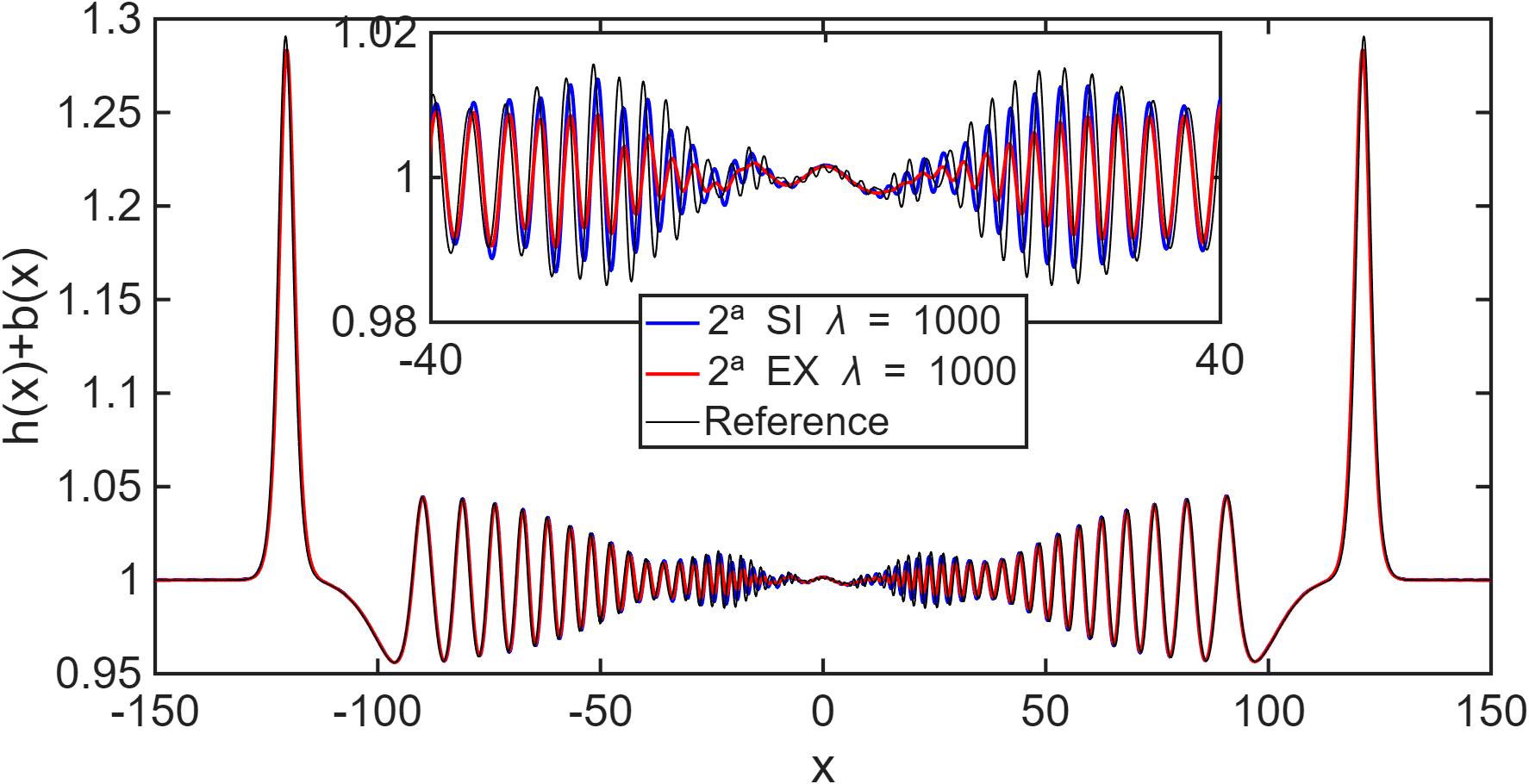}
        \caption{$\lambda = 1000$.}
        \label{ssec:sfig:dispersion:lambda1000_bathy}
    \end{subfigure}    
    \caption{Numerical solutions for dispersive effects using second-order semi-implicit and explicit schemes for different values of $\lambda$, respectively 100, 500 and 1000. The final time is $t = 35$. The reference solution has been obtained with the fourth-order energy preserving scheme for the 
    standard SGN model\eqref{eq:SGN_standard}  developed in \cite{rr25}.}
    \label{ssec:fig:dispersion_bathy}
\end{figure}

Here, $b(x)$ defines a smooth bottom variation and $h(x,0)$ the initial perturbed free-surface above this bathymetry, while $u$ is a small non-zero initial velocity. The computational domain and final time are chosen consistently with previous experiments (e.g., domain $[-150,150]$, final time $t=35$), ensuring that waves interact with the varying bed.

We consider three values of the dispersive parameter $\lambda \in \{100, 500, 1000\}$. For each $\lambda$, the CFL conditions are set analogously to the flat-bottom case: 
\[
\mathrm{CFL}_{\rm IMEX} = 2, 
\qquad
\mathrm{CFL}_{\rm Exp} = 0.4.
\]

Figures \ref{ssec:sfig:dispersion:lambda100_bathy}-\ref{ssec:sfig:dispersion:lambda1000_bathy} display the depth profiles at $t=35$ for semi-implicit and explicit methods. A reference solution obtained   with a fourth order energy preserving scheme  by \cite{rr25} applied to the classical SGN equations on 15000 points is also reported.  The results  highlight how dispersive ripples evolve over uneven bottom and how each scheme captures small-scale oscillations without introducing spurious artifacts or excessive dissipation.

\begin{table}[!ht]
    \centering
    \begin{tabular}{|l|ccc|}
        \hline
        \multicolumn{4}{|c|}{\textbf{CPU time for Dispersive Bathymetry}} \\
        \hline
        Scheme   & $\lambda = 100$ & $\lambda = 500$ & $\lambda = 1000$ \\
        \hline
        IMEX     & 12.73   & 24.78   & 34.63   \\
        Explicit & 36.02   & 72.40   & 109.99    \\
        \hline
    \end{tabular}
    \caption{Computational cost (in seconds) of semi-implicit (IMEX) and explicit schemes for the hyperbolic dispersive model over smooth bathymetry, with $\lambda \in\{100,500,1000\}$, adopting $N=5000$ uniform mesh. The CFL numbers are $\mathrm{CFL}_{\rm IMEX}=2$ and $\mathrm{CFL}_{\rm Exp}=0.4$.}
    \label{tab:dispersive_smooth_bed}
\end{table}

\noindent
Table~\ref{tab:dispersive_smooth_bed} reports the CPU times for semi-implicit and explicit schemes under increasing $\lambda$. Even with smooth bathymetry, the qualitative behavior parallels the flat-bottom case: for $\lambda=100$, the semi-implicit scheme typically outperforms the explicit one due to the larger allowable time step under $\mathrm{CFL}_{\rm IMEX}=2$. As $\lambda$ increases to 500 and 1000, the maximum eigenvalue $\nu_{\max}$ grows, tightening the explicit time-step restriction more severely; hence the explicit scheme’s CPU time rises faster than semi-implicit. The bathymetry variation also affects local wave speeds, but the semi-implicit treatment remains advantageous by handling stiff dispersive terms while accommodating depth-induced speed variability.

\subsection{Solitary waves}

The GN equations   admits exact solitary wave solutions of the form 
\begin{equation}
\begin{split}
h_{\textsf{ex}}  = & h_0(1+\varepsilon \text{sech}^2(\kappa(x-Ct)))\\
u_{\textsf{ex}}= & C \left( 1 -\dfrac{h_0}{h_{\textsf{ex}}}\right)
\end{split}
\end{equation}
with $h_0 $ the depth at rest, $\varepsilon = A/h_0$ the nonlinearity coefficient (amplitude of depth at rest),  $C=\sqrt{gh_0(1+\varepsilon)}$ the soliton
speed, and $\kappa$ a shape parameter given by
$$
\kappa=\sqrt{\dfrac{3\varepsilon}{4 h_0^2(1+\varepsilon)}}
$$

\subsubsection{Accuracy and CPU time}
\label{ssec:accuracy}
To assess the accuracy of the proposed semi-implicit scheme, we consider a solitary wave (soliton) for which an exact solution is available. This allows us to quantify the numerical error and verify the expected convergence rates for both the first- and second-order schemes. The numerical tests are performed for different values of $\lambda$ ($500$, $1000$, $5000$, and $10000$), with results presented in Figures~(a) and (b). The CFL$_{\rm IMEX}$ number is chosen as approximately $2.7$ for most cases, except for $\lambda = 10000$, where it is set to $2.0$. The numerical scheme follows the second version of the implicit discretization, in which $h^2$ is treated implicitly in the fourth equation.
\begin{figure}[!ht]
     \begin{subfigure}[b]{0.49\textwidth}
         \centering
         \includegraphics[trim={2.2cm 0 3cm 0},clip,width=\textwidth]{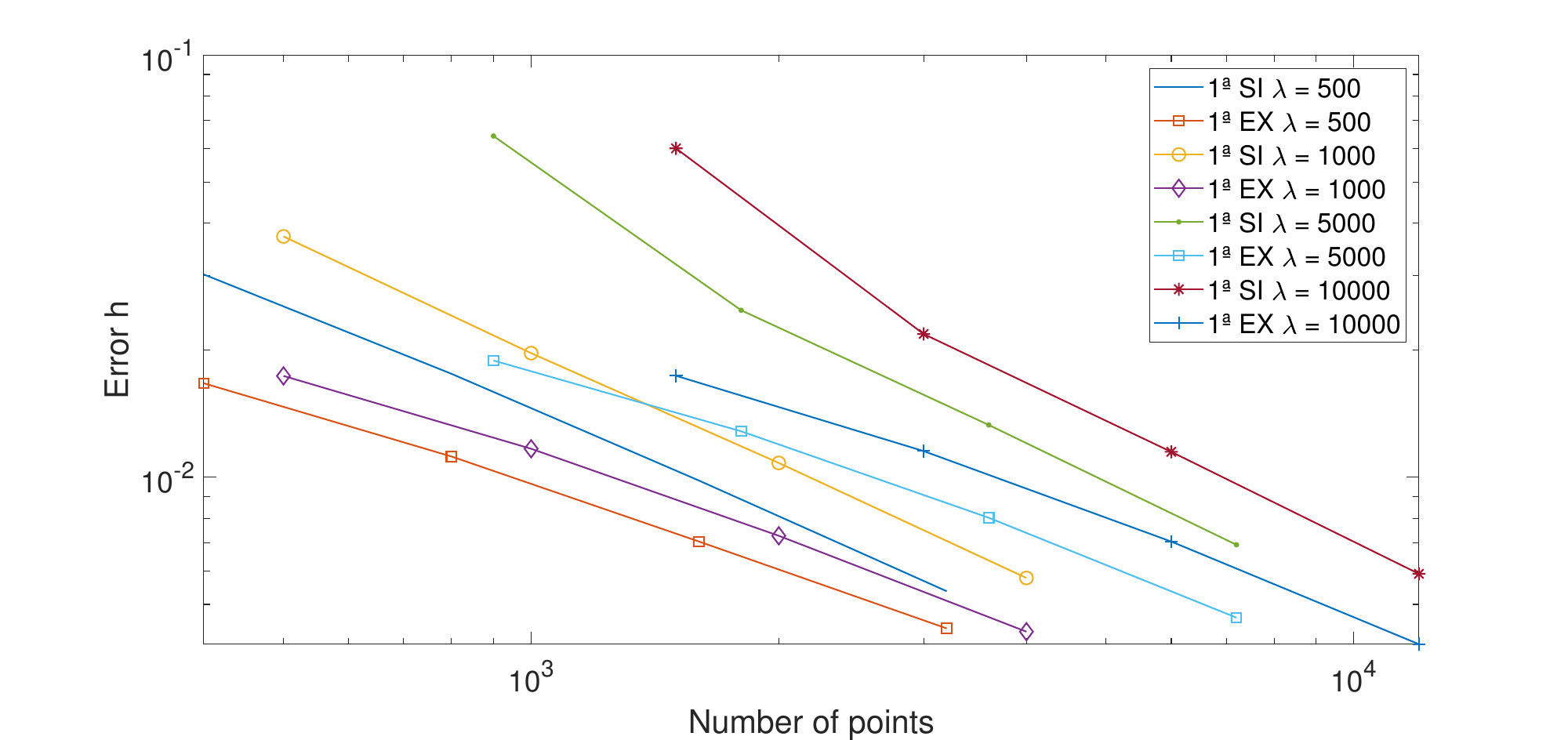}
          \caption{First order grid convergence: for the depth $h$ ($N_{\text{Points}}$-Error, log scale)}
         \label{ssec:sfig:accuracy:h_1}
     \end{subfigure}
     \begin{subfigure}[b]{0.49\textwidth}
         \centering
         \includegraphics[trim={2.2cm 0 3cm 0},clip,width=\textwidth]{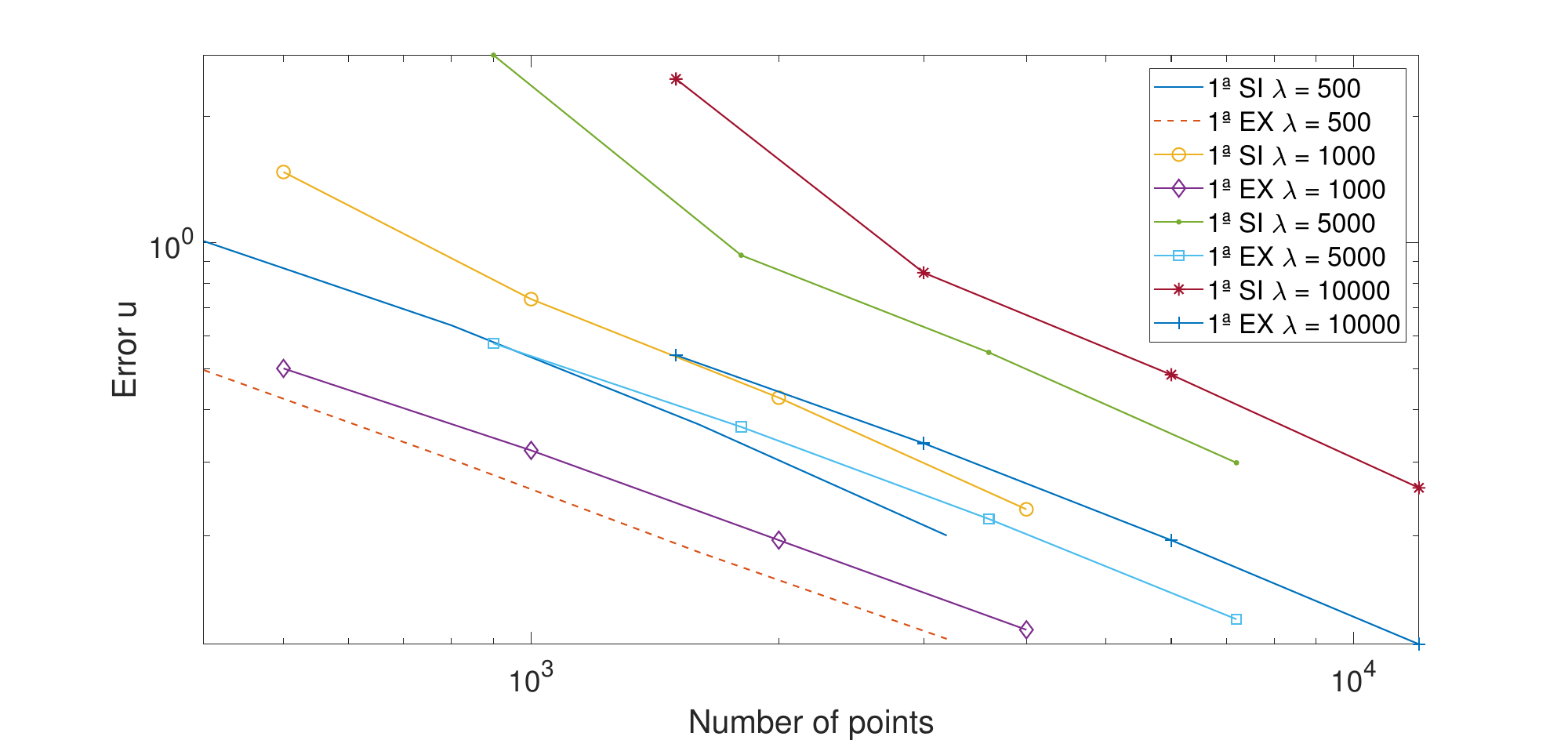}
         \caption{First order grid convergence: for the velocity $u$ ($N_{\text{Points}}$-Error, log scale)}
         \label{ssec:sfig:accuracy:u_1}
     \end{subfigure}
     \\
     \begin{subfigure}[b]{0.49\textwidth}
         \centering
         \includegraphics[trim={2.2cm 0 3cm 0},clip,width=\textwidth]{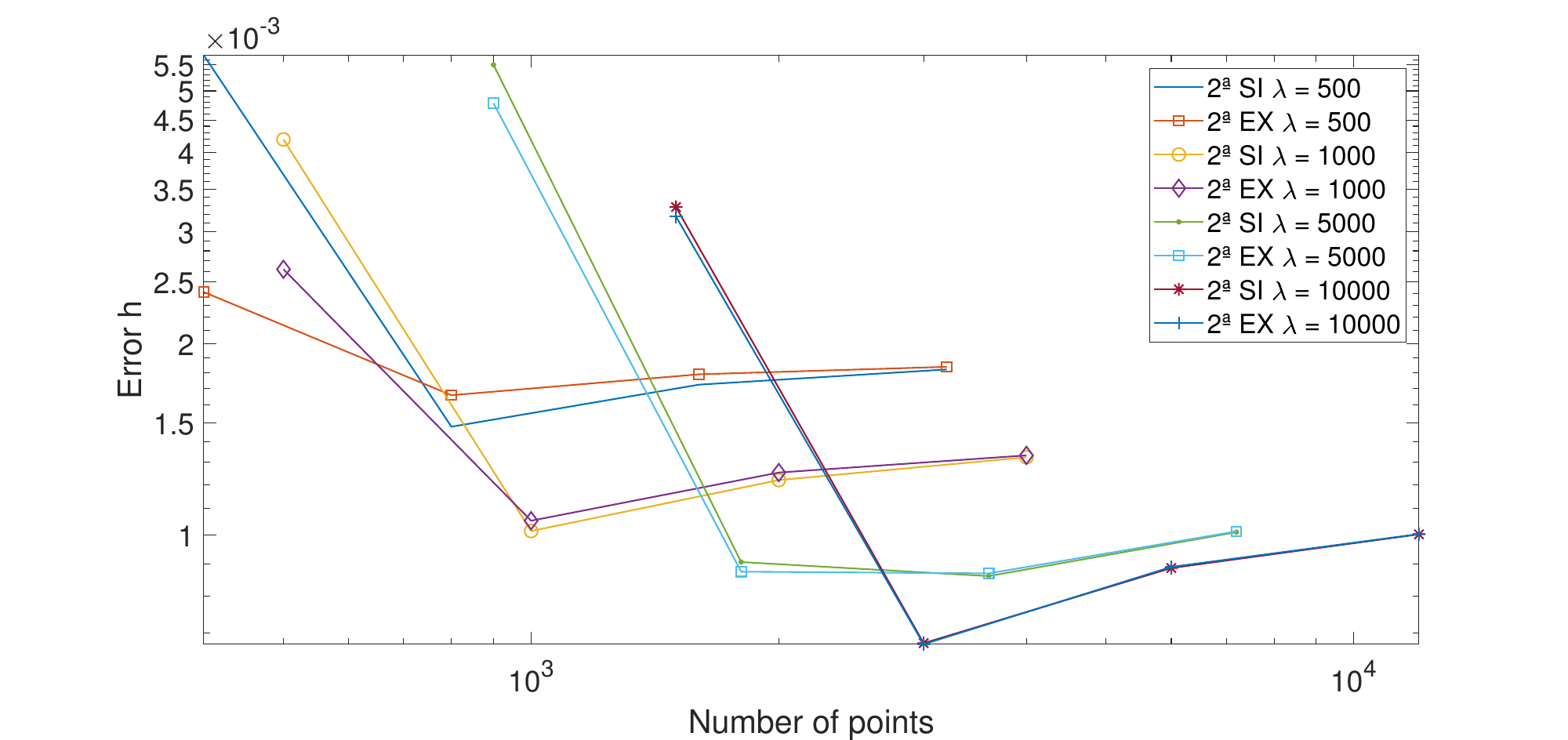}
          \caption{Second order xgrid convergence: for the depth $h$  ($N_{\text{Points}}$-Error, log scale)}
         \label{ssec:sfig:accuracy:h_2}
     \end{subfigure}
     \begin{subfigure}[b]{0.49\textwidth}
         \centering
         \includegraphics[trim={2.2cm 0 3cm 0},clip,width=\textwidth]{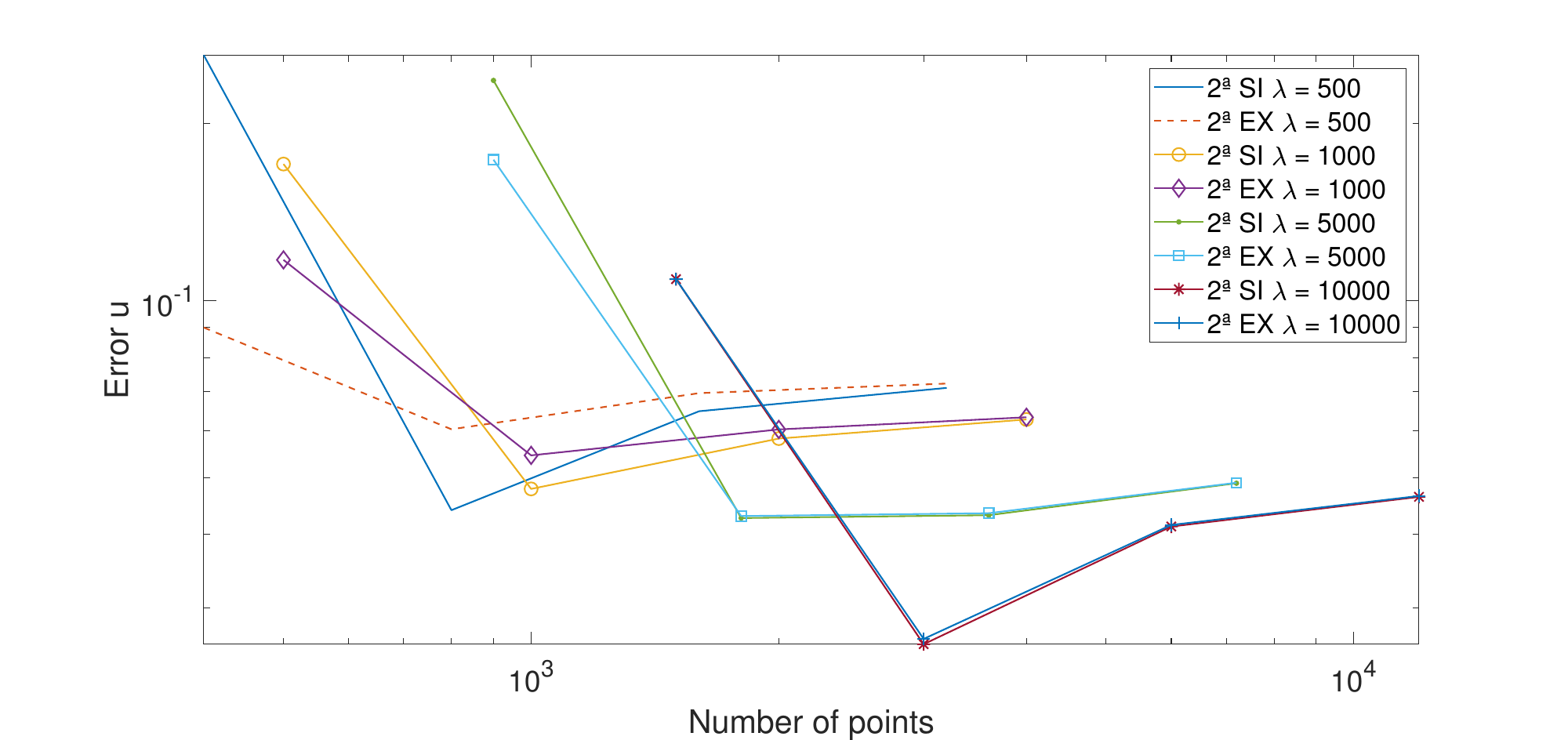}
         \caption{Second order grid convergence: for the velocity $u$ ($N_{\text{Points}}$-Error, log scale)}
         \label{ssec:sfig:accuracy:u_2}
     \end{subfigure}\\
     \begin{center}
     \begin{subfigure}[b]{0.49\textwidth}
         \centering
         \includegraphics[trim={2.2cm 0 3cm 0},clip,width=\textwidth]{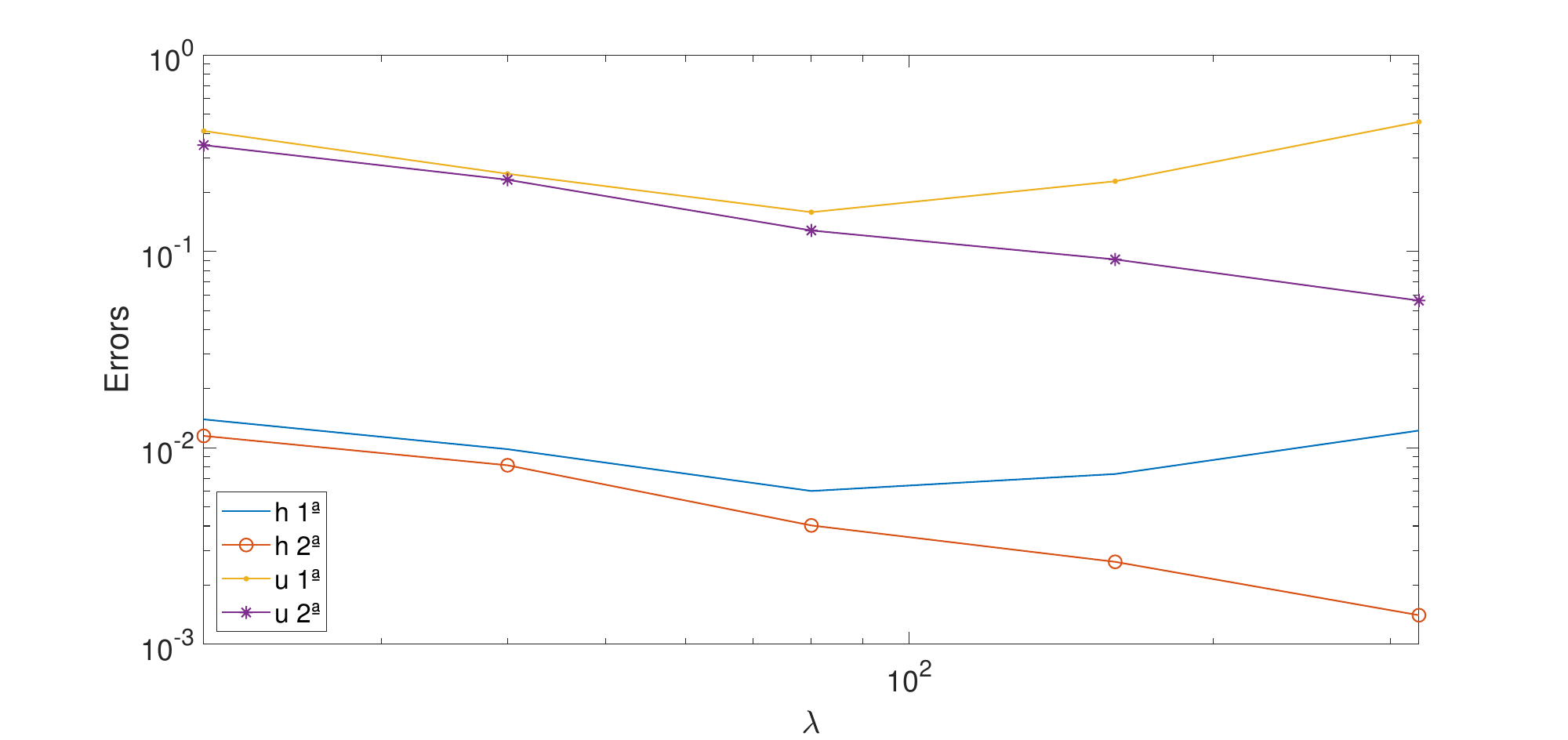}
         \caption{Error convergence vs $\lambda$.}
         \label{ssec:sfig:accuracy:lambda}
     \end{subfigure} 
     \begin{subfigure}[b]{0.49\textwidth}
         \centering
         \includegraphics[trim={2.2cm 0 3cm 0},clip,width=\textwidth]{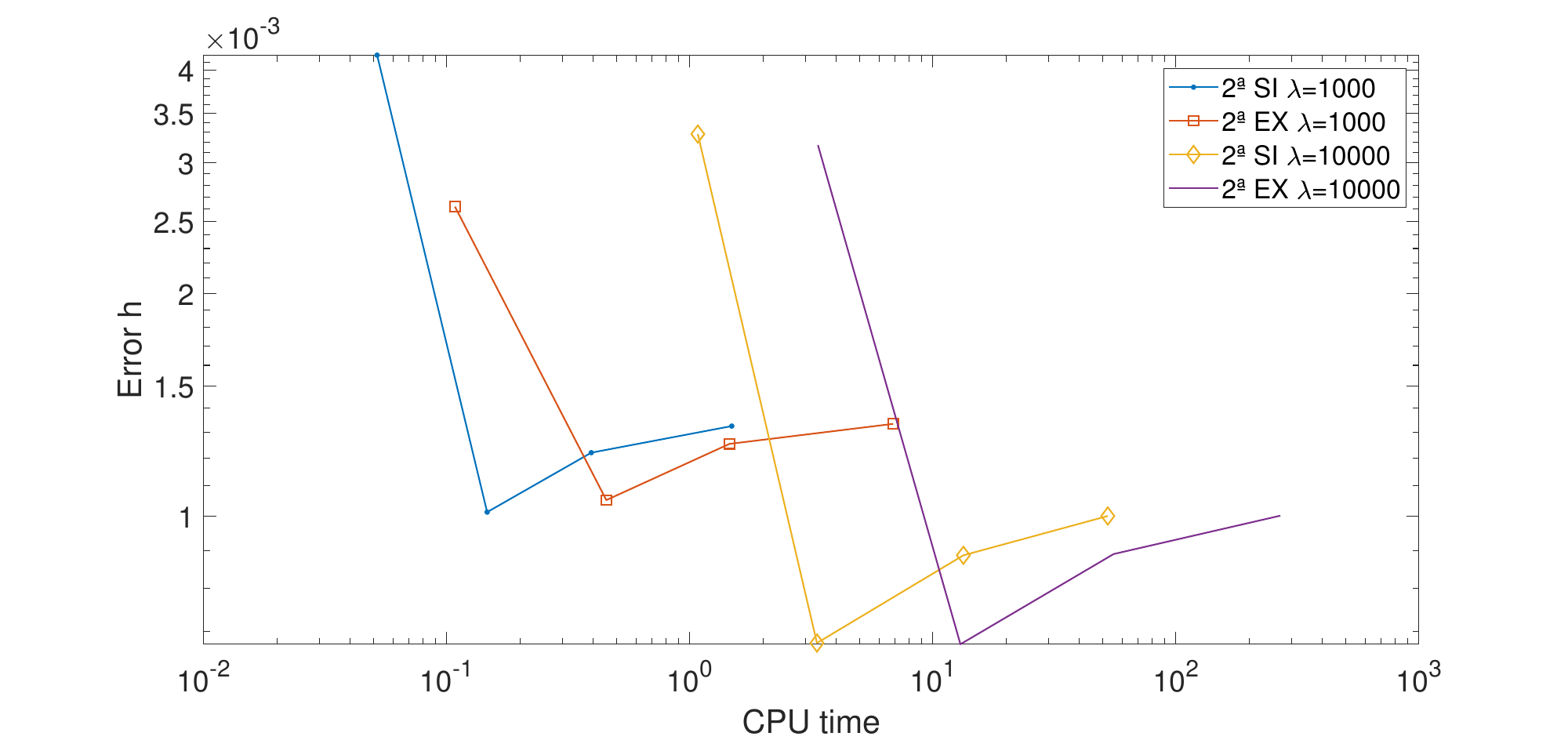}
         \caption{Error vs CPU time expressed in second}
         \label{ssec:sfig:accuracy:CPU}
     \end{subfigure}
     \end{center}
     \caption{Test \ref{ssec:accuracy} L$^1$ errors vs number of points, L$^1$ errors vs $\lambda$ and L$^1$ errors vs CPU time expressed in second. Numerical accuracy analysis: \ref{ssec:sfig:accuracy:h_1}-\ref{ssec:sfig:accuracy:h_2} error comparison for first- and second-order schemes for $h$; \ref{ssec:sfig:accuracy:u_1}-\ref{ssec:sfig:accuracy:u_2} error comparison for first- and second-order schemes for $u$; while \ref{ssec:sfig:accuracy:lambda} error reduction for increasing $\lambda$ and \ref{ssec:sfig:accuracy:CPU} error vs CPU time expressed in second. The final time is $1$.}
     \label{ssec:fig:accuracy}
\end{figure}

The initial conditions for the solitary wave are given by
\begin{align}
    h_{\infty} &= 1, \quad A = 1, \quad \varepsilon = \frac{A}{h_{\infty}}, \\
    k &= \sqrt{\frac{3\varepsilon}{4 h_{\infty}^2 (1+\varepsilon)}}, \quad
    c = \sqrt{g h_{\infty} (1+\varepsilon)},
\end{align}
where $g = 9.81$ is the gravitational acceleration. The initial profiles for water height and velocity are defined as
\begin{align}
    h(x,0) &= h_{\infty} \left( 1 + \varepsilon \operatorname{sech}^2 (k x) \right), \\
    u(x,0) &= c \left( 1 - \frac{h_{\infty}}{h(x,0)} \right),
\end{align}
meanwhile $\eta(x,0) = h(x,0)$ and $w = -h(x,0)u_x(x,0).$

The numerical results, shown in Figures~\ref{ssec:sfig:accuracy:h_1}-\ref{ssec:sfig:accuracy:u_2}, confirm the expected accuracy of the scheme. The first-order method exhibits an error scaling consistent with $\mathcal{O}(\Delta x)$, while the second-order scheme achieves an error proportional to $\mathcal{O}(\Delta x^2)$. These behaviors remain consistent across different values of $\lambda$, demonstrating the robustness of the numerical approach. Additionally, Figure~\ref{ssec:sfig:accuracy:lambda} illustrates the error reduction as $\lambda$ increases, and Figure~\ref{ssec:sfig:accuracy:CPU} the computational cost comparison further validating the semi-implicit formulation. These results highlight the accuracy and efficiency of the proposed scheme for capturing solitary wave dynamics with high fidelity. The exact solution is given by:
\begin{align}
    h(x,t) &= h_{\infty} \left( 1 + \varepsilon \operatorname{sech}^2 (k x-ct) \right), \\
    u(x,t) &= c \left( 1 - \frac{h_{\infty}}{h(x,t)} \right).    
\end{align}

\begin{table}[!ht]
    \numerikNine
    \centering
    \begin{tabular}{|l|cc|}
        \hline \multicolumn{3}{|c|}{\textbf{hSGN - CPU time}} \\
        \hline
        Scheme & Order 1 & Order 2 \\
        \hline
        IMEX  & 1.2502  & 4.0350  \\
        Explicit  & 3.2105  & 16.9138  \\
        \hline
    \end{tabular}
    \caption{Computational cost (in seconds) of first- and second-order semi-implicit and explicit schemes for $\lambda = 1000$ and with $N=2000$. The CFL$_{\rm IMEX} = 2.5$ while the CFL$_{\rm Exp} = 0.4$.}
    \label{ssec:tab:accuracy:cpu}
\end{table}

The computational efficiency of the semi-implicit (IMEX) schemes is evident when compared to their explicit counterparts. While explicit methods require a significantly smaller CFL number (0.4) to ensure stability, the IMEX approach remains stable even with a much larger CFL value (2.5). This results in a dramatic reduction in computational cost, as seen in Table~\ref{ssec:tab:accuracy:cpu}. 

In particular, the first-order semi-implicit method is approximately 2.5 times faster than its explicit counterpart, while the second-order one exhibits a 4-fold reduction in computational sense. These results highlight not only the superior stability properties of the semi-implicit approach but also its remarkable efficiency, making it a compelling choice for large-scale simulations.

\subsection{Favre waves}
\label{ssec:Favre_waves}
We consider a smooth hydraulic jump problem on the interval $[x_a,x_b]=[-50,50]$ under constant gravity $g=9.81$. The still‐water depth is set to $h_0=0.2$, and we perform three experiments with jump amplitudes $\epsilon\in\{1.1,\,1.2,\,1.3\}$.
\begin{figure}[!ht]
     \centering
     \begin{subfigure}[b]{0.473\textwidth}
         \centering         \includegraphics[width=\textwidth]{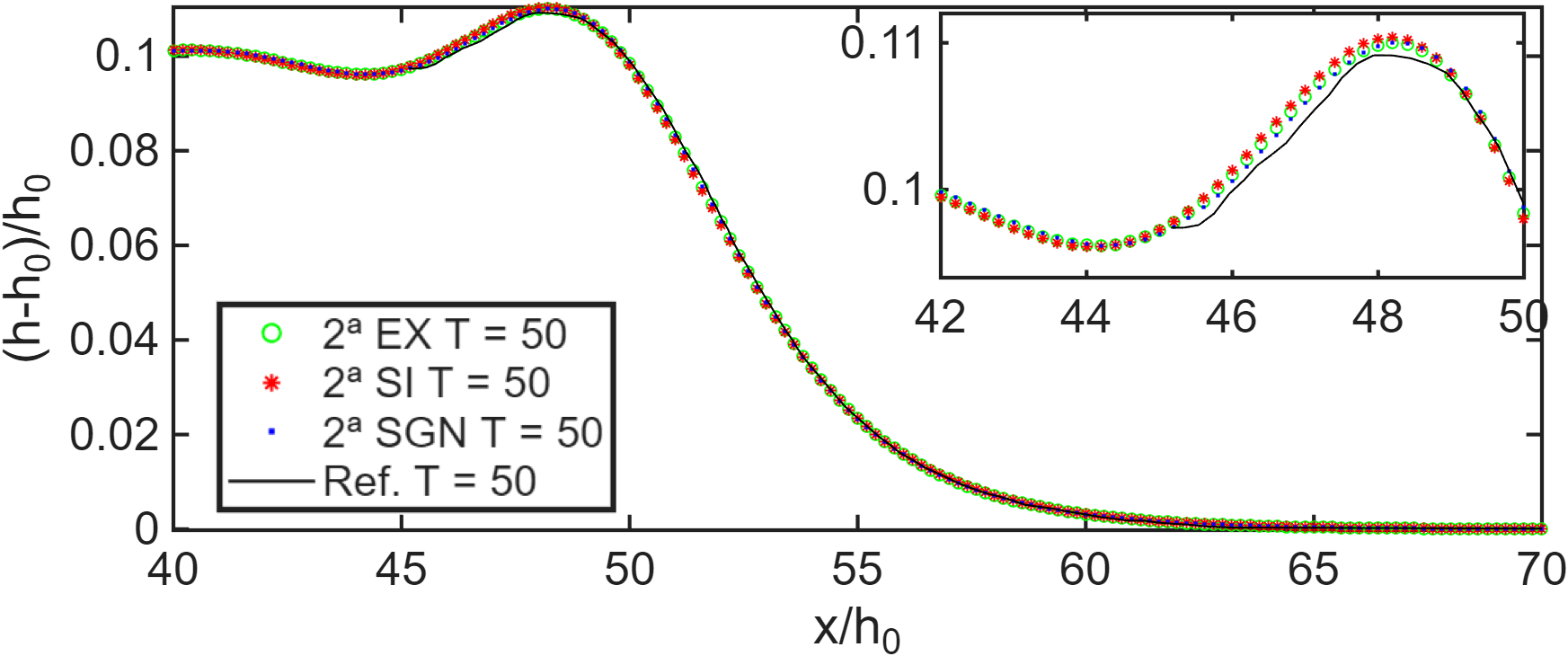}
         \caption{$T = 50,$ $\epsilon = 1.1$ and $\lambda = 100.$}
         \label{sfig:Favre_1_100_1}
     \end{subfigure} 
     \hfill
     \begin{subfigure}[b]{0.473\textwidth}
         \centering         \includegraphics[width=\textwidth]{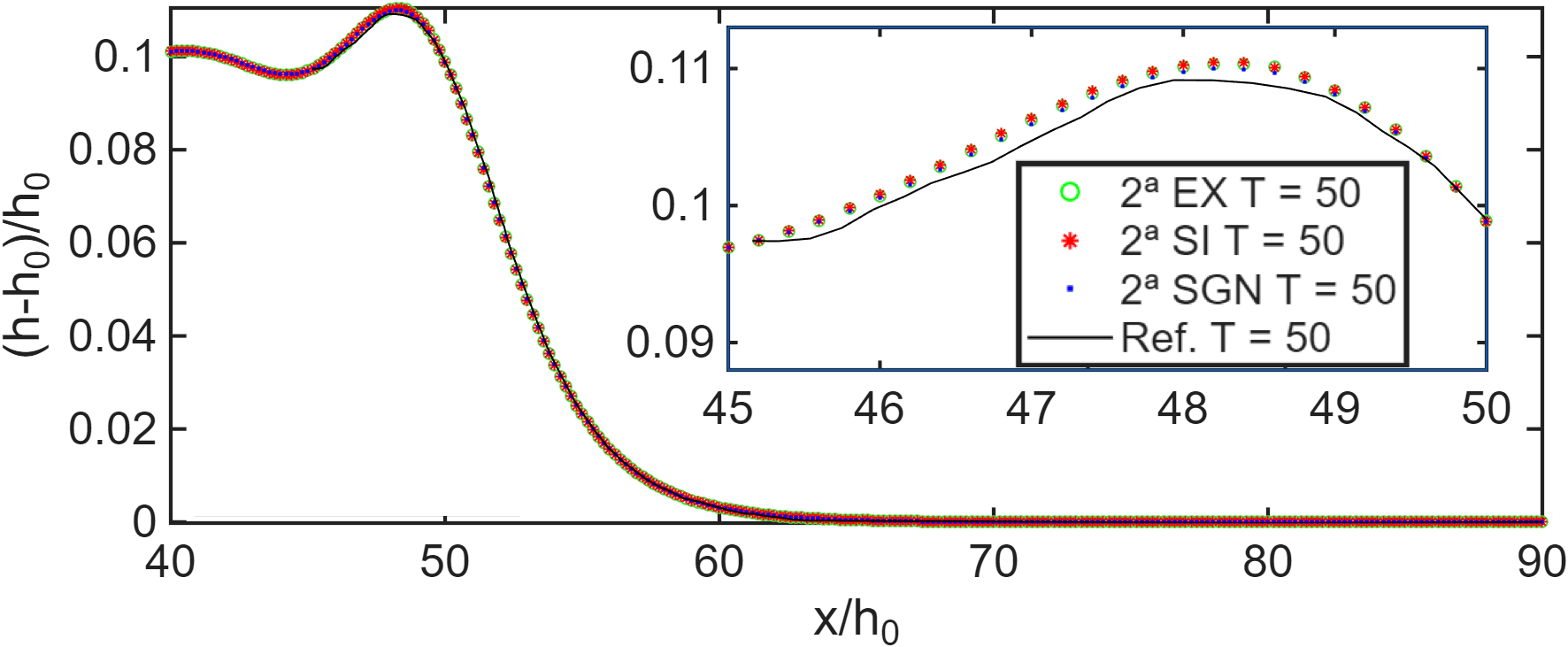}
         \caption{$T = 50,$ $\epsilon = 1.1$ and $\lambda = 500.$}
         \label{sfig:Favre_1_500_1}
     \end{subfigure} 
     \vfill
     \begin{subfigure}[b]{0.473\textwidth}
         \centering         \includegraphics[width=\textwidth]{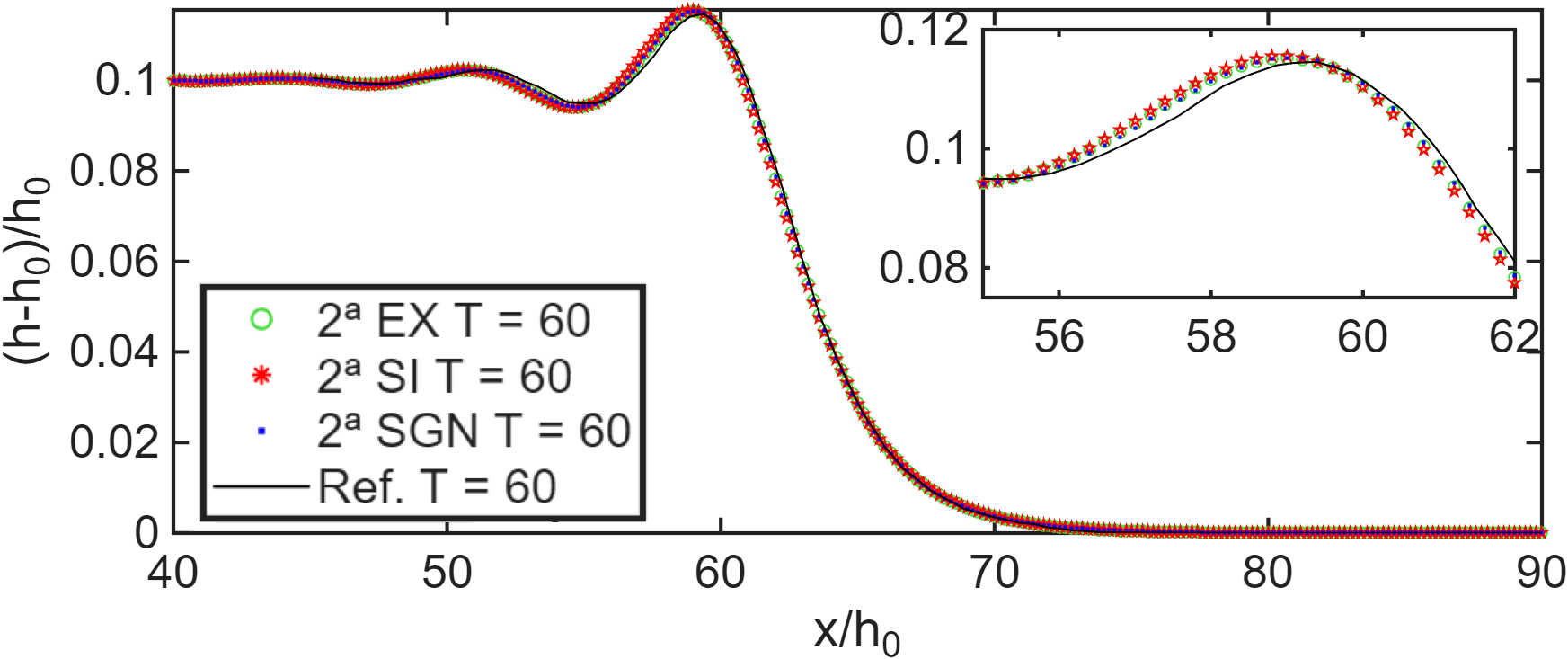}
         \caption{$T = 60,$ $\epsilon = 1.1$ and $\lambda = 100.$}
         \label{sfig:Favre_1_100_2}
     \end{subfigure} 
     \hfill
     \begin{subfigure}[b]{0.473\textwidth}
         \centering         \includegraphics[width=\textwidth]{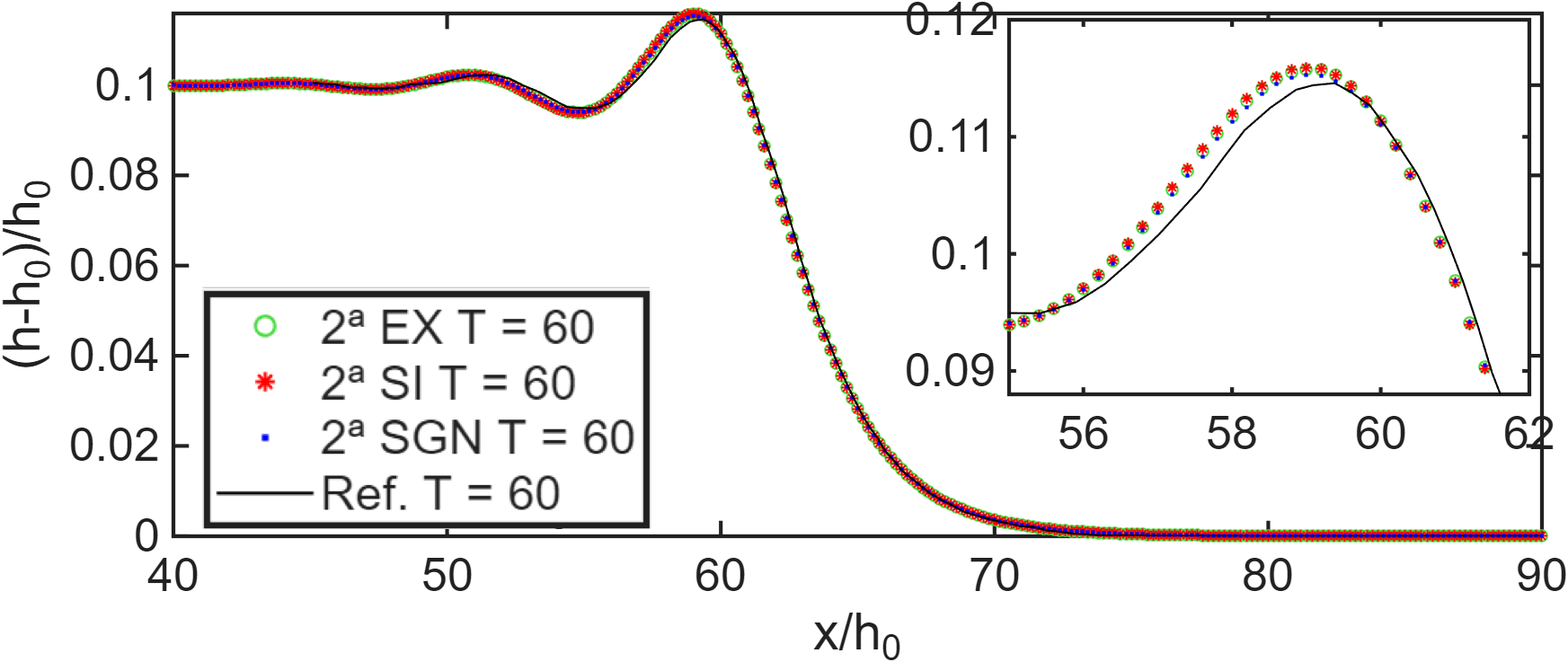}
         \caption{$T = 60,$ $\epsilon = 1.1$ and $\lambda = 500.$}
         \label{sfig:Favre_1_500_2}
     \end{subfigure} 
     \vfill\begin{subfigure}[b]{0.473\textwidth}
         \centering         \includegraphics[width=\textwidth]{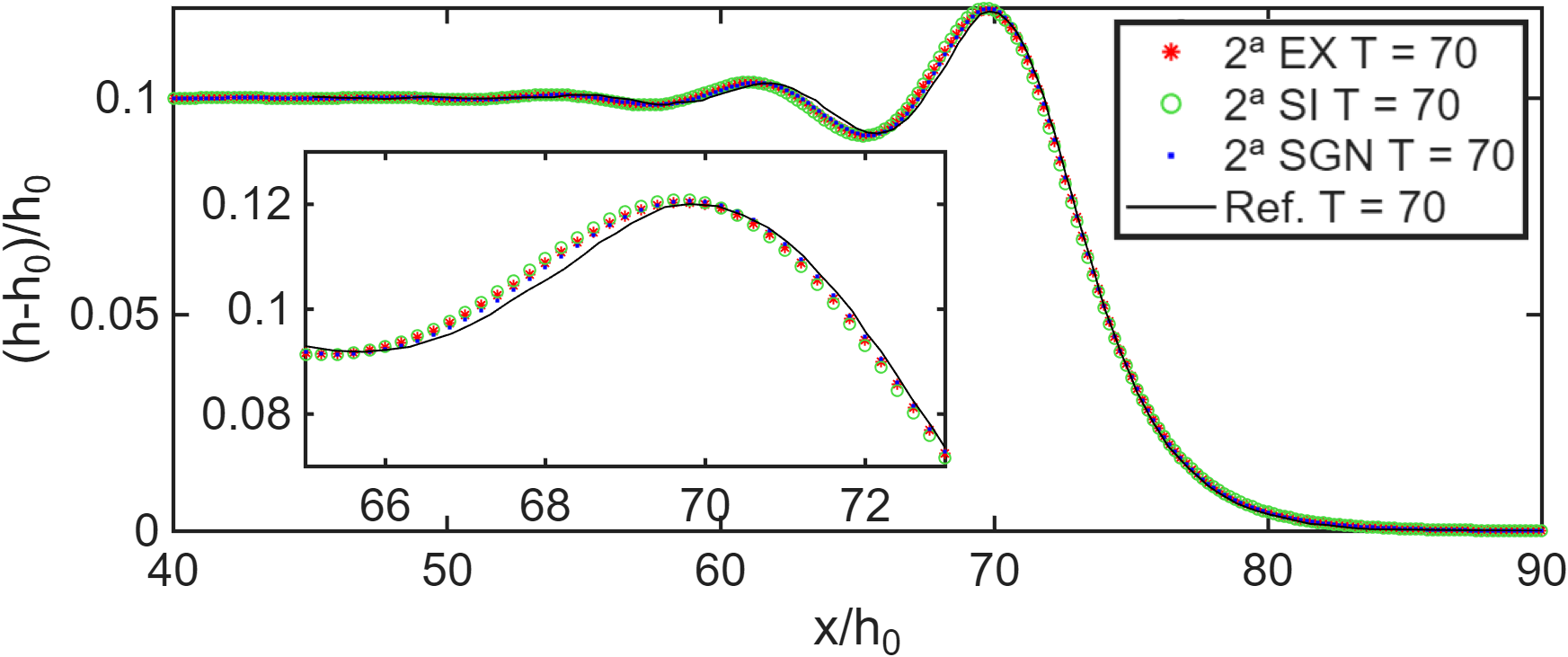}
         \caption{$T = 70,$ $\epsilon = 1.1$ and $\lambda = 100.$}
         \label{sfig:Favre_1_100_3}
     \end{subfigure} 
     \hfill
     \begin{subfigure}[b]{0.473\textwidth}
         \centering         \includegraphics[width=\textwidth]{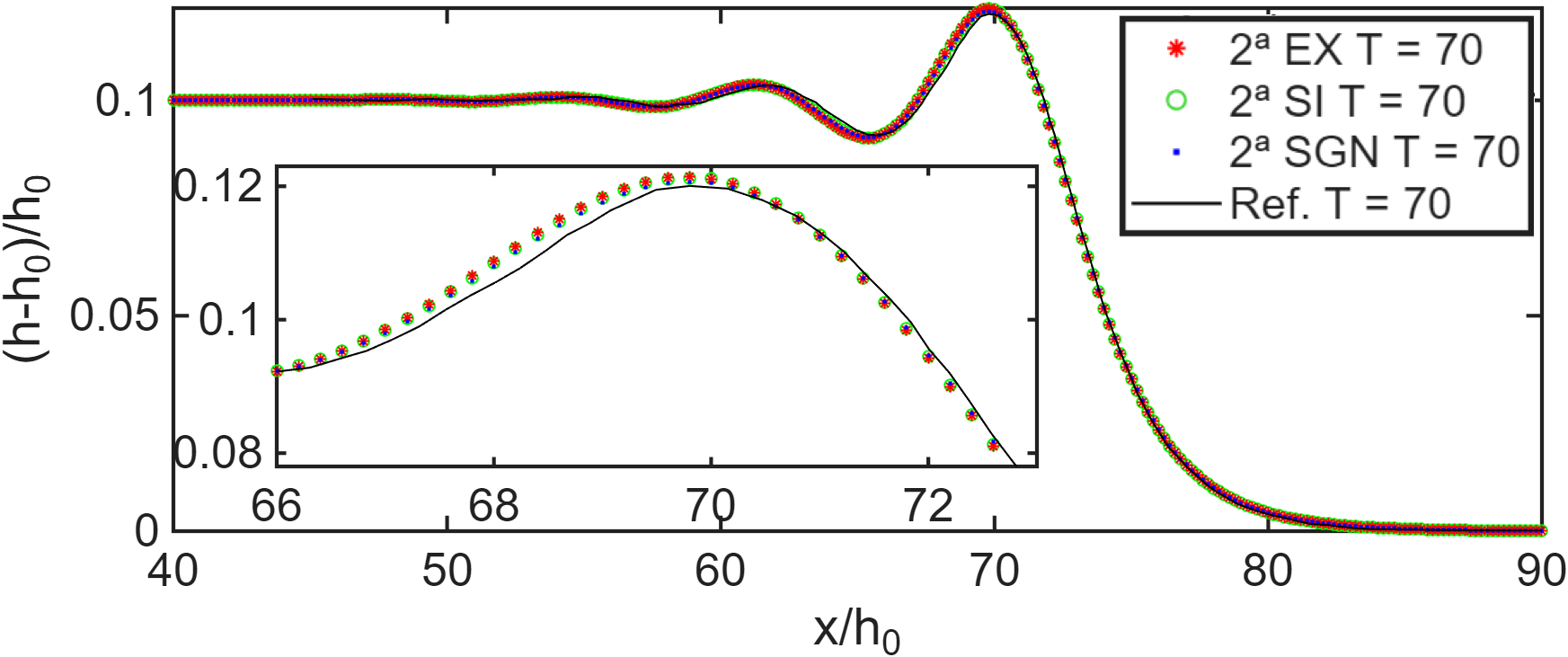}
         \caption{$T = 70,$ $\epsilon = 1.1$ and $\lambda = 500.$}
         \label{sfig:Favre_1_500_3}
     \end{subfigure} 
     \caption{Numerical solutions of Favre waves using second-order semi-implicit and explicit schemes, and second-order scheme for the Serre-Green-Nagdhi model \eqref{eq:SGN_standard}, at  different final times with nonlinearity coefficient $\epsilon=1.1$ and $\lambda \in\{ 100,500\}$. The solutions are plotted with the dimensionless variables $x/h_0$ on the horizontal axis and $(h-h_0)/h_0$ on the vertical axis. Reference: nonlinear potential flow  solution by \cite{wei1995fully}. The CFL numbers are $\mathrm{CFL}_{\rm IMEX}=4$, $\mathrm{CFL}_{\rm EX}=0.4$ and $\mathrm{CFL}_{\rm SGN}=0.9$.}
     \label{fig:favre_eps_1}
\end{figure}
\begin{figure}[!ht]
     \centering
     \begin{subfigure}[b]{0.473\textwidth}
         \centering         \includegraphics[width=\textwidth]{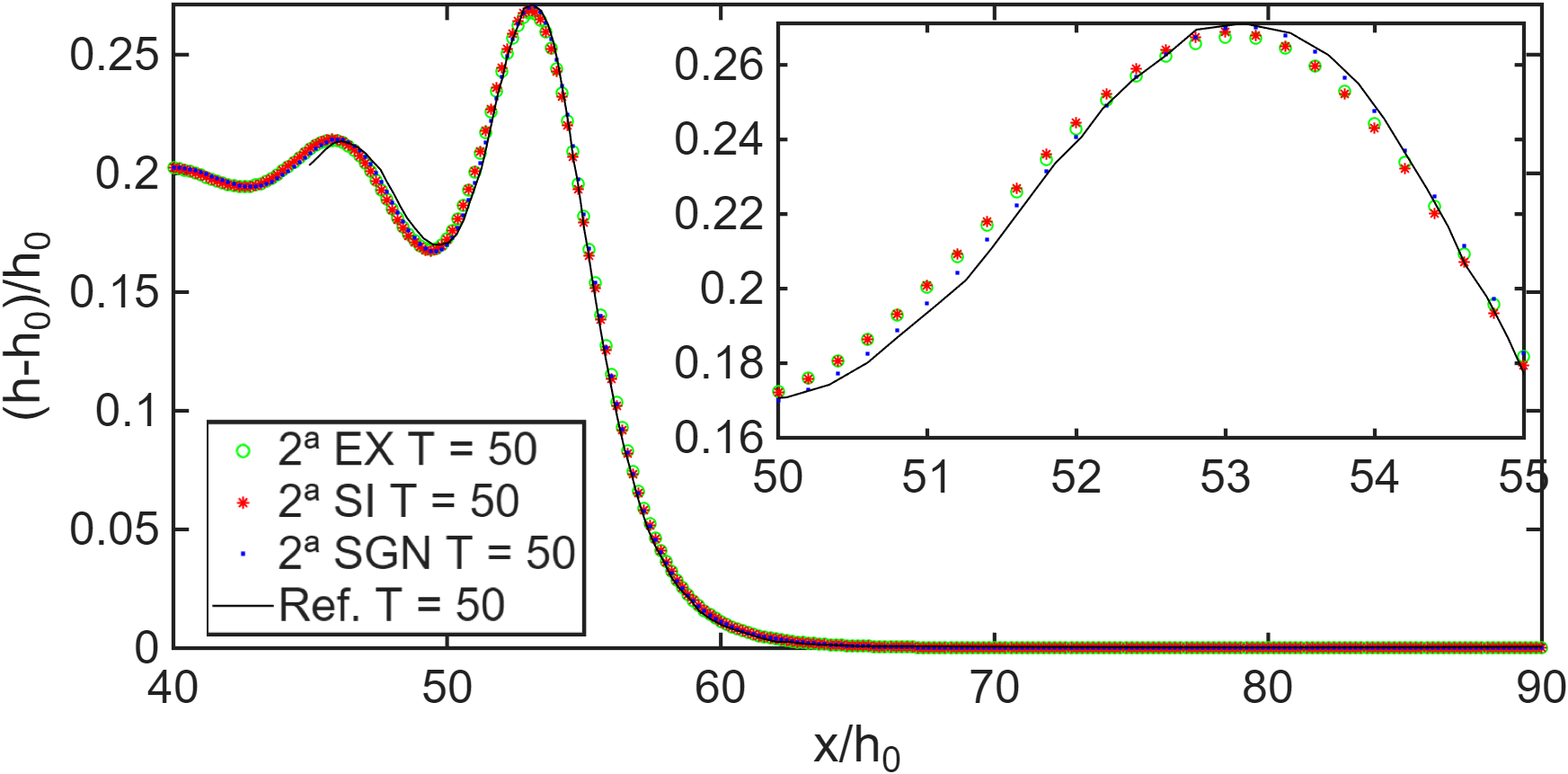}
         \caption{$T = 50,$ $\epsilon = 1.2$ and $\lambda = 100.$}
         \label{sfig:Favre_2_100_1}
     \end{subfigure} 
     \hfill
     \begin{subfigure}[b]{0.473\textwidth}
         \centering         \includegraphics[width=\textwidth]{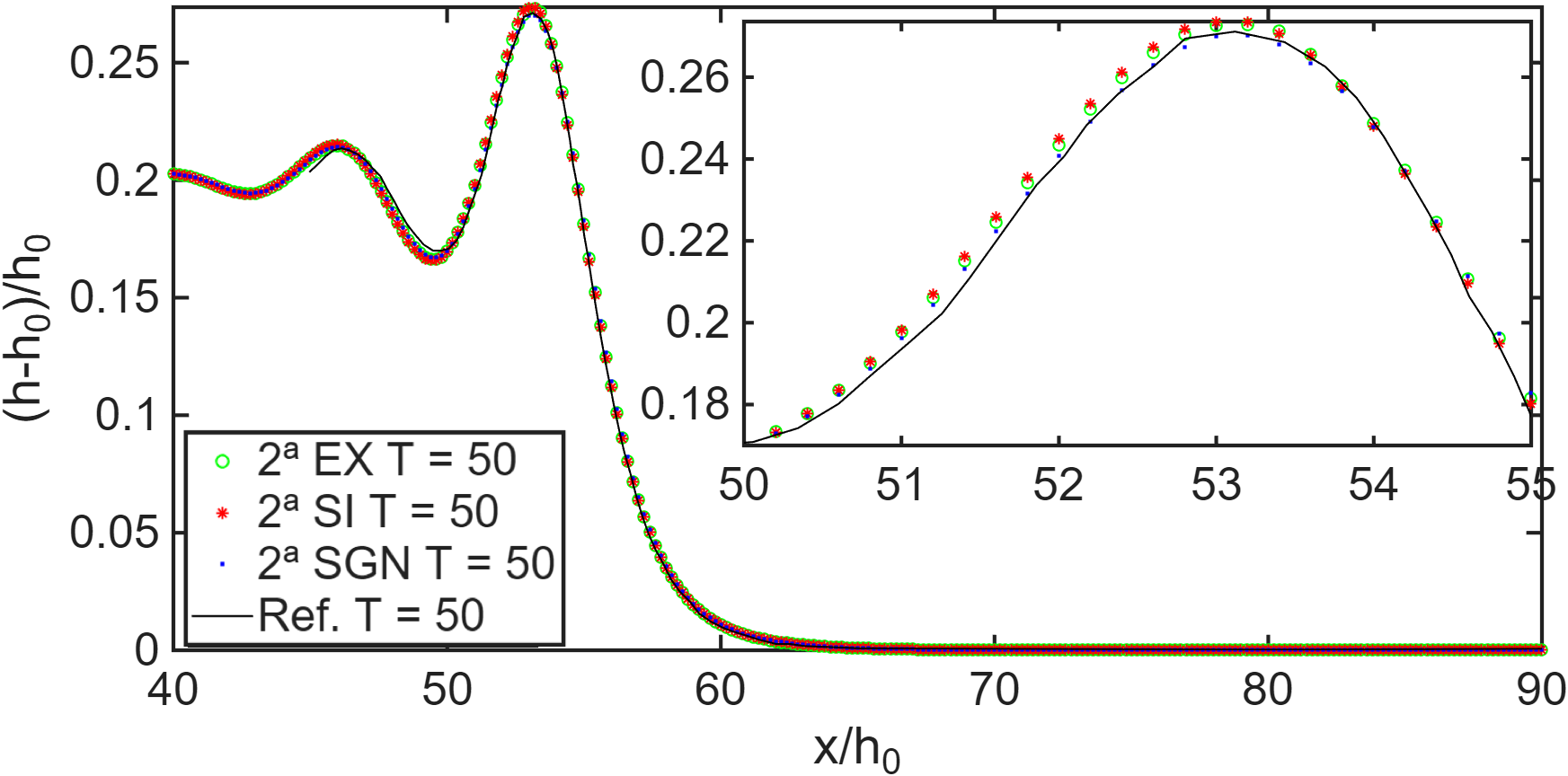}
         \caption{$T = 50,$ $\epsilon = 1.2$ and $\lambda = 500.$}
         \label{sfig:Favre_2_500_1}
     \end{subfigure} 
     \vfill
     \begin{subfigure}[b]{0.473\textwidth}
         \centering         \includegraphics[width=\textwidth]{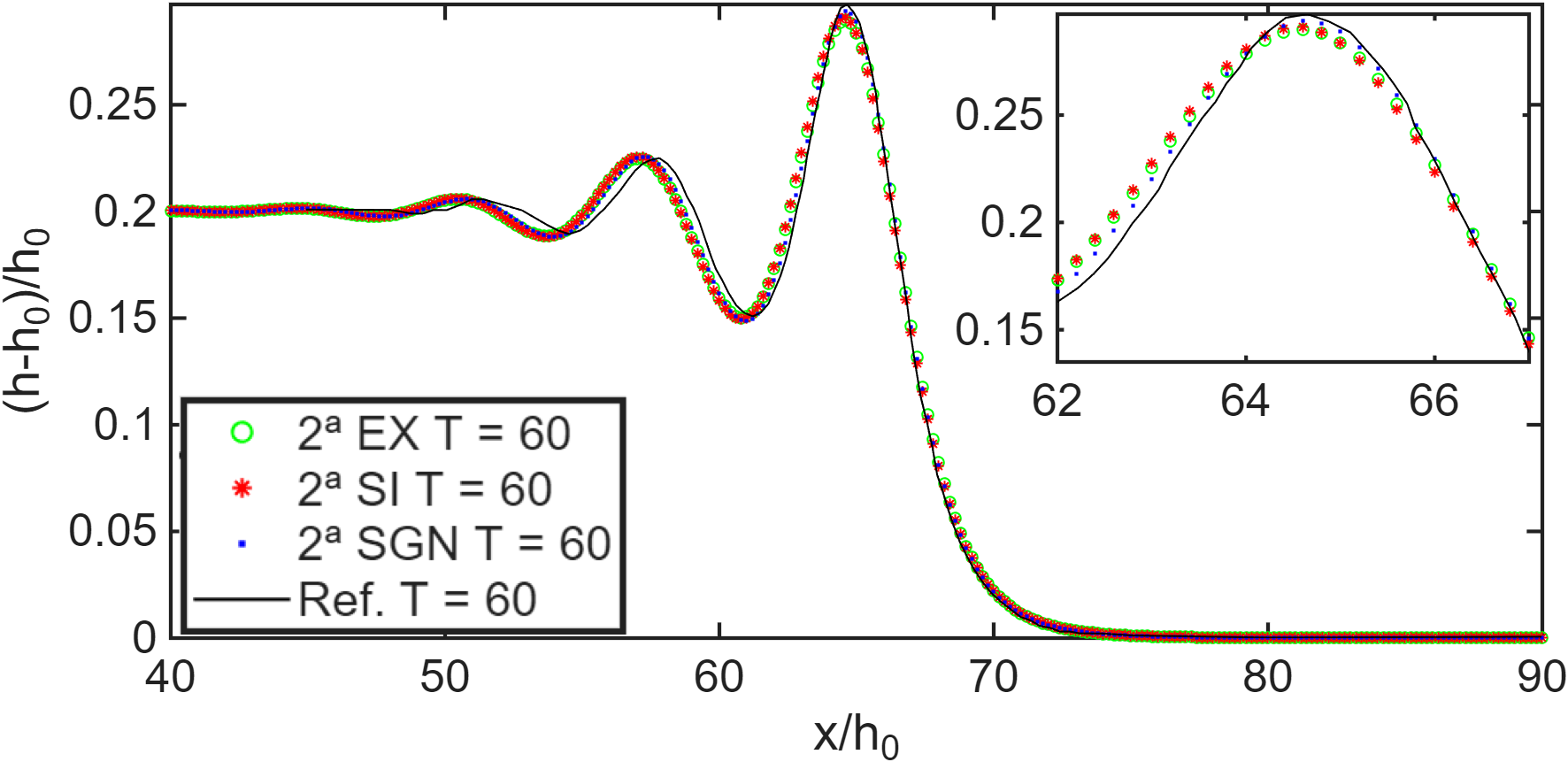}
         \caption{$T = 60,$ $\epsilon = 1.2$ and $\lambda = 100.$}
         \label{sfig:Favre_2_100_2}
     \end{subfigure} 
     \hfill
     \begin{subfigure}[b]{0.473\textwidth}
         \centering         \includegraphics[width=\textwidth]{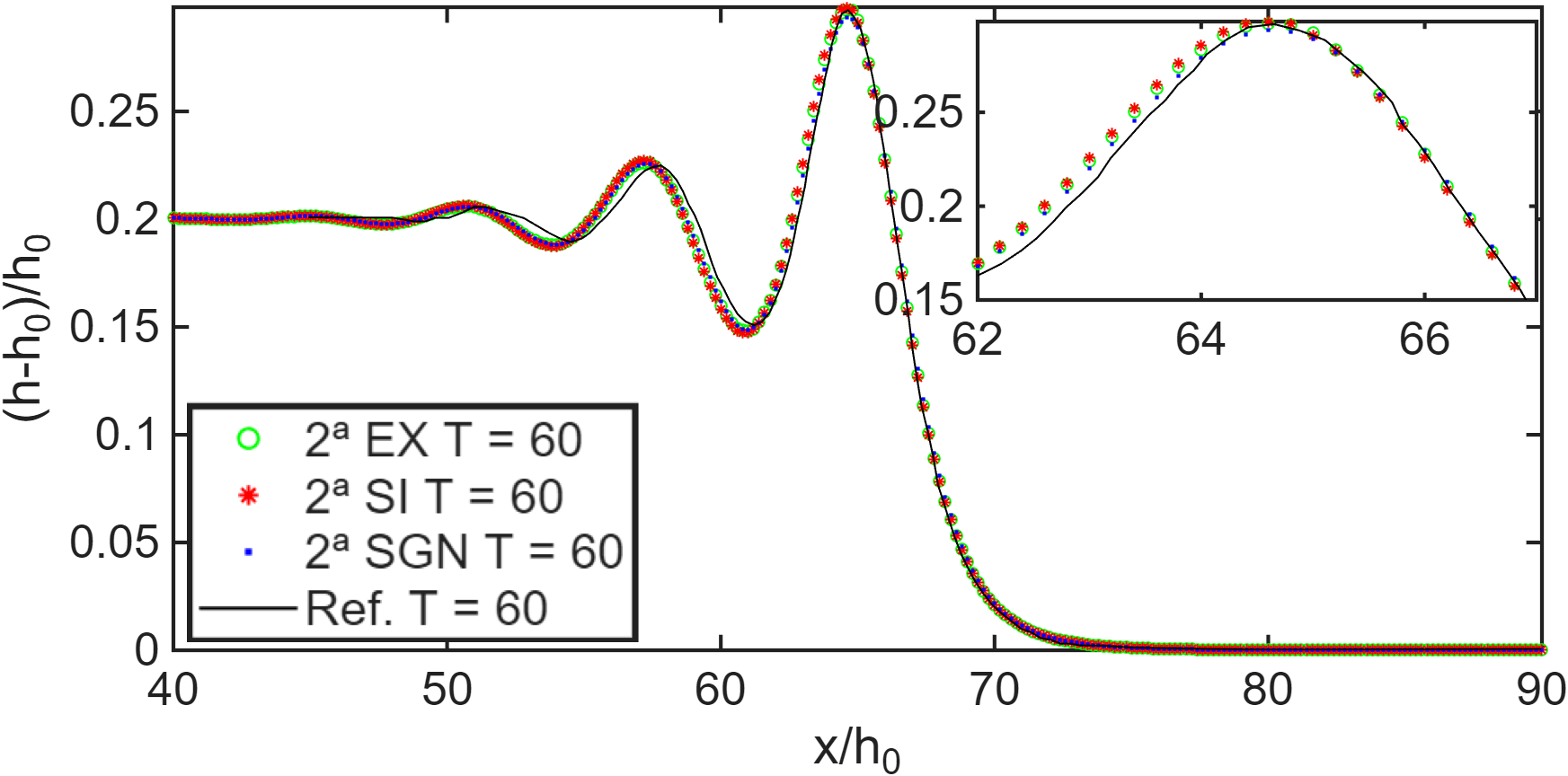}
         \caption{$T = 60,$ $\epsilon = 1.2$ and $\lambda = 500.$}
         \label{sfig:Favre_2_500_2}
     \end{subfigure} 
     \vfill\begin{subfigure}[b]{0.473\textwidth}
         \centering         \includegraphics[width=\textwidth]{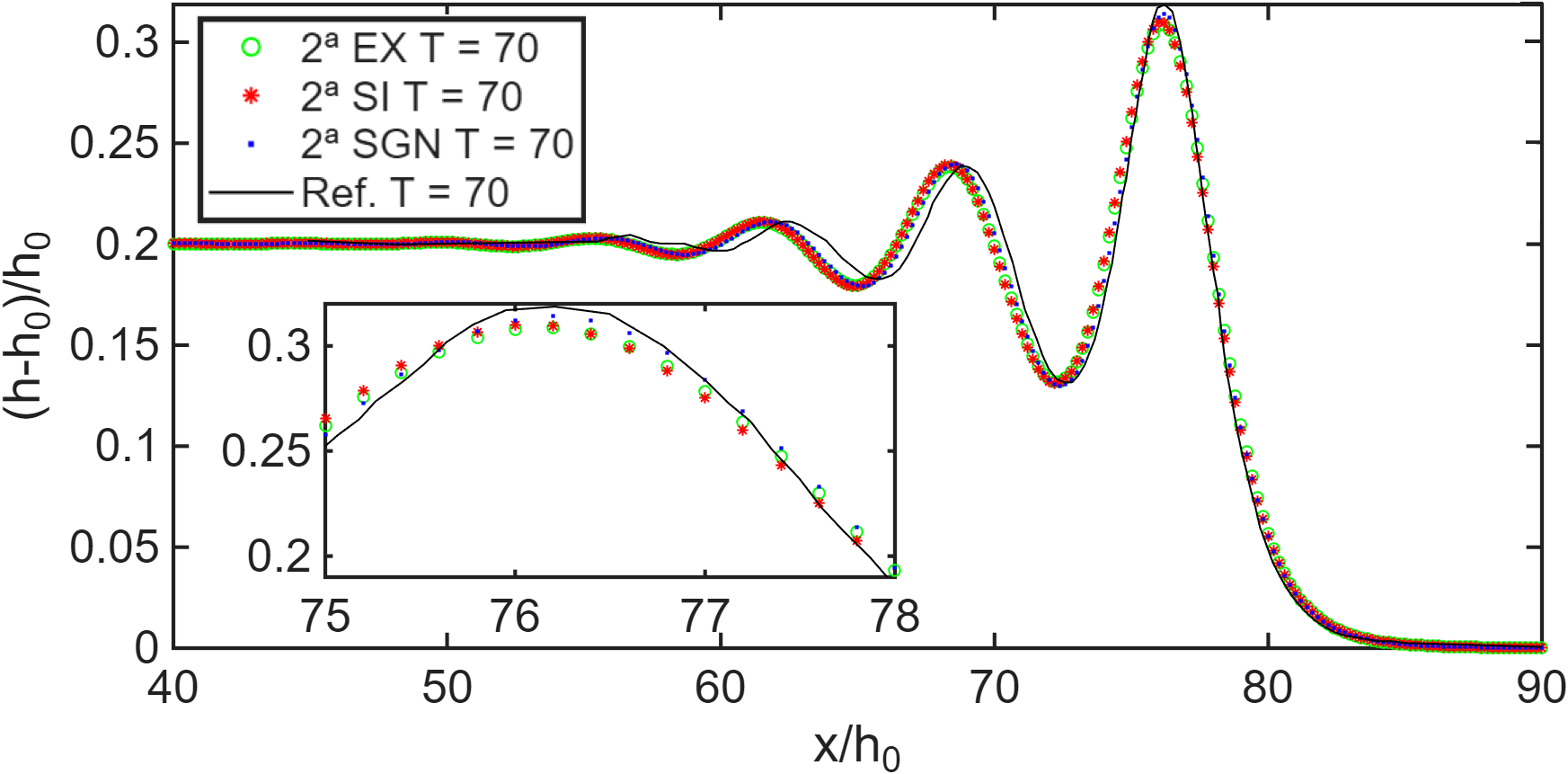}
         \caption{$T = 70,$ $\epsilon = 1.2$ and $\lambda = 100.$}
         \label{sfig:Favre_2_100_3}
     \end{subfigure} 
     \hfill
     \begin{subfigure}[b]{0.473\textwidth}
         \centering         \includegraphics[width=\textwidth]{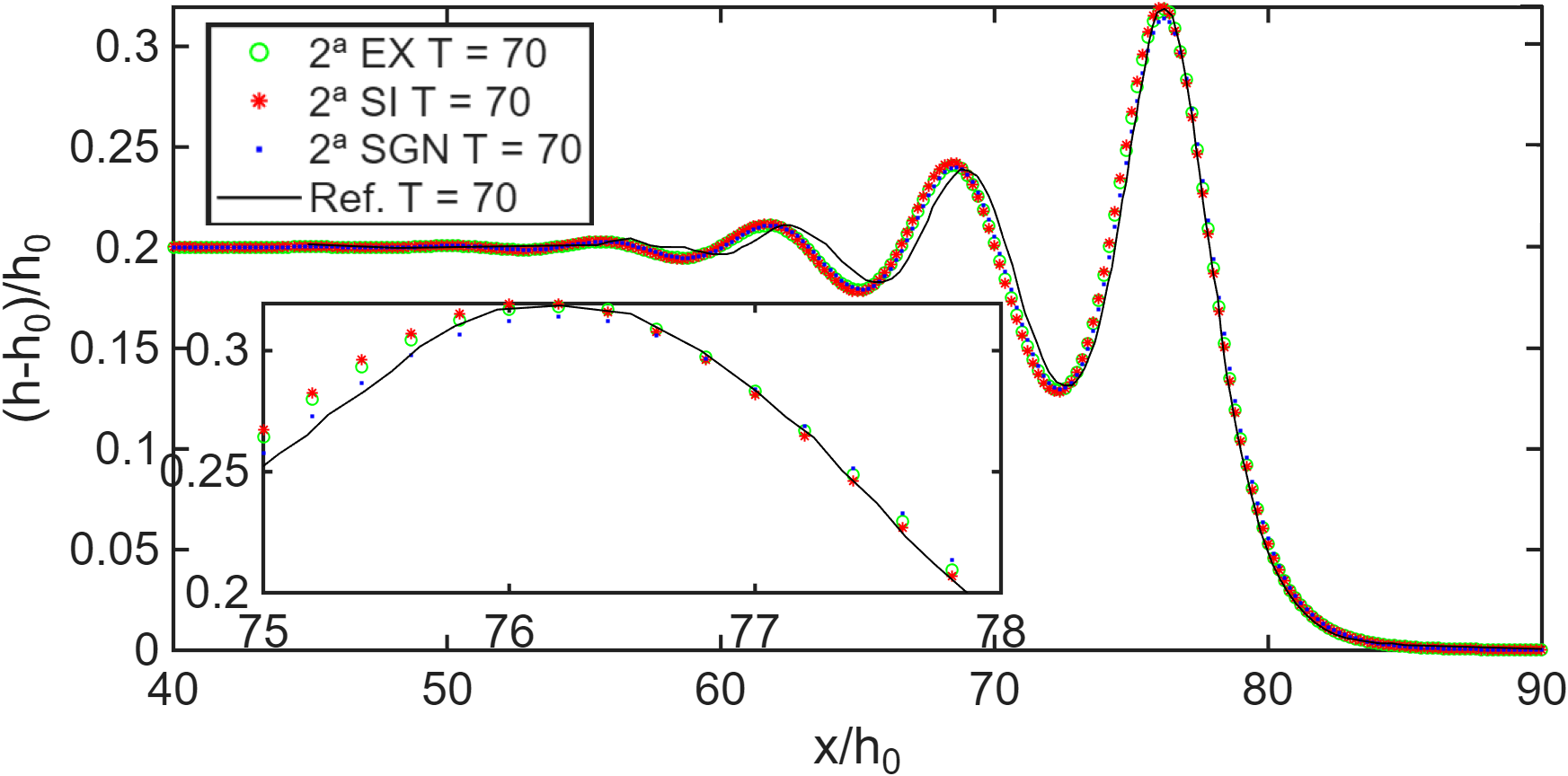}
         \caption{$T = 70,$ $\epsilon = 1.2$ and $\lambda = 500.$}
         \label{sfig:Favre_2_500_3}
     \end{subfigure} 
     \caption{Numerical solutions of Favre waves using second-order semi-implicit and explicit schemes, and second-order scheme for the Serre-Green-Nagdhi model \eqref{eq:SGN_standard}, at  different final times with nonlinearity coefficient $\epsilon=1.1$ and $\lambda \in\{ 100,500\}$. The solutions are plotted with the dimensionless variables $x/h_0$ on the horizontal axis and $(h-h_0)/h_0$ on the vertical axis. Reference: nonlinear potential flow  solution by \cite{wei1995fully}. The CFL numbers are $\mathrm{CFL}_{\rm IMEX}=4$, $\mathrm{CFL}_{\rm EX}=0.4$ and $\mathrm{CFL}_{\rm SGN}=0.9$.}
     \label{fig:favre_eps_2}
\end{figure}
\begin{figure}[!ht]
     \centering
     \begin{subfigure}[b]{0.473\textwidth}
         \centering         \includegraphics[width=\textwidth]{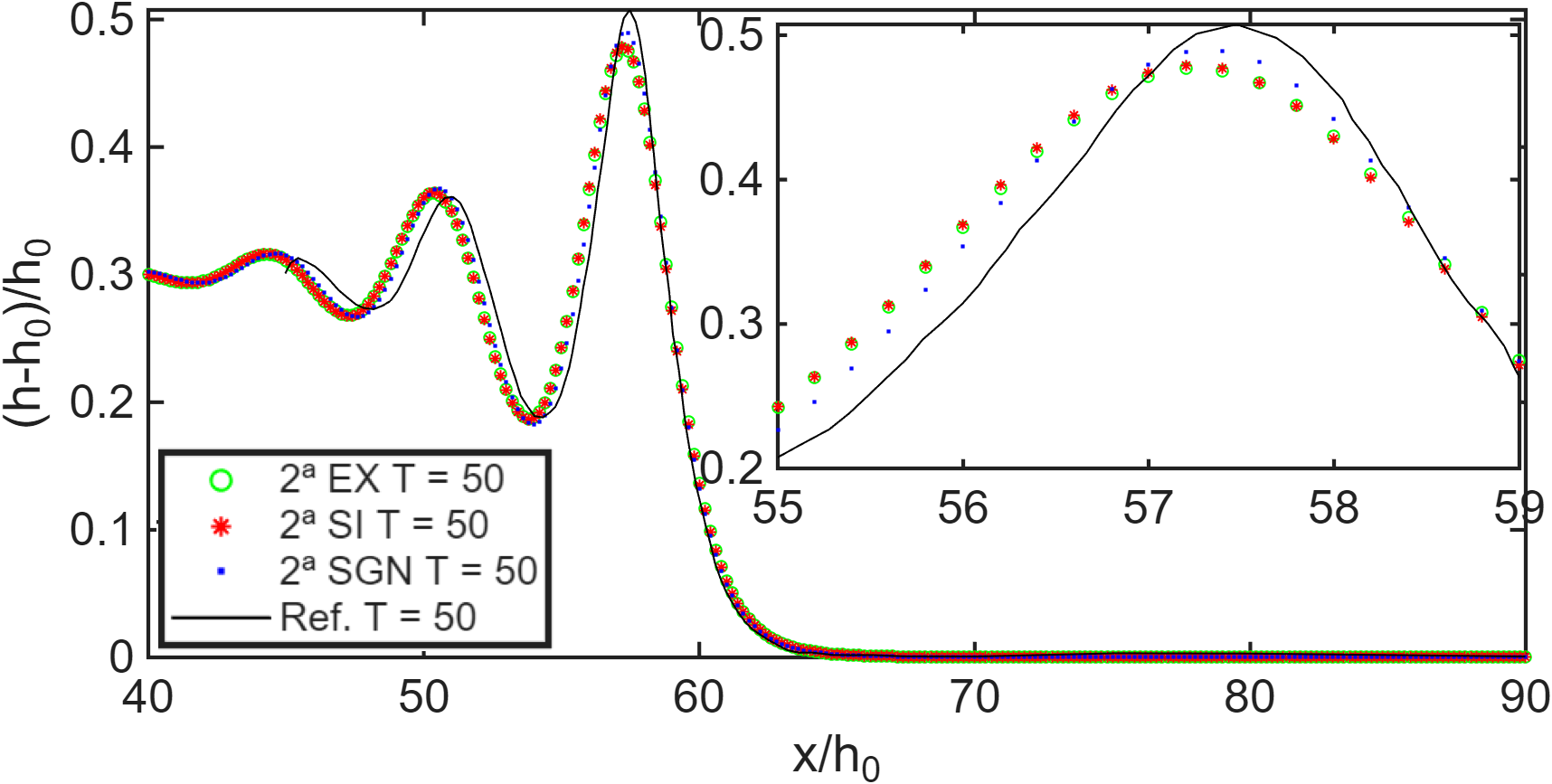}
         \caption{$T = 50,$ $\epsilon = 1.3$ and $\lambda = 100.$}
         \label{sfig:Favre_3_100_1}
     \end{subfigure} 
     \hfill
     \begin{subfigure}[b]{0.473\textwidth}
         \centering         \includegraphics[width=\textwidth]{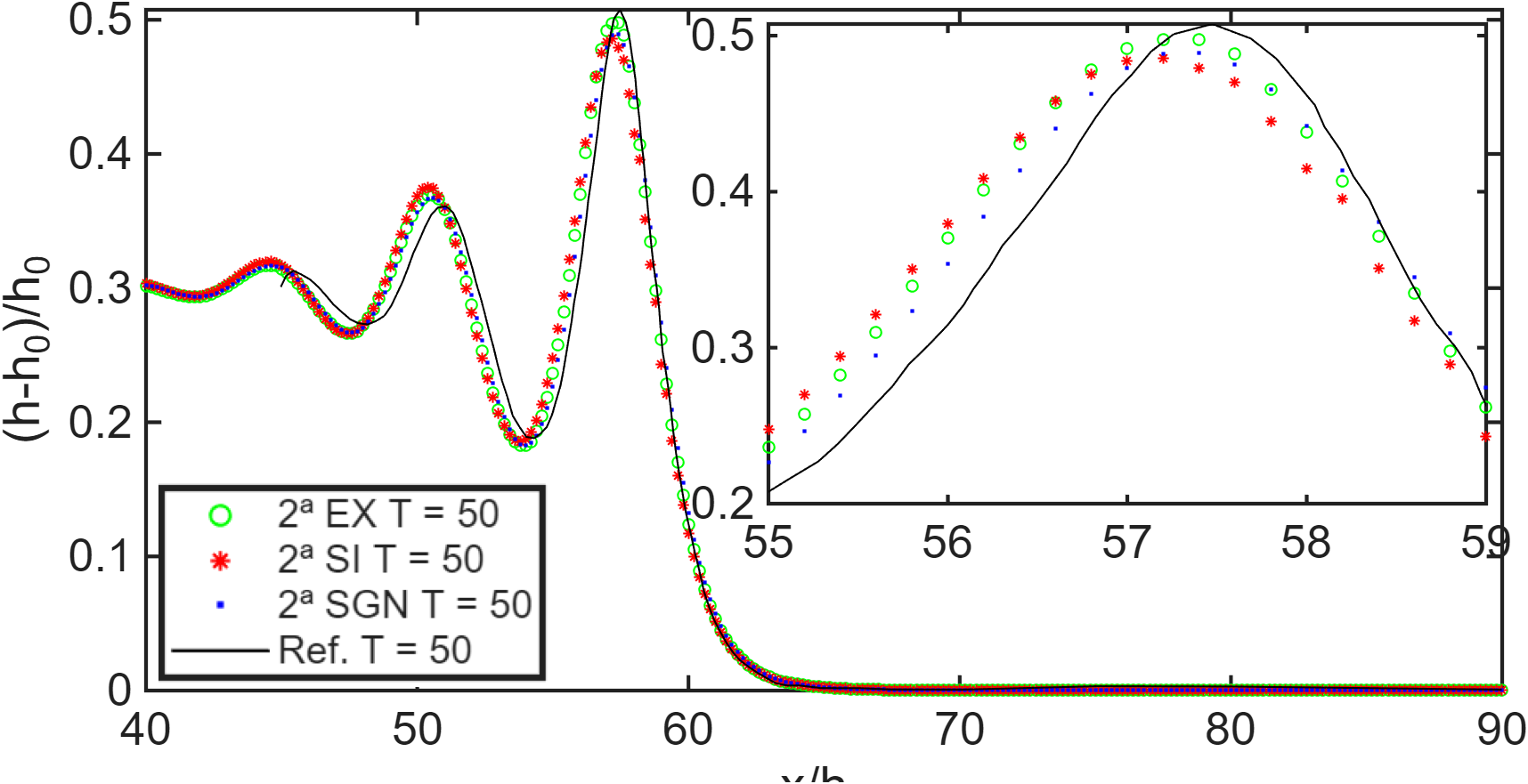}
         \caption{$T = 50,$ $\epsilon = 1.3$ and $\lambda = 500.$}
         \label{sfig:Favre_3_500_1}
     \end{subfigure} 
     \vfill
     \begin{subfigure}[b]{0.473\textwidth}
         \centering         \includegraphics[width=\textwidth]{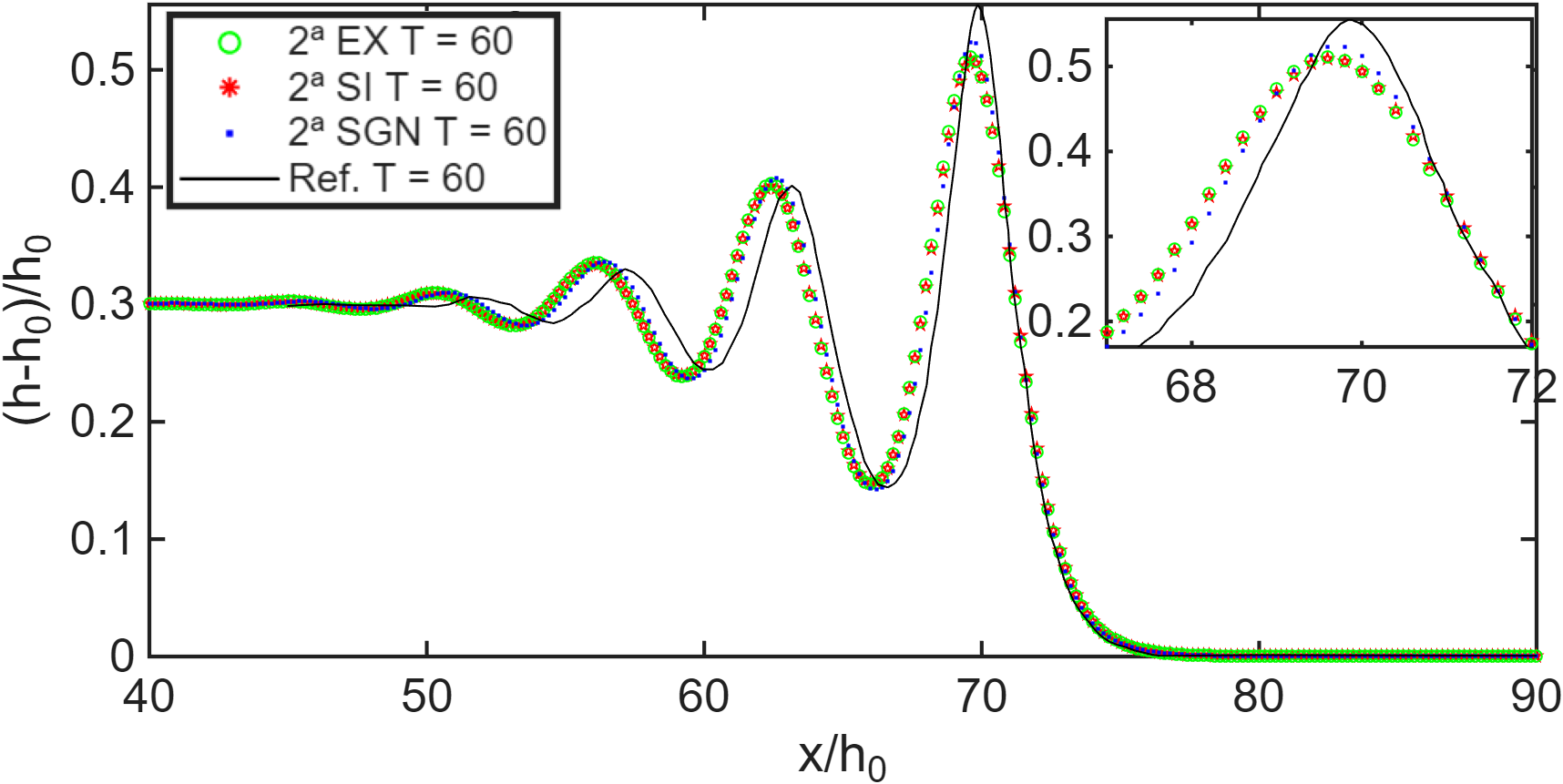}
         \caption{$T = 60,$ $\epsilon = 1.3$ and $\lambda = 100.$}
         \label{sfig:Favre_3_100_2}
     \end{subfigure} 
     \hfill
     \begin{subfigure}[b]{0.473\textwidth}
         \centering         \includegraphics[width=\textwidth]{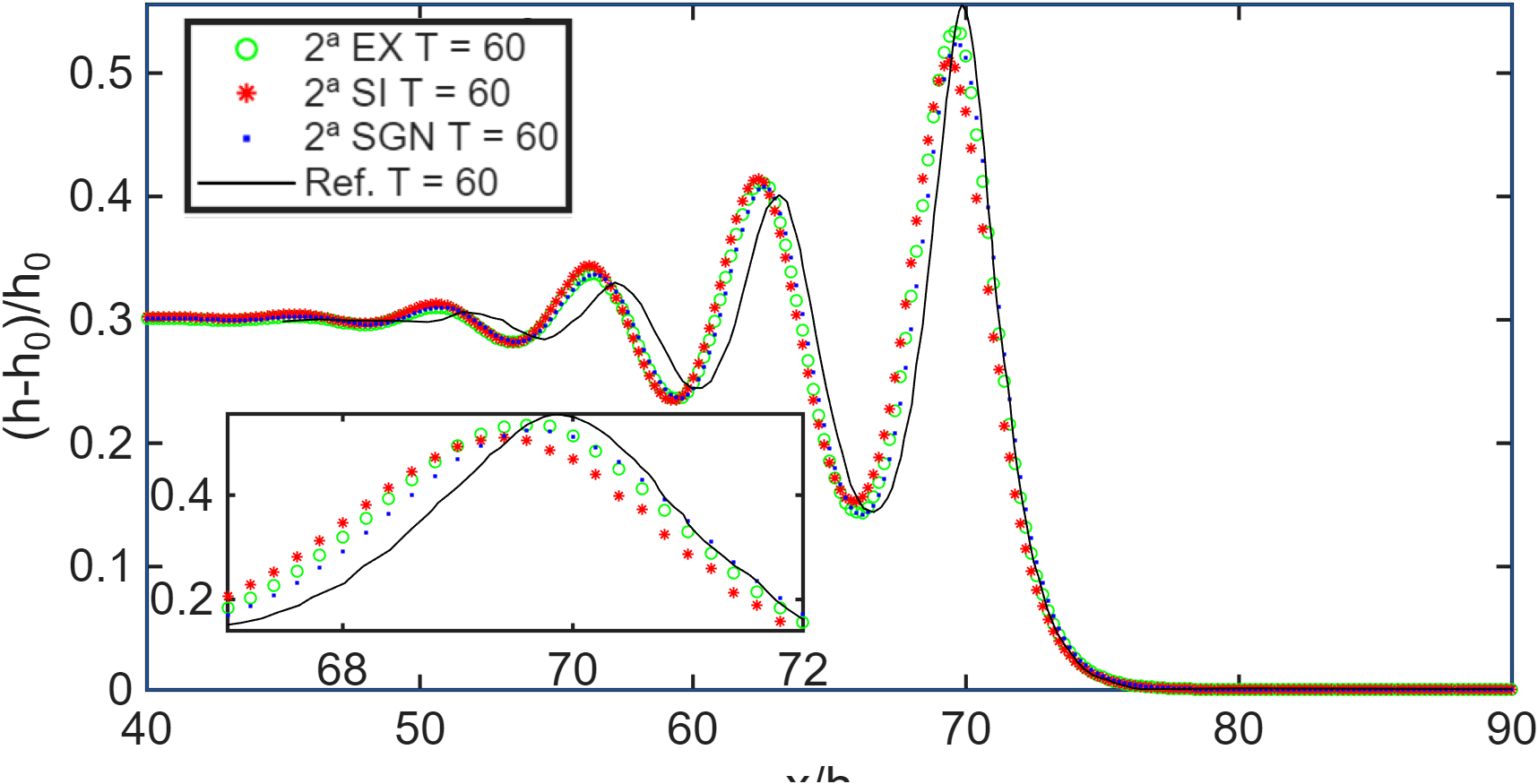}
         \caption{$T = 60,$ $\epsilon = 1.3$ and $\lambda = 500.$}
         \label{sfig:Favre_3_500_2}
     \end{subfigure} 
     \vfill\begin{subfigure}[b]{0.473\textwidth}
         \centering         \includegraphics[width=\textwidth]{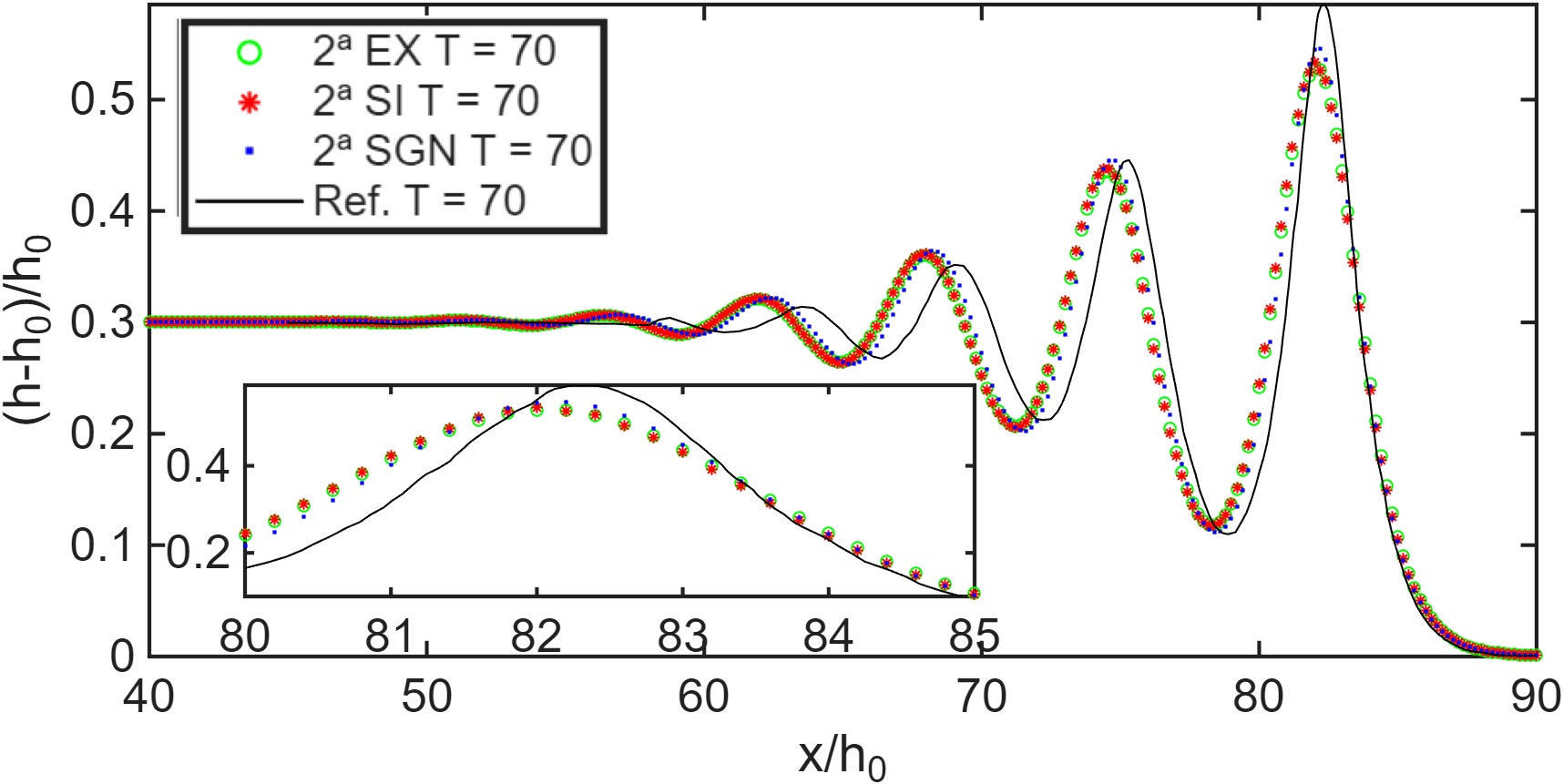}
         \caption{$T = 70,$ $\epsilon = 1.3$ and $\lambda = 100.$}
         \label{sfig:Favre_3_100_3}
     \end{subfigure} 
     \hfill
     \begin{subfigure}[b]{0.473\textwidth}
         \centering         \includegraphics[width=\textwidth]{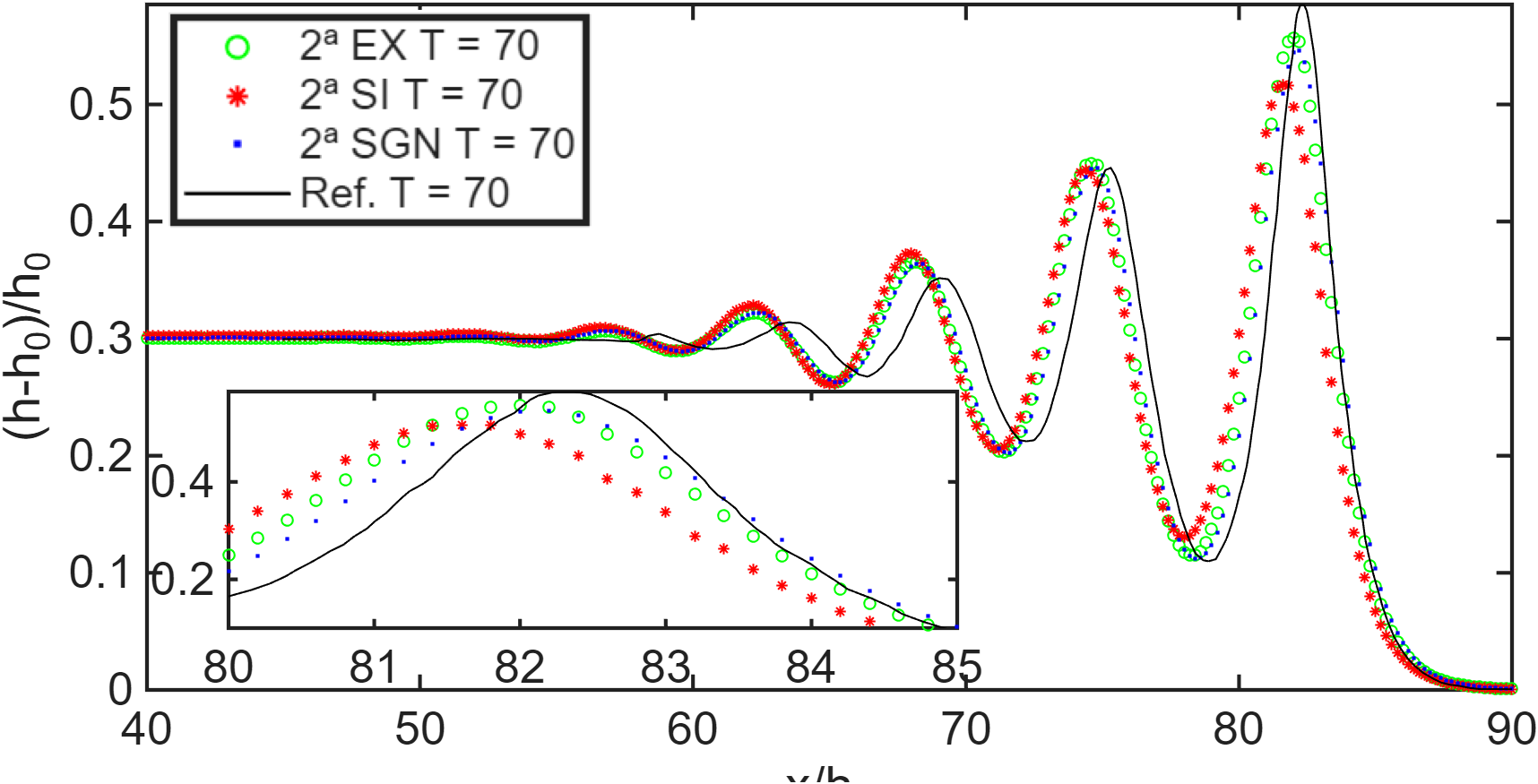}
         \caption{$T = 70,$ $\epsilon = 1.3$ and $\lambda = 500.$}
         \label{sfig:Favre_3_500_3}
     \end{subfigure} 
     \caption{Numerical solutions of Favre waves using second-order semi-implicit and explicit schemes, and second-order scheme for the Serre-Green-Nagdhi model \eqref{eq:SGN_standard}, at  different final times with nonlinearity coefficient $\epsilon=1.1$ and $\lambda \in\{ 100,500\}$. The solutions are plotted with the dimensionless variables $x/h_0$ on the horizontal axis and $(h-h_0)/h_0$ on the vertical axis. Reference: nonlinear potential flow  solution by \cite{wei1995fully}. The CFL numbers are $\mathrm{CFL}_{\rm IMEX}=4$, $\mathrm{CFL}_{\rm EX}=0.4$ and $\mathrm{CFL}_{\rm SGN}=0.9$.}
     \label{fig:favre_eps_3}
\end{figure}

The upstream depth is 
\[
h_1 = \varepsilon\,h_0,
\]
and the corresponding jump velocity is chosen to satisfy the Rankine–Hugoniot condition,
\[
u_1 \;=\;\sqrt{\frac{g\,(h_1+h_0)}{2\,h_0\,h_1}}\,(h_1-h_0).
\]
We introduce a smooth transition of width $\alpha=5\,h_0$ by means of the auxiliary functions
\[
g(x)=1-\tanh\!\bigl(\tfrac{x}{\alpha}\bigr), 
\]
so that the initial water depth and velocity profiles read
\[
h(x)=h_0 + \tfrac12\,(h_1-h_0)\,g(x),
\qquad
u(x)=\tfrac12\,u_1\,g(x),
\qquad
u_x(x)=\tfrac12\,u_1\,g'(x).
\]
The discontinuities in $h$ and $u$ across the jump are smoothed by the hyperbolic tangent profile.

Figures~\ref{fig:favre_eps_1}-\ref{fig:favre_eps_3} present the numerical solutions for the three values of $\epsilon$. In each subplot, the second-order explicit and semi-implicit solutions are shown together with the numerical solution of the Serre-Green-Naghdi model \eqref{eq:SGN_standard} at different final time $t$ where $t(T) = T\sqrt{gh_0}$ with $T = [50,60,70]$ for different values of $\lambda \in\{100,500\}$ adopting $N = 5000$ uniform mesh and the CFL numbers are $\mathrm{CFL}_{\rm IMEX}=4$, $\mathrm{CFL}_{\rm EX}=0.4$ and $\mathrm{CFL}_{\rm SGN}=0.9$. Both schemes, in all configuration, consistently exhibit comparable solution with the standard Serre-Green-Naghdi problem, allowing it to accurately maintain the plateau at $\epsilon-1$, which is a distinctive feature of Favre waves.
The results compare well with the reference solution by \cite{wei1995fully} obtained
with a fully nonlinear potential flow solver.

The effect of the nonlinearity coefficient $\epsilon$ is evident in the results. As $\epsilon$ increases, the waves exhibit more pronounced steepening, leading to sharper gradients and a more defined plateau. The second-order schemes consistently captures these characteristics with high accuracy, highlighting its robustness and ability to resolve and preserve intricate wave structures.

\subsubsection{CPU time benchmarking}
\begin{table}[!ht]
    \centering
    \begin{tabular}{|l|cc|cc|c|}
        \hline
        \multicolumn{6}{|c|}{\textbf{CPU time for Favre waves}} \\
        \hline
          \multicolumn{1}{|c}{\textbf{}} & \multicolumn{2}{c}{\textbf{$\lambda = 500$}} & \multicolumn{2}{c}{\textbf{$\lambda = 100$}} &  \\
        \hline
        $\epsilon$ & SI & EX & SI & EX & SGN \\
        \hline
        \multirow{3}{*}{$0.1$}
            & 7,09  & 173,28 & 0,63 & 10,65 & 2,16 \\
            & 8,50  & 207,51 & 0,76 & 13,00 & 2,45 \\
            & 9,81  & 241,80 & 0,85 & 15,51 & 3,13 \\
        \hline
        \multirow{3}{*}{$0.2$}
            & 7,70  & 180,10 & 1,03 & 10,23 & 2,26 \\
            & 7,83  & 198,22 & 1,34 & 14,79 & 3,27 \\
            & 8,88  & 232,00 & 1,53 & 17,53 & 4,07 \\
        \hline
        \multirow{3}{*}{$0.3$}
            & 7,61  & 162,67 & 2,50 & 10,19 & 2,60 \\
            & 8,54  & 210,23 & 3,03 & 12,68 & 3,33 \\
            & 10,26 & 267,89 & 3,79 & 15,68 & 3,98 \\
        \hline
    \end{tabular}
    \caption{CPU time (in second) for different values of $\lambda$ and $\epsilon$ for the semi-implicit and explicit second-order schemes applied to hyperbolic SGN and explicit second-order scheme applied to original SGN model. The CFL numbers are $\mathrm{CFL}_{\rm IMEX}=4$, $\mathrm{CFL}_{\rm EX}=0.4$ and $\mathrm{CFL}_{\rm SGN}=0.9$. }
    \label{tab:cpu_favre}
\end{table}
In these configurations, we conduct a detailed assessment of the computational cost associated with each scheme and benchmark these costs against those of the canonical SGN model \eqref{eq:SGN_standard}. Specifically, Table~\ref{tab:cpu_favre} and Figures~\ref{fig:favre_eps_1}-\ref{fig:favre_eps_3} illustrate that, for sufficiently small values of the parameter~$\lambda$, the semi-implicit scheme formulated for the hyperbolic SGN model achieves solution accuracy on par with both the original SGN discretization and a high-resolution reference solution. Remarkably, despite this comparable level of fidelity, the semi-implicit approach delivers a substantial speed-up in terms of CPU time. These results underscore not only the robustness of the semi-implicit scheme in capturing the essential dynamics of the standard SGN equations, but also its enhanced computational efficiency, rendering it an attractive and practical alternative for large-scale or real-time applications.

\subsubsection{Comparison with experimental data}\label{sec:long_time}
To validate the performance of the semi‐implicit hyperbolic relaxation scheme against the fully explicit one 
we now consider 
 extended propagation distances until the dimensionless final time 
\[
\tilde t \;=\; t\sqrt{\frac{g}{h_0}} \;=\; 256,
\]
representative of the  experiments by \cite{Treske}.

\begin{figure}[!ht]
  \centering
  \includegraphics[width=0.9\textwidth]{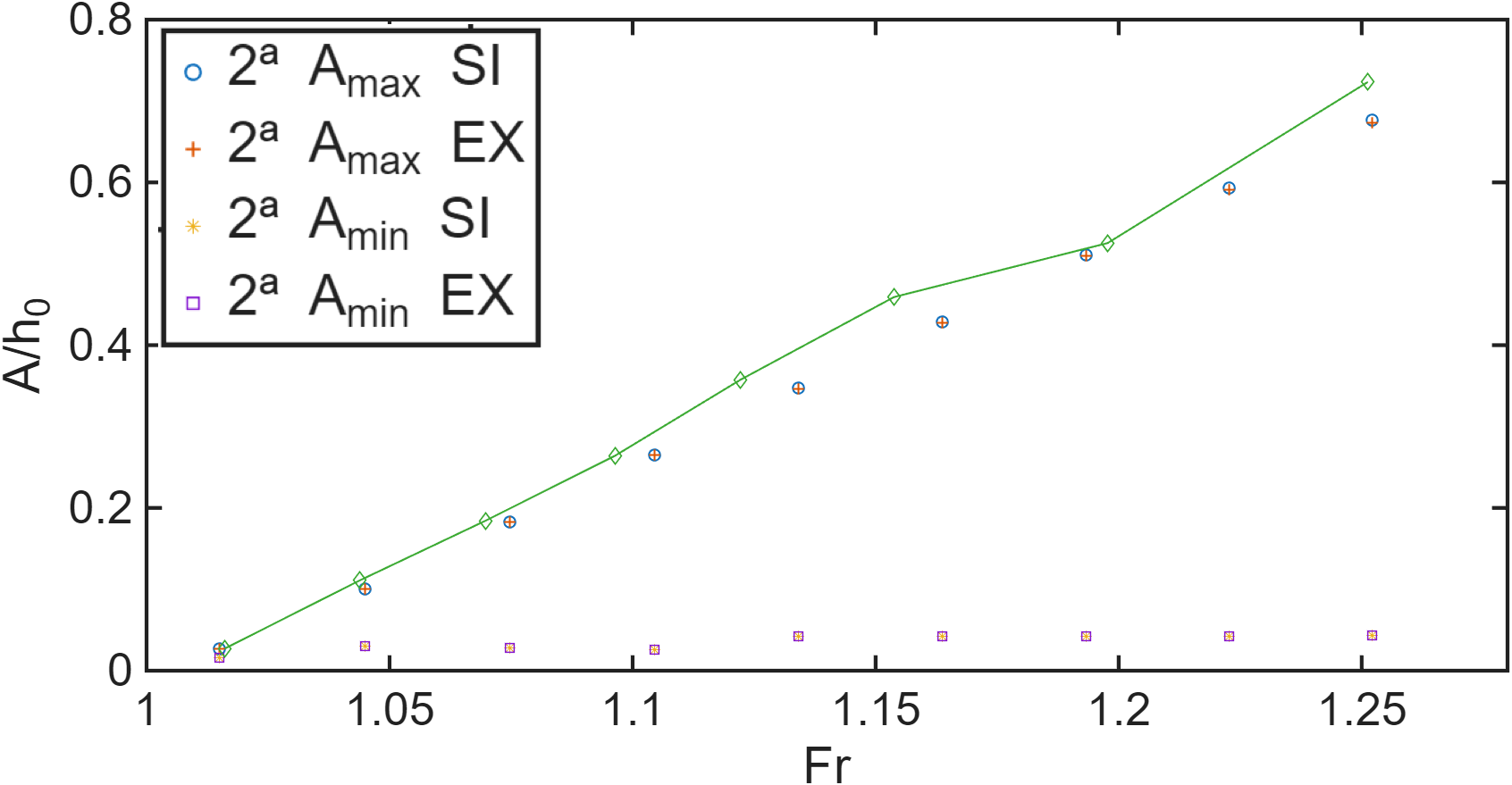}
  \caption{Favre waves: heights of the first peak  (circles) and trough (squares)  at \(\tilde t=256\) as a function of the Froude number.   Solid line with circles:   data  by \cite{Treske}.}
  \label{fig:long_time}
\end{figure}

The computations are now run on  the spatial domain \(x\in[0,73.58]\), initial jump at \(x_0=68.58\), and hyperbolic relaxation parameter \(\lambda=500\). 

Figure~\ref{fig:long_time} presents the free‐surface profile at \(\tilde t=256\), overlaying numerical solutions with the laboratory measurements   by \cite{Treske}. The numerical solutions of both schemes agree closely with the experimental data, with maximum elevation errors below $2\%$ and phase deviations within one grid cell. However, the semi‐implicit relaxation method allows a tenfold increase in the stable time step compared to the explicit solver, leading to a reduction in total CPU time of approximately $70\%$ for this test.  This substantial gain in computational efficiency, combined with the high-fidelity agreement to laboratory measurements, underscores the advantage of the semi‐implicit scheme for long‐time simulations of weakly dispersive free‐surface flows.

\subsubsection{Wave propagation over a shelf}
In this section we examine the evolution of a solitary wave of amplitude $A = 0.2\mathrm{m}$ as it traverses a submerged shelf.  The shelf causes the still‑water depth to vary from the upstream depth $h_{0} = 1\mathrm{m}$ to a shallower downstream depth of $0.5\mathrm{m}$ via a ramp of slope $1:20$. The computational domain extends from $x=0$ to $x=280,\mathrm{m}$, and initially at time $t=0$ the wave crest is positioned at $x_{0}=80\mathrm{m}$.
\begin{figure}[!ht]
\centering
\includegraphics[width=0.95\textwidth]{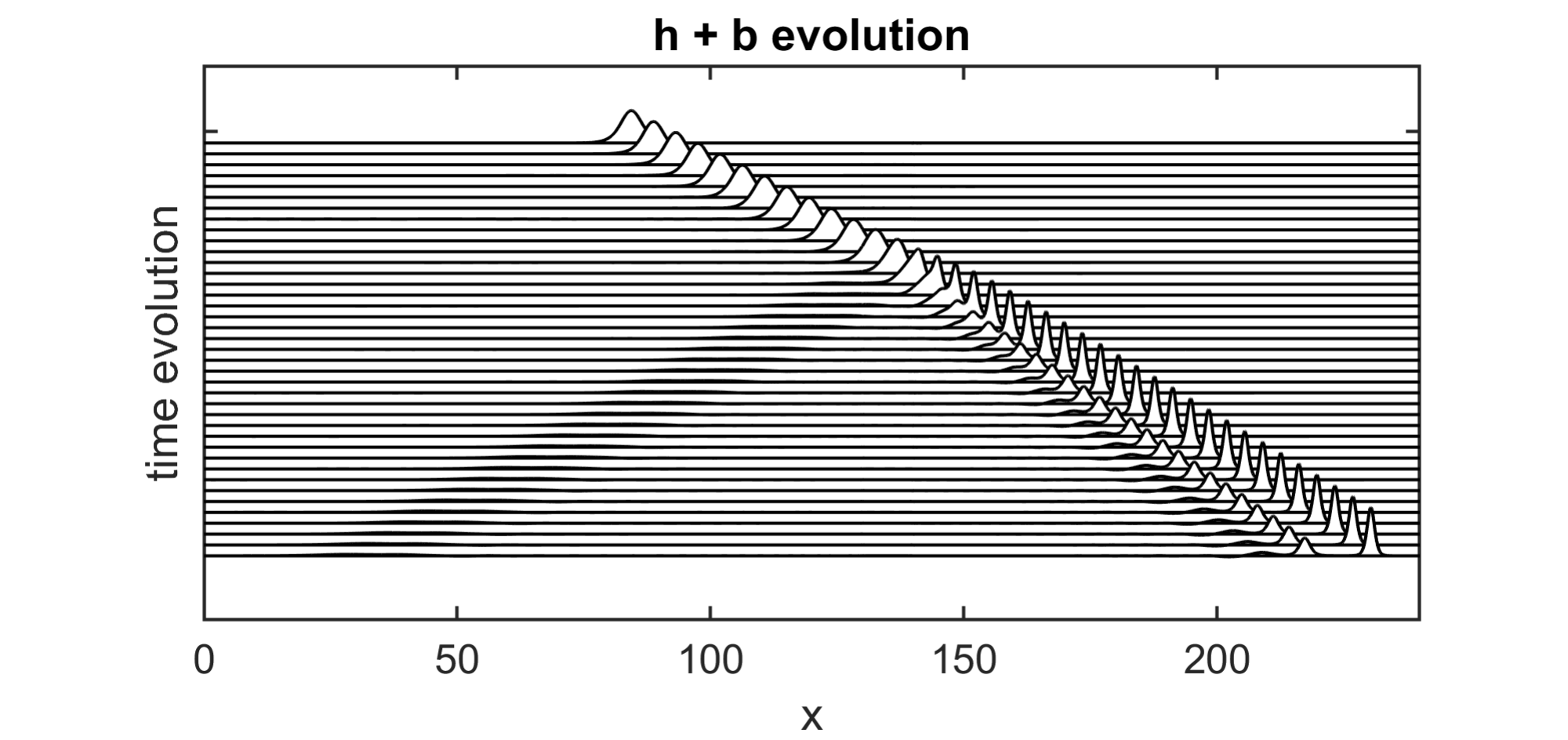}
\caption{Evolution of the solitary wave as it propagates over the submerged shelf, showing reflection and splitting into three pulses for the explicit scheme.}
\label{fig:wave_split_ex}
\end{figure}
\begin{figure}[!ht]
\centering
\includegraphics[width=0.95\textwidth]{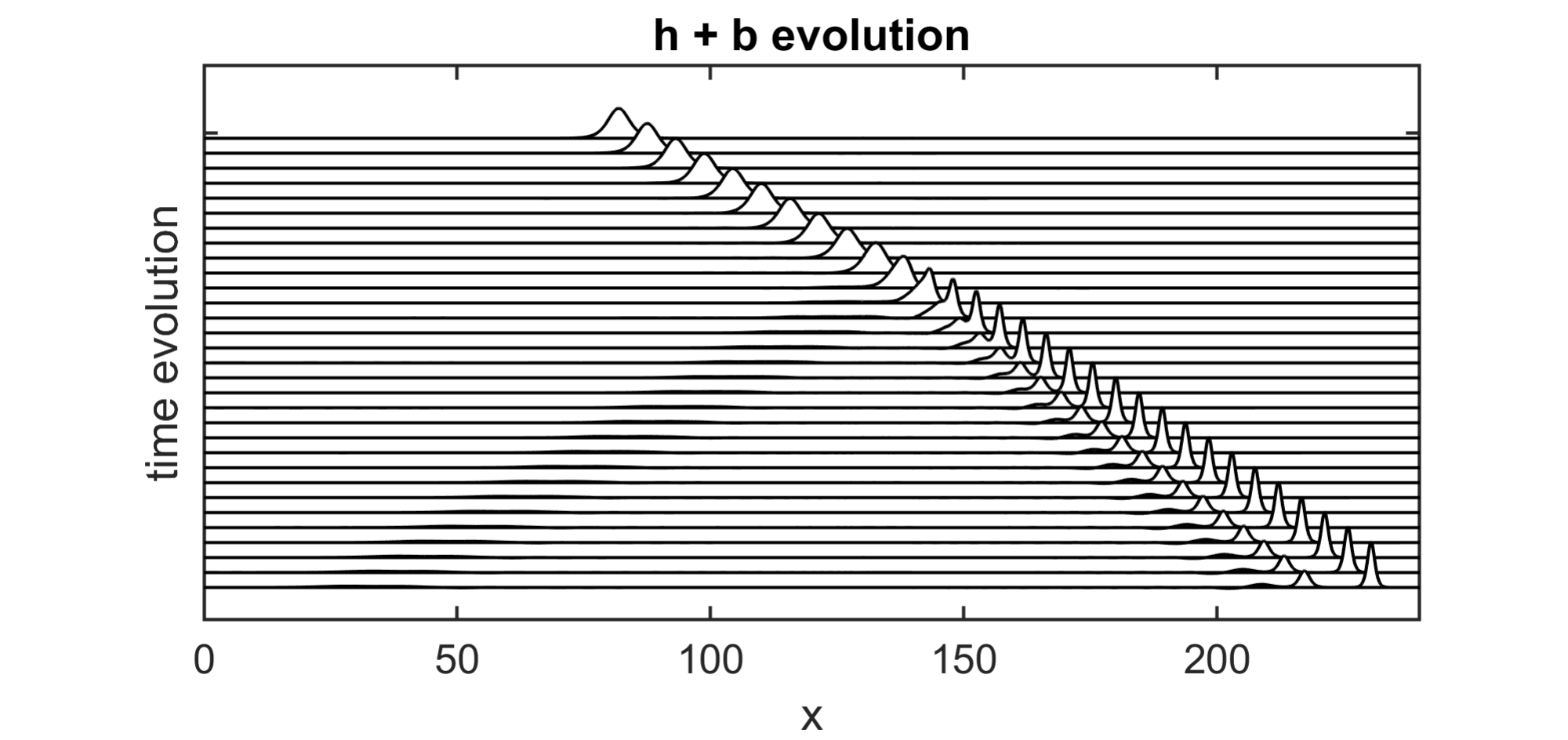}
\caption{Evolution of the solitary wave as it propagates over the submerged shelf, showing reflection and splitting into three pulses for the semi-implicit scheme.}
\label{fig:wave_split_im}
\end{figure}
\begin{figure}[!ht]
\centering
\includegraphics[width=0.95\textwidth]{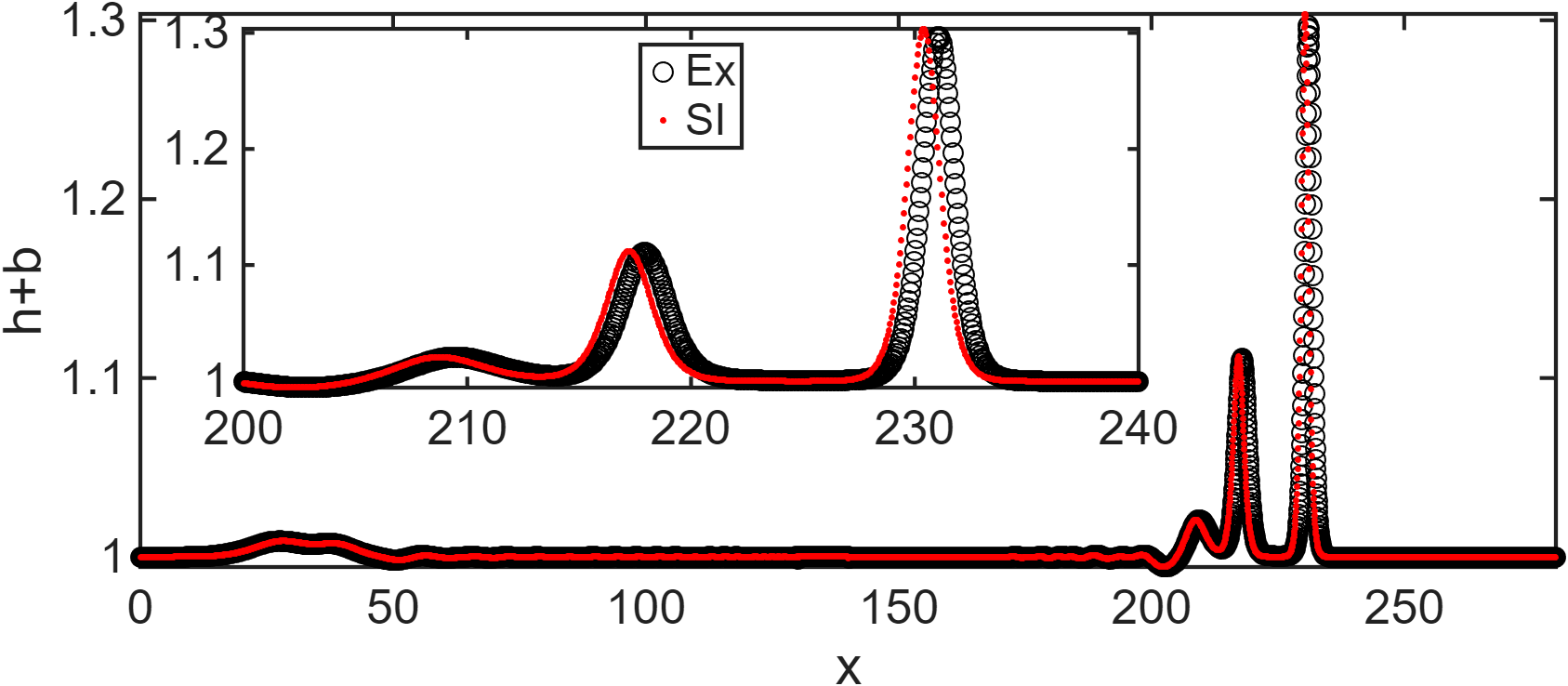}
\caption{Comparison of the solitary wave at the final time as it propagates over the submerged shelf, showing reflection and splitting into three pulses for the semi-implicit and explicit schemes.}
\label{fig:wave_comparison}
\end{figure}

The bathymetric profile is described by a piecewise linear function:
\begin{align}
b(x) &= \begin{cases}
0, & x \le 130,\\
\frac{x - 130}{20}, & 130 < x < 140,\\
0.5, & x \ge 140,
\end{cases}
\end{align}
Accordingly, the still‑water depth is given by $h(x,0) = h_{0} - b(x)$.

Upstream of the shelf ($x \le 130,\mathrm{m}$), we impose an exact solitary wave solution of the enhanced Boussinesq system:
\begin{align}
h(x,0) &= h_{0}\Bigl[1 + \varepsilon\sech^{2}\bigl(k(x - x_0)\bigr)\Bigr],\\
u(x,0) &= c\Bigl(1 - \frac{h_{0}}{h(x,0)}\Bigr).
\end{align}

We discretize the domain using a uniform grid with $N = 5000$ and advance in time with CFL$= 3$ and $\lambda = 500.$

As the incident soliton encounters the shoal, a small reflected wave is generated upstream, and the primary pulse splits into three transmitted pulses downstream, each of increasing amplitude.  Figures~\ref{fig:wave_split_ex},\ref{fig:wave_split_im} and \ref{fig:wave_comparison} illustrate the free‑surface elevation after interaction, clearly showing the separated pulses for the explicit and semi-implicit scheme showing similar qualitative results.

\section{Dingemann's experiment}\label{sec:Dingemann}

The Dingemans flume experiment was originally conceived to investigate nonlinear wave shoaling and transformation over a gently sloping submerged shoal under controlled laboratory conditions.  In the physical setup, a monochromatic wavetrain of small amplitude $A = 0.02\,$m propagates in a flume of constant still‐water depth $h_{0} = 0.8\,$m.  The shoal consists of a smooth upslope rising to $0.6\,$m above the undisturbed bottom, followed by a short plateau and a symmetrical downslope returning to the original depth.  Four capacitance gauges located at 
\[
x = 3.04,\;9.44,\;20.04,\;\text{and}\;26.04\;\mathrm{m}
\]
\begin{figure}[!ht]
  \centering
  \begin{subfigure}[t]{0.48\textwidth}
    \includegraphics[width = \textwidth]{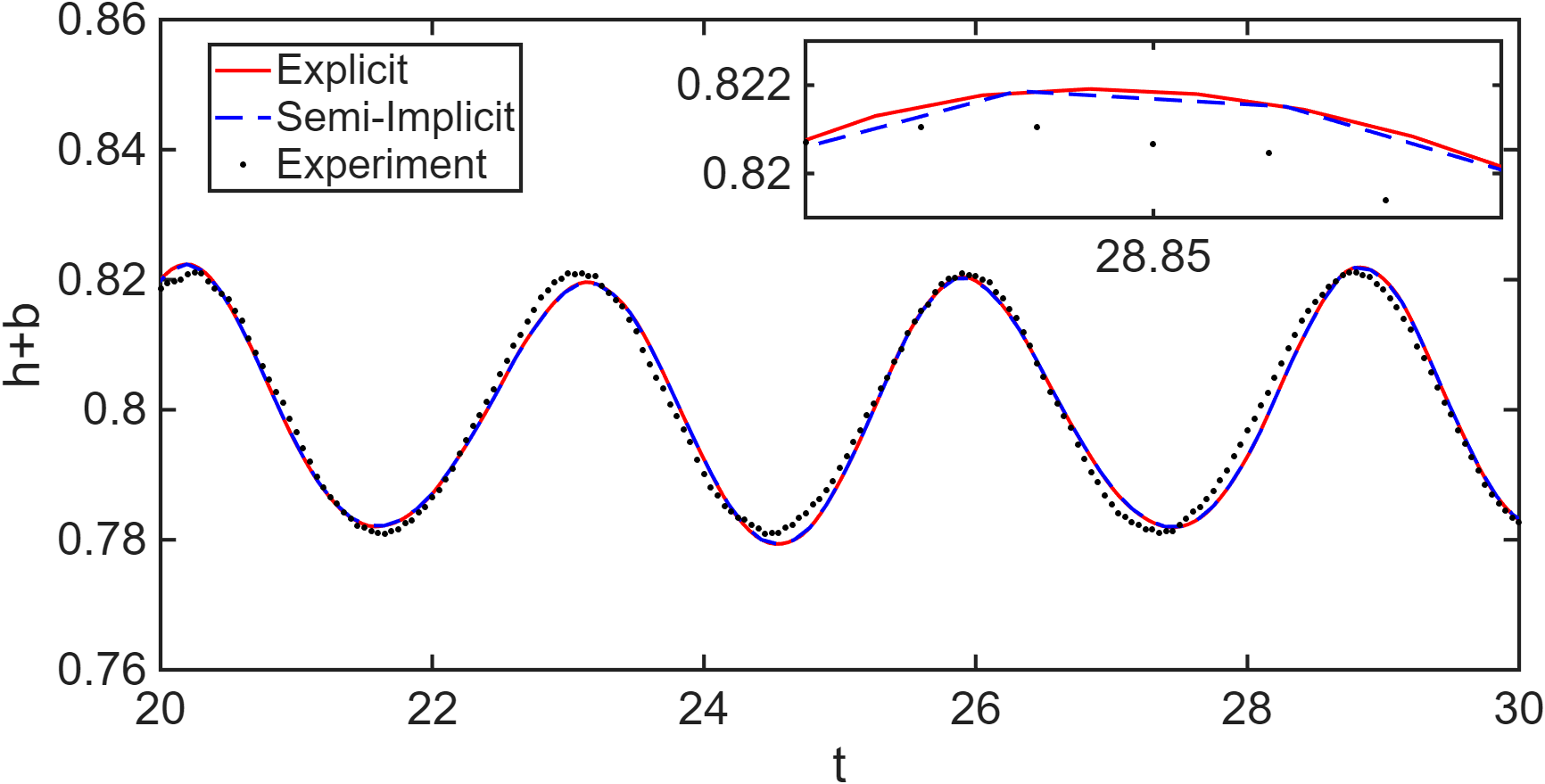}
    \caption{$x=3.04\,$m}
  \end{subfigure}% 
  \hfill
  \begin{subfigure}[t]{0.48\textwidth}
    \includegraphics[width = \textwidth]{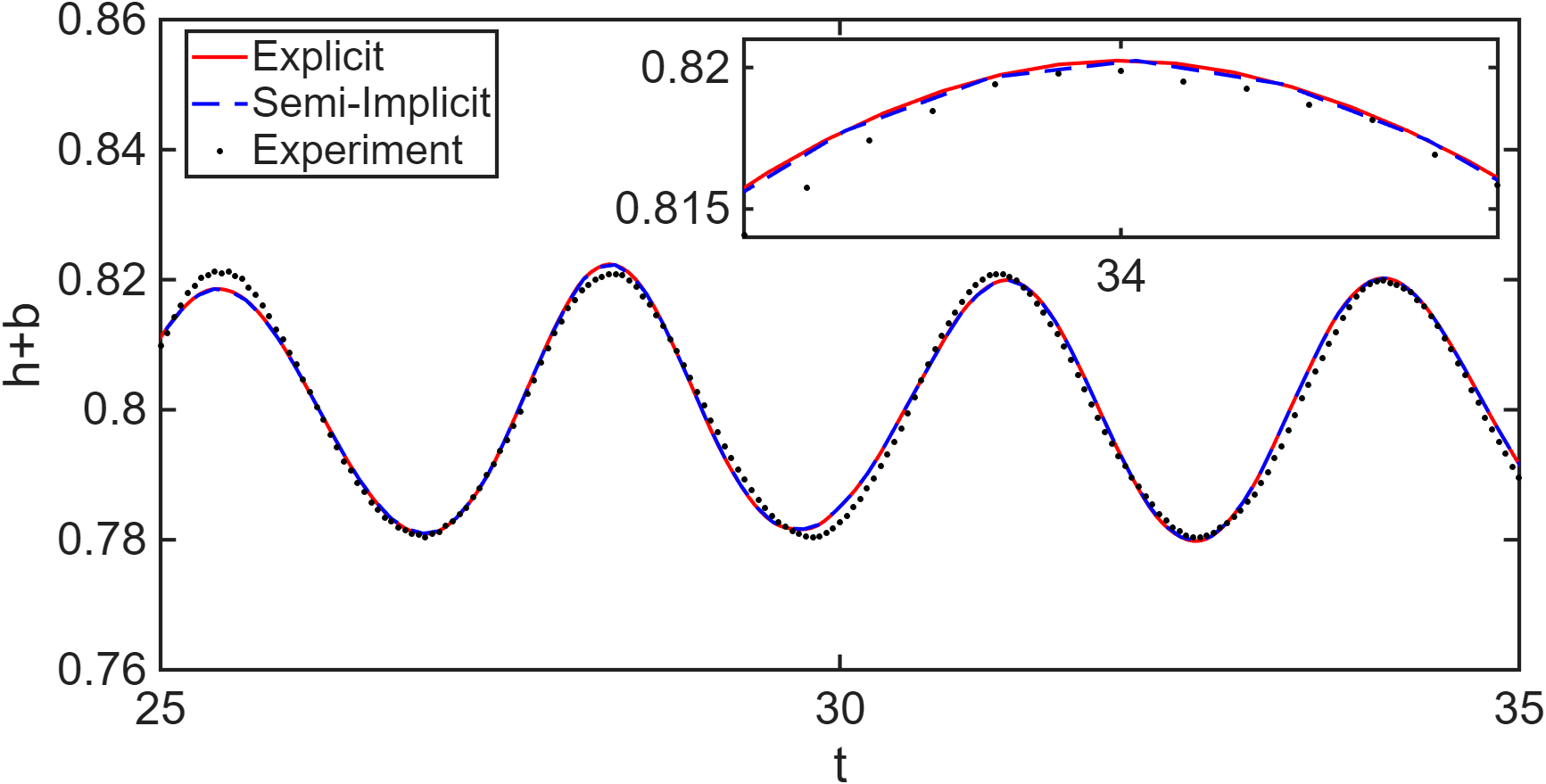}
    \caption{$x=9.44\,$m}
  \end{subfigure} \\
  \begin{subfigure}[t]{0.48\textwidth}
    \includegraphics[width = \textwidth]{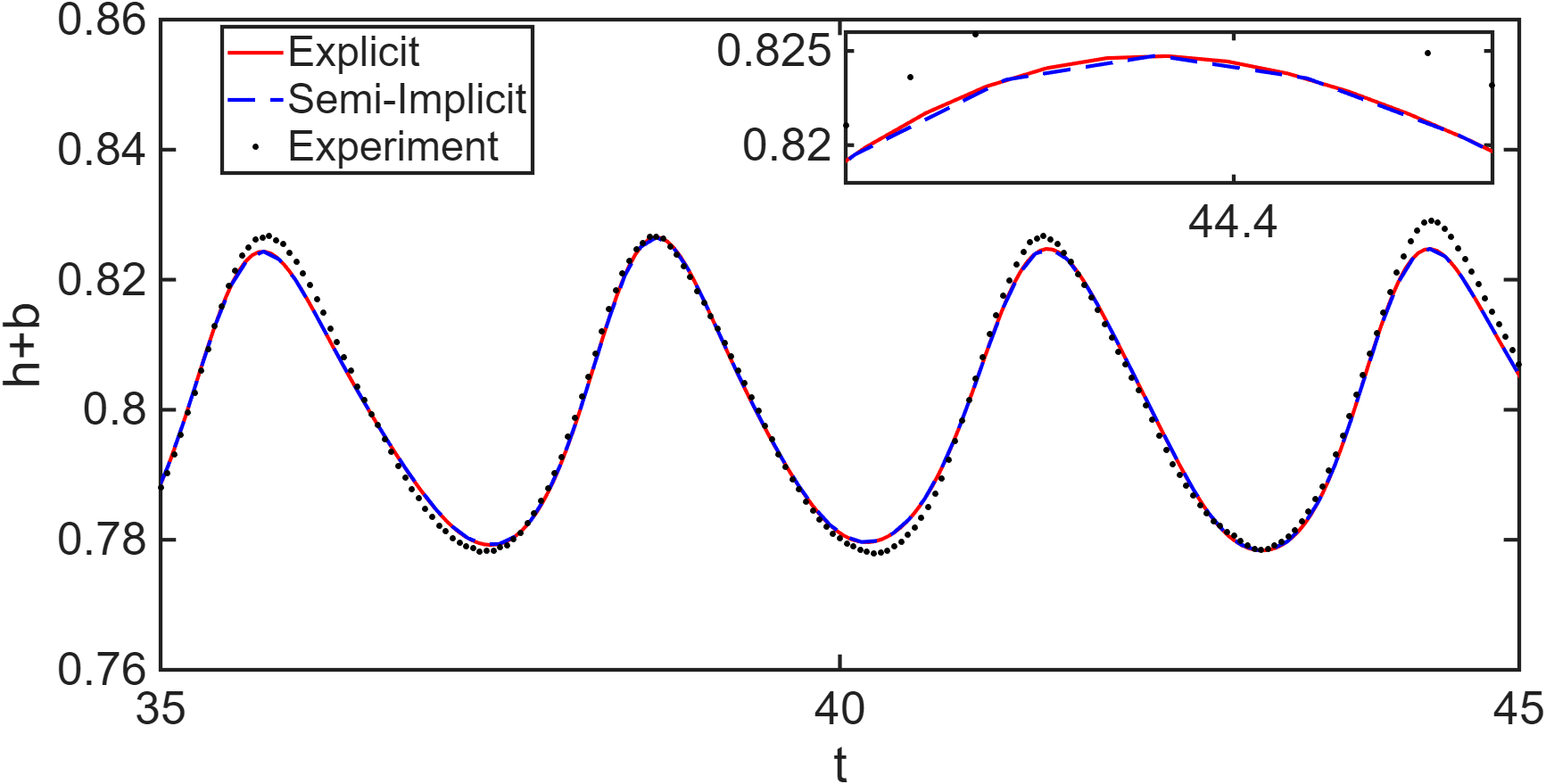}
    \caption{$x=20.04\,$m}
  \end{subfigure}
  \begin{subfigure}[t]{0.48\textwidth}
    \includegraphics[width = \textwidth]{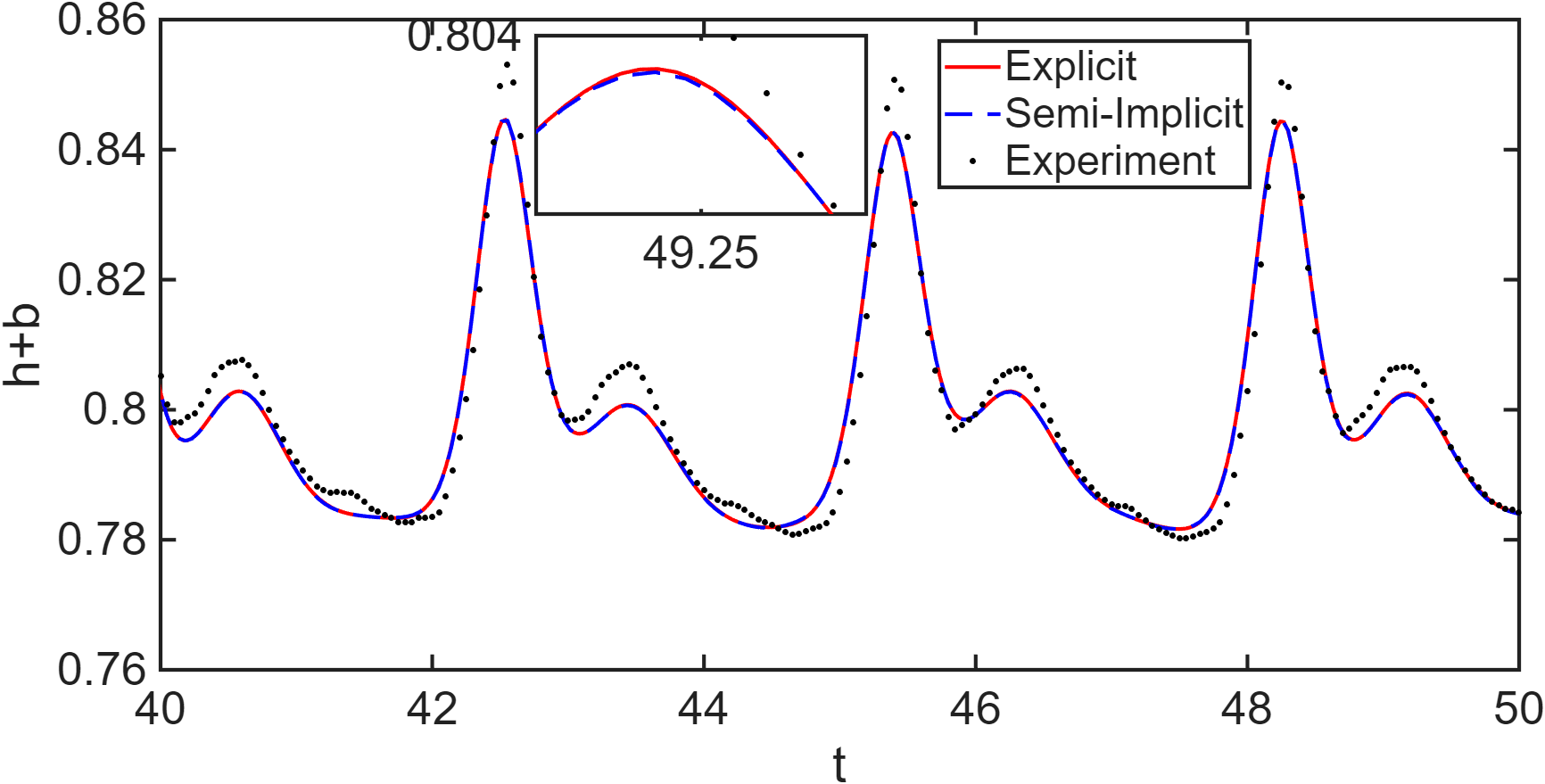}
    \caption{$x=26.04\,$m}
  \end{subfigure}
  % \\
  % \begin{subfigure}[t]{0.48\textwidth}
  %   \includegraphics[width = \textwidth]{Figures/Dingemann/x30.png}
  %   \caption{$x=30.44\,$m}
  % \end{subfigure}%
  % \begin{subfigure}[t]{0.48\textwidth}
  %   \includegraphics[width = \textwidth]{Figures/Dingemann/x36.png}
  %   \caption{$x=37.04\,$m}
  % \end{subfigure}
  \caption{Time‐series comparison at gauge locations.  Semi‐implicit SGN (solid), explicit SGN (dashed), experimental data (symbols).}
  \label{fig:dingemans_results}
\end{figure}
\noindent record the time series of free‐surface elevation over a duration of $70\,$s, providing high‐resolution experimental data for model validation.

To numerically reproduce these measurements, we solve the one‐dimensional hyperbolic SGN equations on the domain $x\in[-140,\,100]\,$m.  The initial condition imposes a sinusoidal surface displacement
\[
z(x,0) = 
\begin{cases}
h_{0} + A\cos(kx), & -\tfrac{34.5\pi}{k} \le x \le -\tfrac{8.5\pi}{k},\\
h_{0}, & \text{otherwise},
\end{cases}
\]
with $k = 0.8406220896381442\,$rad/m determined by solving $g\,k\,\tanh(kh_{0}) = \omega^{2}$ for the given wave frequency and $h_0 = 0.8$m.  The corresponding velocity is initialized as
\[
u(x,0) \;=\; \sqrt{\frac{g}{k}\,\tanh(kh_{0})}\,\frac{z(x,0)-h_{0}}{h_{0}}\,.
\]
The shoal bathymetry $b(x)$ is defined piecewise to match the laboratory profile, and the initial water depth is $h(x,0)=z(x,0)-b(x)$.

Spatial discretization uses a uniform mesh with $N = 15000$ points, and time integration proceeds until $t_{\mathrm{final}}=70\,$s under a $\mathrm{CFL}=3$, with $\lambda = 100$.

Figure~\ref{fig:dingemans_results} presents the free‐surface elevation time series at all six gauges.  Both numerical solutions closely follow the experimental measurements in phase and amplitude.  No systematic bias is observed between the two schemes, confirming that the hyperbolic SGN model faithfully captures the wave dynamics even with the moderate relaxation parameter.

Remarkably, despite advancing with time steps three times larger, the semi‐implicit scheme retains the same high fidelity as the explicit method when reproducing the experimental wave trains.

Figure~\ref{fig:dingemans_results} presents the time series of the free-surface elevation recorded data. Both numerical solutions exhibit an excellent agreement with the experimental measurements, accurately reproducing both the amplitude and phase of the wave trains. No systematic bias can be identified between the two schemes, confirming that the hyperbolic SGN model is able to faithfully capture the main features of the wave dynamics, even when using a moderate value of the relaxation parameter.

Remarkably, even when advancing with time steps three times larger, the semi-implicit scheme maintains the same level of accuracy as the explicit method in reproducing the experimental signals. This result highlights the enhanced stability of the semi-implicit formulation, which allows for a substantial relaxation of the time-step constraint without any loss of physical fidelity. In addition, the computational cost remains comparable, with total CPU times of $CPU_{\mathrm{IMEX}} = 5.1806 \times 10^{3}$ and $CPU_{\mathrm{EX}} = 5.7232 \times 10^{3}$. 
From a performance perspective, the semi-implicit solver completes the simulation in slightly less time than the explicit scheme. This gain arises from the enlarged time‐step $\Delta t$, which reduces the total number of updates by a factor of three while maintaining identical accuracy.  Such a computational advantage is critical for large‐scale or real‐time wave‐prediction applications, where both precision and efficiency are essential.

\section{Conclusions}
In this paper we have designed a semi-implicit finite volume scheme for the solution the Serre-Green-Naghdi (SGN) equations written in first-order hyperbolic form. The stiff part of the system, which includes relaxation terms and pressure fluxes, has undergone an implicit discretization, while the convection sub-system and the remaining terms have been discretized explicitly. The novel scheme has been compared in terms of both accuracy and efficiency against the classical mixed elliptic-hyperbolic solver for the SGN model as well as against a fully explicit method for the hyperbolic SGN system. Several benchmarks have been performed to verify accuracy and robustness, and performing simulations of solitary wave propagation and dispersive wave breaking phenomena. This ultimately allows us to provide some pragmatic considerations for the effective use of semi-implicit schemes in the context of dispersive wave propagation models. 

Further developments will consider the extension of the proposed methodology to multidimensional unstructured meshes in order to capture realistic landscape geometries and the design of high order schemes in space and time to increase even further the accuracy and to reduce the error related to the propagation of dispersive waves.

\section*{Acknowledgements}
This research has received funding from the European Union’s NextGenerationUE – Project: Centro Nazionale HPC, Big Data e Quantum Computing, “Spoke 1” (No. CUP E63C22001000006). E. Macca was partially supported by: GNCS No. CUP E53C24001950001 Research Project "Soluzioni Innovative per Sistemi Complessi: Metodi Numerici e Approcci Multiscala"; PRIN 2022 PNRR “FIN4GEO: Forward and Inverse Numerical Modeling of hydrothermal systems in volcanic regions with application to geothermal energy exploitation”, No. P2022BNB97; PRIN 2022 “Efficient numerical schemes and optimal control methods for time-dependent partial differential equations”, No. 2022N9BM3N - Finanziato dall’Unione europea - Next Generation EU – CUP: E53D23005830006. E. Macca would like to thank the Italian Ministry of Instruction, University and Research (MIUR) to support this research with funds coming from PRIN Project 2022  (2022KA3JBA, entitled “Advanced numerical methods for time dependent parametric partial differential equations and applications”). E. Macca and W. Boscheri are members of the INdAM Research group GNCS. %WB received financial support from the Italian Ministry of University and Research (MUR) with the PRIN Project 2022 No. 2022N9BM3N. 
MR is a member of the CARDAMOM research team,  Inria at University of Bordeaux research center.  

\section*{Conflict of interest}
The authors declare that they have no conflict of interest.

\bibliographystyle{plain}
\bibliography{biblio}
\end{document}